\documentclass[article,11pt]{aiaa-pretty}    % journal, submit, article
\usepackage{amssymb}
\usepackage{amsmath}
\usepackage{booktabs}
 \usepackage{varioref}%  smart page, figure, table, and equation referencing
 \usepackage{wrapfig}%   wrap figures/tables in text (i.e., Di Vinci style)
 \usepackage{threeparttable}% tables with footnotes
 \usepackage{dcolumn}%   decimal-aligned tabular math columns
  \newcolumntype{d}{D{.}{.}{-1}}
 \usepackage{nomencl}%   nomenclature generation via makeindex
  \makeglossary
 \usepackage{subfigure}% subcaptions for subfigures
 \usepackage{graphicx}
 \usepackage{subfigmat}% matrices of similar subfigures, aka small mulitples
 \usepackage{fancyvrb}%  extended verbatim environments
  \fvset{fontsize=\footnotesize,xleftmargin=2em}
 \usepackage{lettrine}%  dropped capital letter at beginning of paragraph
\usepackage{cite}
\usepackage{enumitem}

\newtheorem{remark}{Remark}

\newcommand{\set}[1]{\left\{#1\right\}}
\newcommand{\real}[1]{{\mathbb R}^{#1}}
\newcommand{\bb}{{\boldsymbol b}}

\newcommand{\be}{{\boldsymbol e}}
\newcommand{\bff}{{\boldsymbol f}}

\newcommand{\bu}{{\boldsymbol u}} 

\newcommand{\bv}{{\boldsymbol v}}
\newcommand{\bw}{{\boldsymbol w}}
\newcommand{\bx}{\boldsymbol x}

\newcommand{\bcheb}{{\bold{cheb}}}

\newcommand{\bxf}{{\bx(\cdot)}}  % f for function
\newcommand{\buf}{{\bu(\cdot)}}  % f for function
\newcommand{\bB}{{\boldsymbol B}}

\newcommand{\bzero}{{\bf 0}}

\newcommand{\bX}{{\mbox{\boldmath $X$}}}
\newcommand{\bV}{{\mbox{\boldmath $V$}}}
\newcommand{\bU}{{\mbox{\boldmath $U$}}}

\newcommand{\bnu}{{\mbox{\boldmath $\nu$}}}
\newcommand{\bmu}{{\mbox{\boldmath $\mu$}}}

\newcommand{\blam}{{\mbox{\boldmath $\lambda$}}}

\newcommand{\bomega}{\mbox{\boldmath$\omega$}}

\newcommand{\bOmega}{\mbox{\boldmath$\Omega$}}
\newcommand{\bLam}{\mbox{\boldmath$\Lambda$}}
\newcommand{\bW}{{\mbox{\boldmath $W$}}}

\begin{document}
%=====================================================
\title{Universal Birkhoff Method for Computing Extremals in the Elliptic Restricted Three-Body Problem}
%=======================================================
\author{Michael J. Dixon\thanks{Corresponding author; Commander, USN, Department of Mechanical and Aerospace Engineering; michael.dixon@nps.edu. Currently, Assistant Professor, United States Naval Academy, Annapolis, MD.} \ and \
Isaac M. Ross\thanks{Distinguished Professor, Department of Mechanical and Aerospace Engineering.}
\\
\textit{Naval Postgraduate School, Monterey, CA 93943}
}

%=========================================================================
\abstract{
The computation of finite-thrust extremal arcs in the elliptic restricted three-body, trajectory optimization problem is considered. Libration point orbits are approximated to near-machine-precision using a fast Fourier transform of the sampled values of the state vector at Chebyshev-Gauss-Lobatto points. Checkable optimality conditions are derived by applying Pontryagin's Principle to minimum-time and time-constrained minimum-propellant problems.  These necessary conditions include criteria for optimal departure and arrival points.  For propellant consumption, a recently-developed computational model is employed.  This model is agnostic to the specific impulse of the propellant and varies as an inverse-quadratic of a cosine term.
Candidate optimal solutions are generated by combining the universal Birkhoff theory for trajectory optimization with the fast, guess-free spectral algorithm. The extremality of the Birkhoff-computed solution is validated against the Hamiltonian minimization condition and the transversality conditions. It is shown that the Birkhoff-theoretic spectral algorithm can generate verifiable extremals without any assistance of initialization from dynamical systems theory.
}		
%=========================================================================

\maketitle{}

\footnote{A version of this paper was presented at the 2025 AAS/AIAA Space Flight Mechanics Meeting, Kauai, HI, Jan. 19-23, 2025, Paper AAS-25-197.}

\section{Introduction}\label{sec:intro}
Optimal orbital maneuvers in the three-body problem are unlike their two-body counterparts because of the existence of a vast number of low-energy pathways in the former case \cite{gmez_study_1993, belbruno_capture_2004, koon_dynamical_2011, parker_lowenergy_2014}.  The existence of these low-energy trajectories has been discovered using the tools of dynamical systems theory\cite{koon_dynamical_2011,parker_lowenergy_2014, howell_application_1997}.
A conceptual mapping of these results to optimal control theory \cite{longuski_optimal_2014, clarke_nonsmooth_1998, bryson_applied_1978, ross_primer_2015} suggests that there is a limitless number of minimal Hamiltonian trajectories for low-energy transfers.  Here, the Hamiltonian refers to the Pontryagin Hamiltonian\cite{clarke_nonsmooth_1998, ross_primer_2015} and not the Hamiltonian describing the dynamics of the three-body problem \cite{belbruno_capture_2004,koon_dynamical_2011}.  In principle, the low-energy trajectories offer mission designers a vast number of choices of minimal impulse solutions.  Some engineering challenges in using this body of results are as follows:
\begin{enumerate}
\item Transforming a solution from an idealized case, such as the circular restricted three-body problem (CR3BP), to a high-fidelity model requires an optimization process that does not necessarily replicate the low-energy result \cite{park_leveraging_2022}.
\item The impulsive maneuvers assumed in the idealized cases must be re-optimized for finite-burn practical maneuvers, particularly for lower values of maximum thrust\cite{finite-burn-1966,finite-burn-2004}.  This re-optimization is typically performed using computational optimal control techniques \cite{mingotti_low-energy_2009,chupin_transfer_2018,pritchett_impulsive_2018,finite-burn-2004}.
\item It is not easy to produce minimum-time and time-bound solutions using dynamical systems theory\cite{senent_low-thrust_2005, bonnard_time-minimum_2016, caillau_minimum_2012}.
\item The process of utilizing a solution obtained from dynamical systems theory to a practical construction of an optimal trajectory is not trivial \cite{parker_lowenergy_2014,park_leveraging_2022,pritchett_impulsive_2018,yoon_minimum-fuel_2023, singh_exploiting_2021,beolchi_multiple-arc_2024}.
\end{enumerate}
In view of these engineering challenges, there has been an increasing interest\cite{ chupin_low-thrust_2017, chupin_transfer_2018, du_transfer_2022, du_low-thrust_2023, kelly_orthogonal_2023, singh_exploiting_2021}
in using computational optimal control theory to assist trajectory optimization performed via dynamical systems theory.
A frequently-used\cite{mingotti_optimal_2011,chupin_transfer_2018, kelly_orthogonal_2023} optimization criterion is propellant consumption.  A number of algorithms\cite{senent_low-thrust_2005, ozimek_low-thrust_2010, kelly_orthogonal_2023, chupin_low-thrust_2017, chupin_transfer_2018, perez-palau_fuel_2018} 
use the concept of indirect shooting\cite{ross_primer_2015,conway_survey_2012, caillau_algorithmic_2023} with assistance from dynamical systems theory.  Because an indirect shooting algorithm is plagued with the ``curse of sensitivity'' \cite{ross_primer_2015},  results from dynamical systems theory are extremely powerful in generating a starting point (i.e. a good initial guess of the solution) for the optimization algorithm. While collocation methods \cite{conway_survey_2012,caillau_algorithmic_2023} are purportedly more robust to initial guess sensitivity, hybrid approaches have emerged that combine elements of dynamical systems theory with collocation or indirect shooting. These have been demonstrated in several studies \cite{senent_low-thrust_2005,singh_exploiting_2021,mingotti_combined_2007,chupin_transfer_2018,ozimek_low-thrust_2010, parrish_low-thrust_2016,pritchett_low-thrust_2017, pritchett_impulsive_2018}.

The three-body problem is an established, mathematically-hard problem \cite{ams-1913, moser-1971}. One of the challenges in the restricted three-body trajectory optimization problem is the location of the optimal arrival and departure points\cite{senent_low-thrust_2005, mingotti_combined_2007, bonnard_time-minimum_2016}.  This is because the libration point orbits are not specified in terms of element sets as in a classical two-body problem; rather, they are determined numerically through various well-established methods\cite{ferrari_periodic_2018,koon_dynamical_2011,howell_application_1997}.  One approach to address the optimality of the arrival problem is to simply insert the spacecraft along a stable manifold\cite{senent_low-thrust_2005, ozimek_low-thrust_2010, kelly_orthogonal_2023}.  In principle, this idea solves the problem if:
\begin{enumerate}
\item[i)] the cost functional is limited to propellant consumption,
\item[ii)] the transfer time is unlimited, and
\item[iii)] the dynamics are governed by the CR3BP.
\end{enumerate}
In fact, a large body of trajectory optimization problems considered in the literature use the dynamical model of the CR3BP\cite{pritchett_low-thrust_2017, pritchett_impulsive_2018, oshima_global_2017,kelly_orthogonal_2023, chupin_transfer_2018, chupin_low-thrust_2017, caillau_minimum_2012, caillau_minimum_2012-1, mingotti_combined_2007, senent_low-thrust_2005}.
Migrating results from a CR3BP to
the more general elliptic restricted three-body problem (ER3BP) is not straightforward\cite{park_leveraging_2022, campagnola_subregions_2008, celletti_dynamics_2024}.  In fact, the ER3BP is widely considered\cite{ferrari_periodic_2018, szebehely_theory_1967, broucke_periodic_1969, broucke_stability_1969, campagnola_subregions_2008, celletti_dynamics_2024, hyeraci_method_2010} to be significantly harder than the CR3BP. This is, in part, because the dynamics of the ER3BP are non-autonomous\cite{szebehely_theory_1967,broucke_periodic_1969,broucke_stability_1969}.  Non-autonomous dynamical systems generate difficult questions with respect to periodicity, stability, etc.\cite{broucke_periodic_1969,broucke_stability_1969,campagnola_subregions_2008, celletti_dynamics_2024}.  It turns out that some of these questions can be answered\cite{broucke_periodic_1969,broucke_stability_1969} through the use of the Nechvile transformation\cite{szebehely_theory_1967,dixon-ross-TN_2025}; however, the results are subtle\cite{dixon-ross-TN_2025} and substantially more intricate than those of the CR3BP\cite{ferrari_periodic_2018,campagnola_subregions_2008, celletti_dynamics_2024}.  For instance, in the case of the ER3BP, one must explicitly account for clock times and asymmetries\cite{szebehely_theory_1967, broucke_periodic_1969, broucke_stability_1969, campagnola_subregions_2008, celletti_dynamics_2024, du_low-thrust_2023} resulting from a nonzero value of the eccentricity.  These points will be more apparent in the sections to follow.  Given these additional challenges, it is apparent that trajectory optimization in the ER3BP is substantially more intricate\cite{parker_lowenergy_2014, hiday-johnston_transfers_1994, bonnard_time-minimum_2016, neelakantan_two-impulse_2022, du_low-thrust_2023, hyeraci_method_2010, scantamburlo_interplanetary_2022} than a corresponding one in the CR3BP. As indicated in \cite{park_leveraging_2022}, solutions to the ER3BP are more desirable than those found in the CR3BP because it is simpler to migrate the results from the former (than the latter) to a full ephemeris model.  Hence, a significant advantage of solving trajectory optimization problems in the ER3BP is a closer approximation to a practical optimal trajectory.  In view of this, we use the ER3BP as a starting point for formulating trajectory optimization problems.  Clearly, our processes will apply to the CR3BP as well.

The process we choose to formulate trajectory optimization problems in the ER3BP parallels those of the classical situations in a two-body problem in that we consider the end-to-end problem of transferring a spacecraft from one libration-point orbit to another. The libration-point orbits are endpoint manifolds in an optimal-control setting\cite{ross_primer_2015,clarke_nonsmooth_1998}.  Propellant consumption is modeled using the new $\Delta{prox}$ formulation introduced in \cite{dixon-ross-TN_2025}.  This model is based on an $L^1$-functional\cite{ross_how_2004, ross_space_2006,ross_primer_2015} that is mapped to propellant consumption via the Nechvile transformation.  Similar to the $L^1$-functional, the $\Delta{prox}$ model is agnostic to the specific-impulse of the propellant and hence has broader applicability than the use of a mass flow-rate equation.  Minimum-time and time-bounded, minimum-propellant cost functionals are considered.  The necessary conditions for optimality are derived by an application of Pontryagin's Principle\cite{ross_primer_2015}.  A suite of checkable conditions are extracted out of these necessary conditions for verifying the extremality of a computed solution. Solutions are computed using the universal Birkhoff theory for trajectory optimization\cite{ross_universal_2023, ross_universal_2023-1, proulx_implementations_2023,ross_hessianTN_2025,sandia_millionPt_2025}.  Computational results for the suite of problems considered detected no singular arcs\cite{bonnard_time-minimum_2016,singular-arc-2010-CDC}.  All computed solutions were found to be verifiably bang-bang or bang-off-bang extremal arcs.

Unlike many of the prior techniques used for trajectory optimization either in the CR3BP or the ER3BP, the universal Birkhoff method is neither a direct nor an indirect method\cite{ross_universal_2023-1,ross_enhancements_2020}, although it can be used in either form\cite{ross_universal_2023-1}.  When used as neither, i.e., in its unique fast spectral form\cite{ross_universal_2023-1, proulx_implementations_2023, ross_enhancements_2020,gong_spectral_2008}, it has the following computational advantages\cite{ross_universal_2023, ross_universal_2023-1, proulx_implementations_2023}:
\begin{enumerate}
\item Because its condition number is mesh-independent\cite{ross_universal_2023-1}, the universal Birkhoff method can generate accurate solutions in a fast and stable manner even when the number of grid points is over a million\cite{sandia_millionPt_2025}.
\item Under proper, non-canonical scaling\cite{scaling,ross_universal_2023-1}, it can operate guess-free\cite{ross_enhancements_2020,ross:guess-free,LRO-CSM}.  Consequently, there is no need to seed a Birkhoff method with an initial solution from dynamical systems theory.
\item With the use of an appropriately weighted Lagrangian\cite{ross_universal_2023-1} and scaling parameters\cite{scaling}, a Birkhoff code can be produced\cite{ross_enhancements_2020,sandia_millionPt_2025} to autonomously generate and solve\cite{LRO-CSM} the Hamiltonian necessary conditions\cite{ross_primer_2015} resulting from an application of Pontryagin's Principle.
\item There is no Gibbs phenomenon if the controls have jump discontinuities\cite{ross_universal_2023-1,proulx_implementations_2023,ross_enhancements_2020}. Hence, there is no need to split the time intervals as in the case of ``$hp$'' pseudospectral knotting methods\cite{gong_spectral_2008,auto-knots,ross_pseudospectral_2004}.  This point will be obvious in Section~\ref{sec:numerics}.
\end{enumerate}
These advantages of the universal Birkhoff method suggests that it should be able to overcome many of the issues associated with the collection of ad hoc and diverse techniques mentioned earlier.  Overall, this paper provides a new approach for computing verifiable extremal trajectories in the ER3BP that require no initial guess or assistance from dynamical systems theory. High-order numerical propagation of the control solution is performed to independently validate the feasibility of the orbital transfers.

\section{Details of a Class of Trajectory Optimization Problems in the ER3BP}\label{sec:OCPs}
The main objective of this section is to formulate the precise details of a class of optimal control problems in the ER3BP.
\subsection{Dynamical Model for the ER3BP}

\begin{figure}[hbt!]
\centering
\includegraphics[width=0.7\columnwidth]{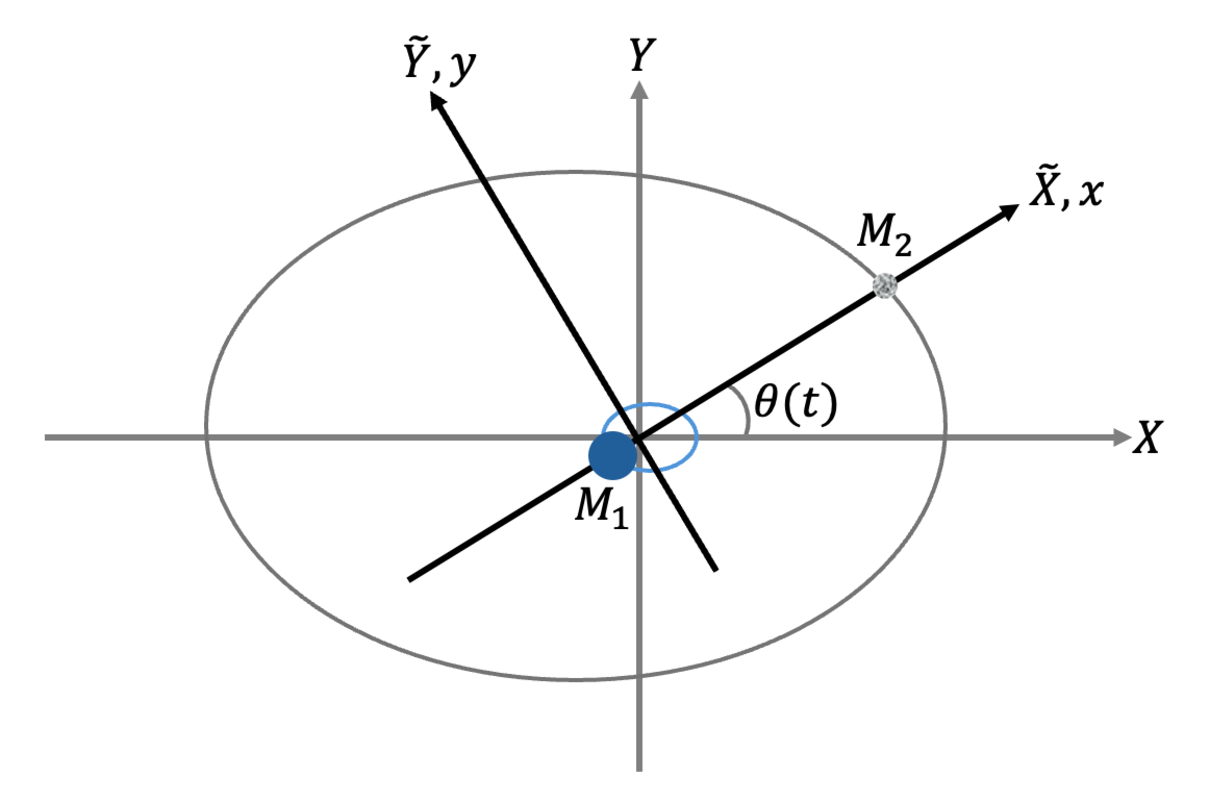}
\caption{ER3BP inertial, rotating, and Nechvile frames}
\label{fig:ER3BP_CoordFrame}
\end{figure}

The ER3BP makes assumptions similar to those in the CR3BP in that the spacecraft mass has no influence on the rotation of the two primaries, and the primary motion is constrained to a plane.  The ER3BP is obviously more general than the CR3BP because the two primaries move in elliptical orbits. The equations of motion (EOM) for a massless particle in the ER3BP are well established\cite{szebehely_theory_1967}.  As discussed in \cite{dixon-ross-TN_2025}, there are several ways to add a control force to these EOM. In this paper, we use the engine-agnostic dynamical model (Model~$C$ in \cite{dixon-ross-TN_2025}) given by
\begin{equation}
    \begin{aligned}
        x''(\theta) - 2y'(\theta) &= \partial_x\Psi(x, y, z, \theta) + \frac{ u_x}{(1+e\cos{\theta})^3}\,
         \\
        y''(\theta) + 2x'(\theta) &= \partial_y\Psi (x, y, z, \theta) + \frac{ u_y}{(1+e\cos{\theta})^3} \,
         \\
        z''(\theta) &= \partial_z\Psi(x, y, z, \theta) + \frac{ u_z}{(1+e\cos{\theta})^3} \,
    \end{aligned}
    \label{eq: DynamicsEOM}
\end{equation}
The coordinate system for \eqref{eq: DynamicsEOM} is shown in Fig.~\ref{fig:ER3BP_CoordFrame} where $XYZ$ is an inertial frame, $\tilde{X}\tilde{Y}\tilde{Z}$ is the Earth-Moon rotating frame, and $xyz$ is the Nechvile frame, which can be described as a dimensionless, rotating, and pulsating coordinate frame.  In Fig.~\ref{fig:ER3BP_CoordFrame}, the $Z$-, $\tilde{Z}$-, and $z$-axes are determined by the right-hand rule.
The EOM are given in the Nechvile frame with the true anomaly, $\theta$, as the independent variable. Hence, the prime notation is used in \eqref{eq: DynamicsEOM} to describe the derivatives instead of a dot.  In the Nechvile frame, the locations of the primaries and the collinear Lagrange points remain fixed along the $x$-axis.  The function $\Psi$ in \eqref{eq: DynamicsEOM} is the pseudo-potential function given by
\begin{multline}
    \Psi(x,y,z, \theta) := \frac{1}{1+e\cos{\theta}}\left[\frac{1}{2}(x^2+y^2-z^2e\cos{\theta}) \right.\\
    \left. + \frac{1-\mu}{\rho_1(x,y,z)} + \frac{\mu}{\rho_2(x,y,z)} \right]
    \label{eq: PseudoPotentialFunction}
\end{multline}
where $\rho_1$ and $\rho_2$ are the distance functions between the spacecraft and the two primaries such that
\begin{equation}
    \begin{aligned}
        \rho_1(x, y, z) &:= \sqrt{(x+\mu)^2 + y^2 + z^2} \\
        \rho_2(x, y, z) &:= \sqrt{(x-1+\mu)^2 + y^2 + z^2} \,,
    \end{aligned}
\end{equation}
$\mu$ is the mass ratio of the primaries, and $e$ is the eccentricity of the two-body system. The controls $u_x$, $u_y$ and $u_z$ are the components of a thrust-force acceleration that is scaled by $G(M_1 + M_2)/p^2$, where $p$ is the semilatus rectum of elliptic orbit of the primaries of mass $M_1$ and $M_2$.  See \cite{dixon-ross-TN_2025} for complete details.  In the discussions to follow, we use the term ``time unit" (TU) synonymously with the independent variable of true anomaly in the Nechvile frame. This minor abuse of terminology aides in the discussion of state and control time histories instead of state and control ``true anomaly" histories.
\begin{remark}
In \eqref{eq: DynamicsEOM}, the controls are divided by a cubic cosine term.  Obviously, this term is unity if $e=0$. Ignoring this term for $e \ne 0$ even for small values of eccentricity (e.g., $e=0.0549$) can generate wildly different trajectories.  Even impulsive $\Delta V$ maneuvers are not immune to this effect.  See \cite{dixon-ross-TN_2025} for further details.
\end{remark}
%------------------------

\subsection{The $\Delta{prox}$-Functional Model for Propellant Consumption}
Because propellant consumption is a key criterion for space maneuvers, it is typically modeled using the mass-flow-rate equation
\begin{equation}\label{eq:rocket}
\dot m = - \frac{T}{I_{sp} g_0}
\end{equation}
where, $T$ is the thrust force produced by a thruster, $I_{sp}$ is the specific impulse of the propellant and $g_0$ is gravitational acceleration at sea level. One problem in using \eqref{eq:rocket} for modeling propellant consumption is that the resulting analysis becomes quite limited to the numerical choice of $I_{sp}$.  For instance, by choosing a high value of $I_{sp}$, one can generate a small amount of propellant consumption.  Conversely, a smaller value of $I_{sp}$ generates a larger consumption of propellant.   To remove this limitation, it was suggested in \cite{ross_how_2004} and \cite{ross_space_2006} that the $L^1$-functional of the control acceleration be used as an $I_{sp}$-agnostic metric for propellant consumption. Because the control acceleration is a vector, it is possible to define an infinite set of $\ell^p$-classes of the $L^1$ functional given by\cite{ross_space_2006,ross_how_2004,dixon-ross-TN_2025}
\begin{equation}\label{eq:proxByDef}
\Delta{prox}_p:=  \int_{t_0}^{t_f} \|\bu(t)\|_p  \,dt
\end{equation}
for $p=1, 2, \ldots$. The notation $\Delta prox$ in \eqref{eq:proxByDef} is used to clarify the fact that the the $L^1$-norm does not compute the actual propellant but a proxy that accurately measures the $I_{sp}$-agnostic propellant consumption\cite{ross_how_2004}.  Using the Nechvile transformation and scaling the result by  $\sqrt{G(M_1 + M_2)/p^3}$, it can be shown that \eqref{eq:proxByDef} transforms to\cite{dixon-ross-TN_2025}
\begin{equation}
    \Delta{prox}_p =  \int_{\theta_0}^{\theta_f} \frac{\|\bu(\theta)\|_p}{(1+e\cos{\theta})^2} \,d\theta
    \label{eq:DeltaProx4ER3BP}
\end{equation}
As shown in \cite{ross_how_2004}, using $p = 1, 2$ and $\infty$ provide a direct connection to certain thruster configurations.  In this context, the ``$L^1$-$\ell^1$'' proxy functional is given by
\begin{multline}
        \Delta{prox}_{1} :=\int_{\theta_0}^{\theta_f} \frac{\|\bu(\theta)\|_1}{(1+e\cos{\theta})^2} \,d\theta \\
        := \int_{\theta_0}^{\theta_f}\frac{\left(|u_x(\theta)| + |u_y(\theta)| + |u_z(\theta)|\right)}{(1+e\cos{\theta})^2}\,d\theta
    \label{eq:L1-l1}
\end{multline}
Similarly, the ``$L^1$-$\ell^2$'' proxy functional is given by
\begin{multline}\label{eq:L1-l2}
        \Delta{prox}_{2} :=\int_{\theta_0}^{\theta_f} \frac{\|\bu(\theta)\|_2}{(1+e\cos{\theta})^2} \,d\theta \\
        := \int_{\theta_0}^{\theta_f}\frac{\sqrt{u_x^2(\theta) + u_y^2(\theta) + u_z^2(\theta)}}{(1+e\cos{\theta})^2}\,d\theta
\end{multline}
As noted in \cite{ross_primer_2015,ross_space_2006,dixon-ross-TN_2025} the $\Delta{prox}_{2}$ formulation is not quadratic.  Squaring it does not generate a quadratic functional.  Quadratic cost functionals are popular \cite{mingotti_combined_2007, mingotti_low-energy_2009,  chupin_low-thrust_2017 ,chupin_transfer_2018, du_transfer_2022, du_low-thrust_2023} because they are perceived to be easier to manage from a mathematical/computational perspective. Unfortunately, taking a quadratic functional can result in solutions that are markedly different than minimum-propellant solutions\cite{ross_space_2006,dixon-ross-TN_2025}.  Even worse, the quadratic functional can consume as much as $50\%$ excess propellant in the generic case\cite{ross_how_2004} and over $100\%$ more in the case of the ER3BP\cite{dixon-ross-TN_2025}.

A mathematical and computational challenge in directly using any of the $\Delta{prox}$ functionals is that they are nonsmooth\cite{clarke_nonsmooth_1998, chupin_low-thrust_2017}. In fact, the point of nondifferentiablity is exactly at the ``worst'' point, namely the most desirable ``point'' of $\bu = \bzero$\cite{ross_how_2004,ross_primer_2015}.  That is, from a propellant-savings point of view, the zero-thrust arcs $\theta \mapsto \bu(\theta) = \bzero$ are most desirable. Yet, at exactly all of these desirable points, the integrands of both \eqref{eq:L1-l1} and \eqref{eq:L1-l2} are nondifferentiable.  In other words, the nondifferentiability at zero (in the integrand) is spread over an entire arc and not isolated to a single point.  Fortunately, the nondifferentiablity can be easily removed for both equations using suitable transformations\cite{ross_how_2004,ross_primer_2015}.  In the case of \eqref{eq:L1-l1}, this transformation is given by replacing $u_x, u_y$ and $u_z$ in \eqref{eq: DynamicsEOM} according to
\begin{subequations}\label{eq:L1-to-smooth}
    \begin{align}
u_x&:= u_x^a - u_x^b &  u_x^a \ge 0,  u_x^b \ge 0\\
u_y&:= u_y^a - u_y^b  &  u_y^a \ge 0,  u_y^b \ge 0\\
u_z&:= u_z^a - u_z^b  &  u_z^a \ge 0,  u_z^b \ge 0
\end{align}
\end{subequations}
while simultaneously rewriting $\Delta{prox}_{1}$ as
\begin{multline}\label{eq:prop-11-smooth}
    \Delta{prox}_{1} = \\
     \int_{\theta_0}^{\theta_f}\frac{\left({u_x^a(\theta) + u_x^b(\theta) + u_y^a(\theta) + u_y^b(\theta) + u_z^a(\theta) + u_z^b(\theta)}\right)}{(1+e\cos{\theta})^2}\,d\theta
\end{multline}
A clear advantage of \eqref{eq:prop-11-smooth} is that there is no need to use a continuation or homotopy approach\cite{mingotti_low-energy_2009, chupin_transfer_2018, du_transfer_2022} wherein a sequence of problems are solved to gradually refine the solution from a quadratic cost to the more desirable correct solution.

In the case of \eqref{eq:L1-l2}, the nonsmooth-to-smooth transformation involves keeping $u_x, u_y$ and $u_z$ in \eqref{eq: DynamicsEOM} unchanged while rewriting $\Delta{prox}_{2}$ as\cite{ross_how_2004,ross_primer_2015}
\begin{equation}\label{eq:u4-integrand}
    \Delta{prox}_{2} = \int_{\theta_0}^{\theta_f}\frac{u_4(\theta)}{(1+e\cos{\theta})^2} \,d\theta
\end{equation}
where $u_4(\theta)$ is constrained by an algebraic equation and an inequality constraint given by\cite{ross_primer_2015}
\begin{equation}\label{eq:u4-constraint}
u_x^2(\theta) + u_y^2(\theta) + u_z^2(\theta)- u_4^2(\theta) = 0, \quad u_4(\theta) \ge 0
\end{equation}
Note that \eqref{eq:u4-integrand} and \eqref{eq:u4-constraint} are based on squaring the square-root inside the integral in \eqref{eq:L1-l2}; hence, the result is not a transformed quadratic functional.
\begin{remark}
Through the use of \eqref{eq:L1-to-smooth}--\eqref{eq:u4-constraint}, it is obvious that both of the $\Delta{prox}_1$ and $\Delta{prox}_2$ functionals can be instantiated in one shot without homotopy.  See also \cite{dixon-ross-TN_2025} for an outsized impact of the quadratic cosine term in \eqref{eq:DeltaProx4ER3BP}.
\end{remark}
To limit the scope of the numerical discussions, we use \eqref{eq:prop-11-smooth} in Section~\ref{sec:numerics}.

\subsection{Algebraic Modeling of Boundary Conditions}
\begin{figure}[h!]
\centering
\includegraphics[width=0.7\columnwidth]{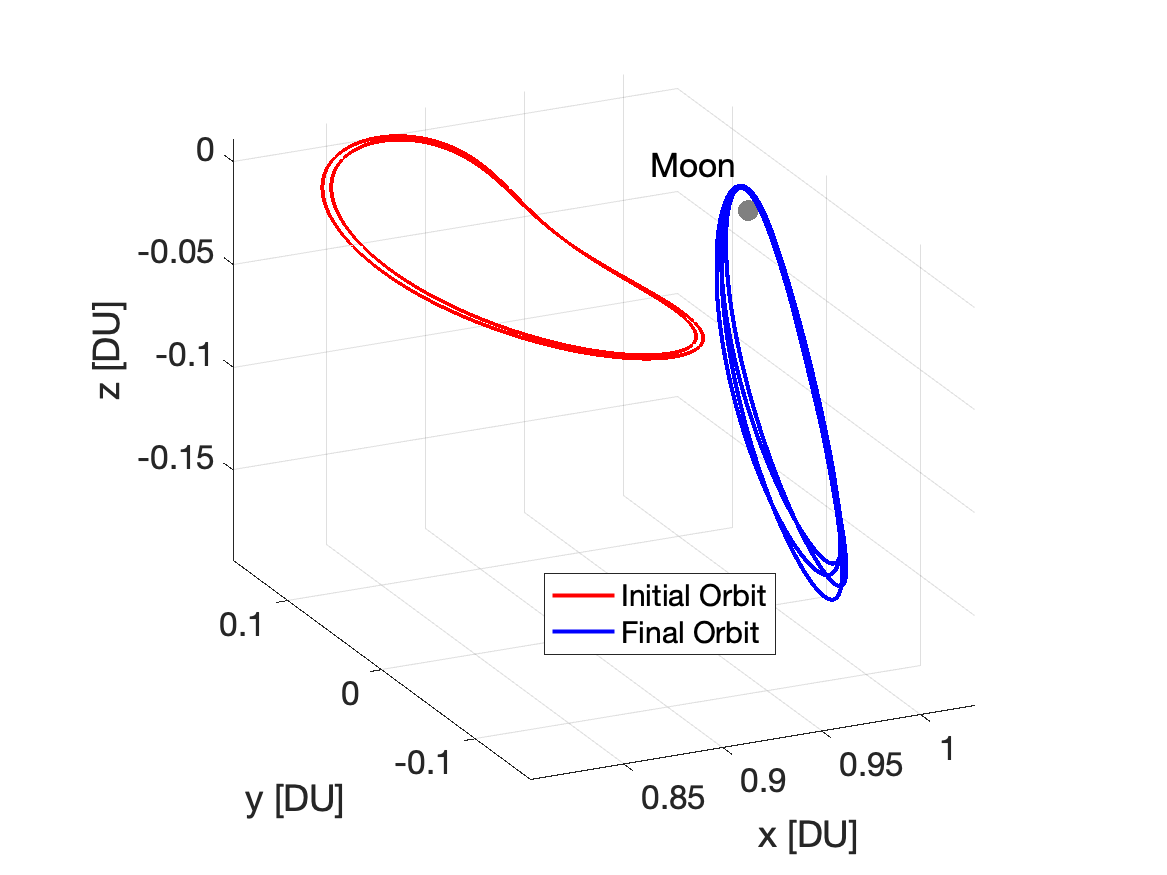}
\caption{Sample initial and final orbits for trajectory optimization in the ER3BP.}
\label{fig:ER3BP_BoundaryOrbits}
\end{figure}

The geometry of typical boundary conditions considered in this paper is shown in Fig. \ref{fig:ER3BP_BoundaryOrbits}.
%
%%%%%%%%%%%%%%%%%%%%%%%%%%%%%%%%%%%%%%%%%%%%%%%%%%%%%%%%%%%%%%%%%%%%%%%%%%%%%
%%%%%%%%%%%%%%%%%%%%%%%%%%%%%%%%%%%%%%%%%%%%%%%%%%%%%%%%%%%%%%%%%%%%%%%%%%%%%
%
The initial condition is any point on a 2:1 resonant $L_1$ Lyapunov orbit while the final point is any point on
a 4:1 resonant southern near-rectilinear halo orbit (NRHO)\cite{Farq-1968,NRHO-2017}.
In other words, the orbits shown in Fig.~\ref{fig:ER3BP_BoundaryOrbits} are endpoint manifolds from an optimal control perspective\cite{ross_primer_2015,clarke_nonsmooth_1998}.  These orbits
were computed through a differential correction process similar to the methods used in other studies \cite{ferrari_periodic_2018,celletti_dynamics_2024}.

\begin{figure}[h!]
\centering
\includegraphics[width=0.7\columnwidth]{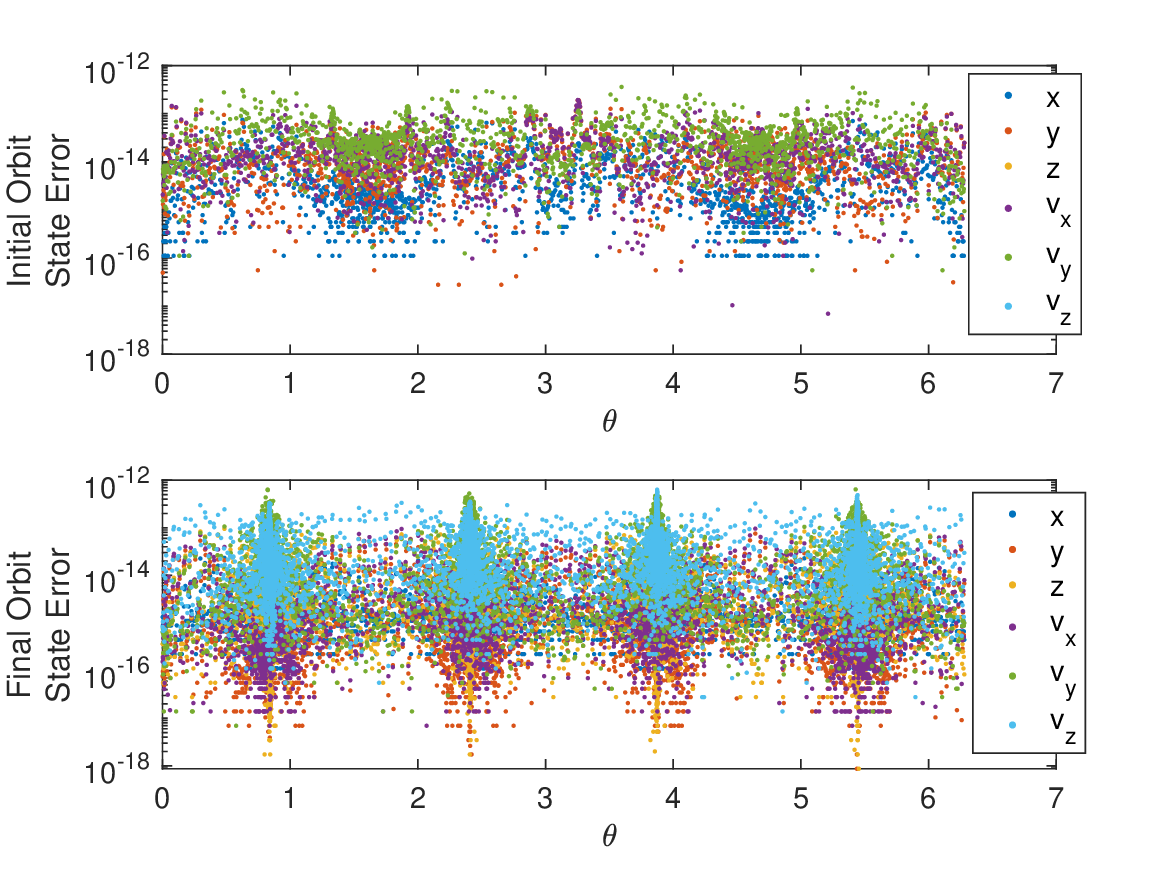}
\caption{Difference between Chebyshev and Runge-Kutta solutions for the initial and final boundary orbits shown in Fig.~\ref{fig:ER3BP_BoundaryOrbits}.}
\label{fig:ER3BP_Chebyshev_OrbitStateError}
\end{figure}

Unlike trajectory optimization in the two-body problem, the orbits shown in Fig.~\ref{fig:ER3BP_BoundaryOrbits} are not expressible in terms of integrals of motion (like semi-major axis or inclination); however, they are computable to high precision by combining a high-order Runge-Kutta (RK) method with a shooting method\cite{ferrari_periodic_2018}.
To utilize the tools of optimal control theory, it is necessary to express these orbits in terms of algebraic equations according to\cite{ross_primer_2015,longuski_optimal_2014}
\begin{subequations}
\begin{align}
    \be_1(\bx_0,\, \theta_0) &= \mathbf{0} \label{eq:e1=0}\\
    \be_2(\bx_f,\, \theta_f) &= \mathbf{0} \label{eq:e2=0}
    \end{align}
\end{subequations}
where $\bx := (x, y, z, v_x, v_y, v_z) \in \real{6}$ is the state variable and $\be_1$ and $\be_2$ are algebraic functions of $\bx$ and $\theta$ whose zeros generate the orbits shown in Fig.~\ref{fig:ER3BP_BoundaryOrbits}.  To achieve this task, we sample the values of the state vector at the Chebyshev-Gauss-Lobatto (CGL) points\cite{trefethen_approximation_2020,proulx_implementations_2023,boyd_chebyshev_2000} to construct a Lagrange polynomial that ``approximates'' the RK solution.  Deferring the details of this procedure for the moment, Fig.~\ref{fig:ER3BP_Chebyshev_OrbitStateError} shows the difference, indicated by ``state error,'' between the Chebyshev approximation and the RK solution.
%%%%%%%%%%%%%%%%%%%%%%%%%%%%%%%%%%%%%%%%%%%%%%%%%%%%%%%%%%%%%%%%%%%%%%%%%%%%%
%%%%%%%%%%%%%%%%%%%%%%%%%%%%%%%%%%%%%%%%%%%%%%%%%%%%%%%%%%%%%%%%%%%%%%%%%%%%%
%
Obviously, these ``errors'' are close to machine precision.  In other words, according to Fig.~\ref{fig:ER3BP_Chebyshev_OrbitStateError}, the Chebyshev solution is no less accurate than the RK solution given that both are numerically computed. In principle, this result should not be surprising because a Chebyshev polynomial is close to the ``best'' approximating polynomial\cite{trefethen_approximation_2020,boyd_chebyshev_2000}.  Furthermore, a Chebyshev polynomial can be very easily generated  using a Lagrange interpolant over $[0, 2\pi]$ using a CGL grid given by
\begin{equation}\label{eq:CGL-grid}
\theta_j = \pi \cos(\pi j/N) + \pi, \quad j = 0, 1, \ldots, N
\end{equation}
In sharp contrast to Legendre node points\cite{fahroo_advances_2008,ross_review_2012}, be it Legendre-Gauss, Legendre-Gauss-Radau or Legendre-Gauss-Lobatto, the CGL grid is given in closed form and requires a single line of code. These attractive computational features form the basis of the Chebyshev pseudospectral (PS) method\cite{fahroo:cheb-jgcd,cheb-costate}.

The CGL interpolant $y^N(\theta)$ of a known function $[0, 2\pi] \ni \theta \mapsto y(\theta)$  can be written in its equivalent spectral form as\cite{trefethen_approximation_2020,boyd_chebyshev_2000}
\begin{equation}\label{eq:ChebFunN}
y^N(\theta) = \sum_{j=0}^{N}a_j\,T_j(\theta/\pi - 1 )
\end{equation}
where $T_j$ is the Chebyshev polynomial of degree $j$ and $a_j, j = 0, \ldots, N$ are the spectral coefficients.  It can be shown\cite{boyd_chebyshev_2000,gong_spectral_2008} that that the value of $\abs{a_N}$ is a good estimate of the error
$$\max_{0 \le \theta \le 2\pi}\abs{y(\theta) - y^N(\theta)}$$
The coefficients $a_j, j = 0, 1, \ldots, N$ can be computed via the Fast Fourier Transform (FFT) of the sampled values at the CGL points \cite{trefethen_approximation_2020}.  Then, based on an error tolerance and the value of $\abs{a_N}$, $N$ can be quickly determined to generate accurate approximations such as the one shown in Fig.~\ref{fig:ER3BP_Chebyshev_OrbitStateError}.  This process forms the basis of the spectral algorithm\cite{gong_spectral_2008} for the Chebyshev PS method\cite{fahroo:cheb-jgcd,cheb-costate,ross_review_2012}. A more sophisticated version of this procedure is implemented in Chebfun\cite{t_a_driscoll_chebfun_2014}, a MATLAB${}^{\text{TM}}$ tool that automates a wide variety of mathematical operations on the basis of the excellent approximation properties of a Chebyshev polynomial.  In fact, Chebfun was used to generate Fig.~\ref{fig:ER3BP_Chebyshev_OrbitStateError}.

\subsubsection{Results for the 2:1 Resonant $L_1$ Lyapunov Orbit }
\begin{figure}[hbt!]
\centering
\subfigure[Lyapunov orbit position state variables]{%
    \includegraphics[width=0.7\columnwidth]{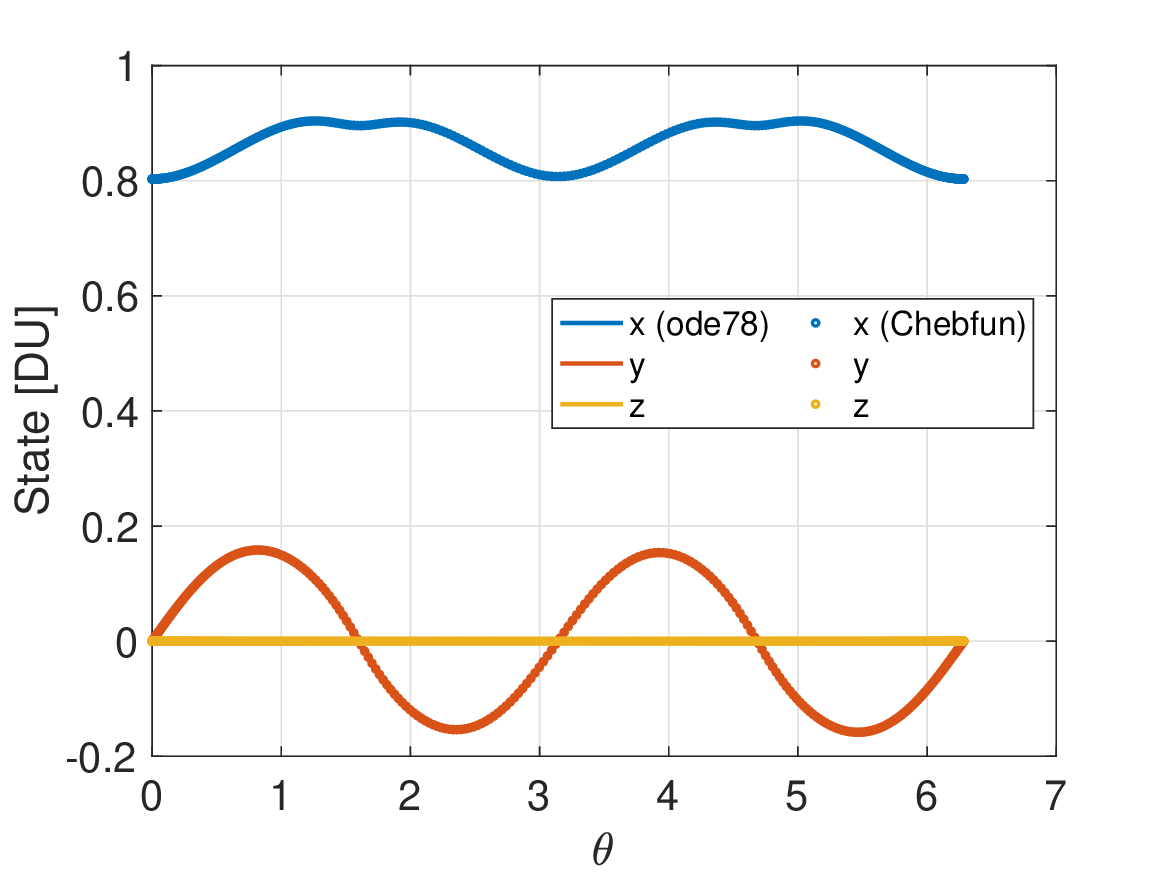}
    \label{fig:ChebyshevToODE_InitialOrbitPosStates}
}%
\\
\subfigure[Lyapunov orbit velocity state variables]{%
    \includegraphics[width=0.7\columnwidth]{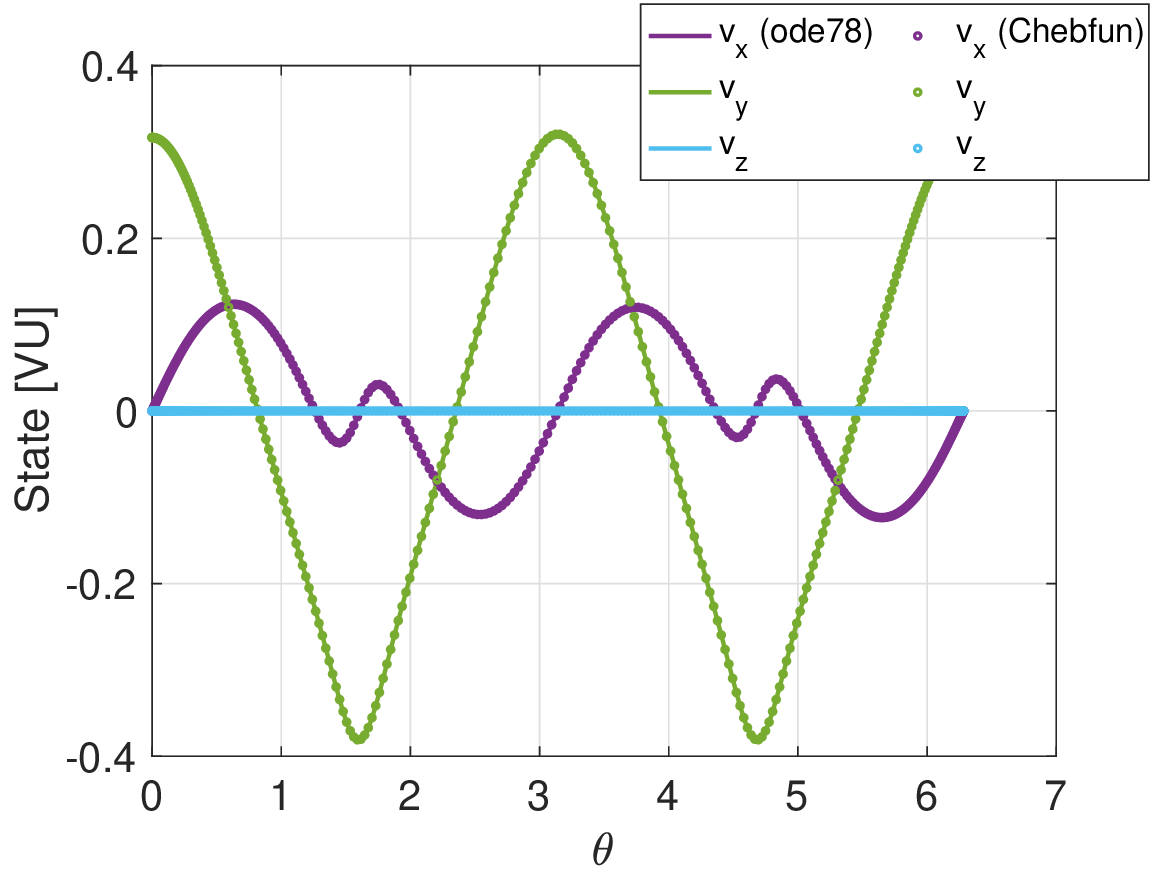}
    \label{fig:ChebyshevToODE_InitialOrbitVelStates.eps}
}%
\caption{Chebyshev and RK78 solution for the 2:1 resonant $L_1$ Lyapunov orbit shown in Fig.~\ref{fig:ER3BP_BoundaryOrbits}.}
\label{fig:ChebyshevToODE_L1Ly}
\end{figure}

Figure~\ref{fig:ChebyshevToODE_L1Ly} shows a Chebyshev solution obtained by Chebfun overlaid with an RK 7th-8th order (RK78) solution.  Based on the results shown in Fig.~\ref{fig:ER3BP_Chebyshev_OrbitStateError}, it is no surprise that the two solutions shown in Fig.~\ref{fig:ChebyshevToODE_L1Ly} are indistinguishable.
%
%%%%%%%%%%%%%%%%%%%%%%%%%%%%%%%%%%%%%%%%%%%%%%%%%%%%%%%%%%%%%%%%%%%%%%%%%%%%%%
%%%%%%%%%%%%%%%%%%%%%%%%%%%%%%%%%%%%%%%%%%%%%%%%%%%%%%%%%%%%%%%%%%%%%%%%%%%%%%
%
However, in the case of the Chebfun solution, the $L_1$ Lyapunov orbit has an analytical expression given in the form of \eqref{eq:ChebFunN}.  We represent this solution as
\begin{equation}\label{eq:x0=cheb}
\bx_0 = \bcheb_0(\theta_0), \quad 0 \le \theta_0 \le 2\pi
\end{equation}
As a result, \eqref{eq:e1=0} can now be represented as
\begin{equation}\label{eq:e1==cheb}
\be_1(\bx_0, \theta_0) = \bx_0 - \bcheb_0(\theta_0) = \bzero, \quad 0 \le \theta_0 \le 2\pi
\end{equation}
Implicit in \eqref{eq:x0=cheb} and hence \eqref{eq:e1==cheb} is a value of $N$. In the case of Fig.~\ref{fig:ChebyshevToODE_L1Ly}, $N = 308$ was sufficient to generate the accuracy depicted in Fig.~\ref{fig:ER3BP_Chebyshev_OrbitStateError}.  Note that although the $L_1$ Lyapunov orbit is 2:1 resonant, it is not periodic over $[0, \pi]$.  This point is apparent by a closer inspection of Fig.~\ref{fig:ER3BP_BoundaryOrbits}.  As a result, it would be incorrect to generate \eqref{eq:x0=cheb} over $[0, \pi]$.

\subsubsection{Results for the 4:1 Resonant Southern NRHO }
\begin{figure}[hbt!]
\centering
\subfigure[NRHO position state variables]{%
    \includegraphics[width=0.7\columnwidth]{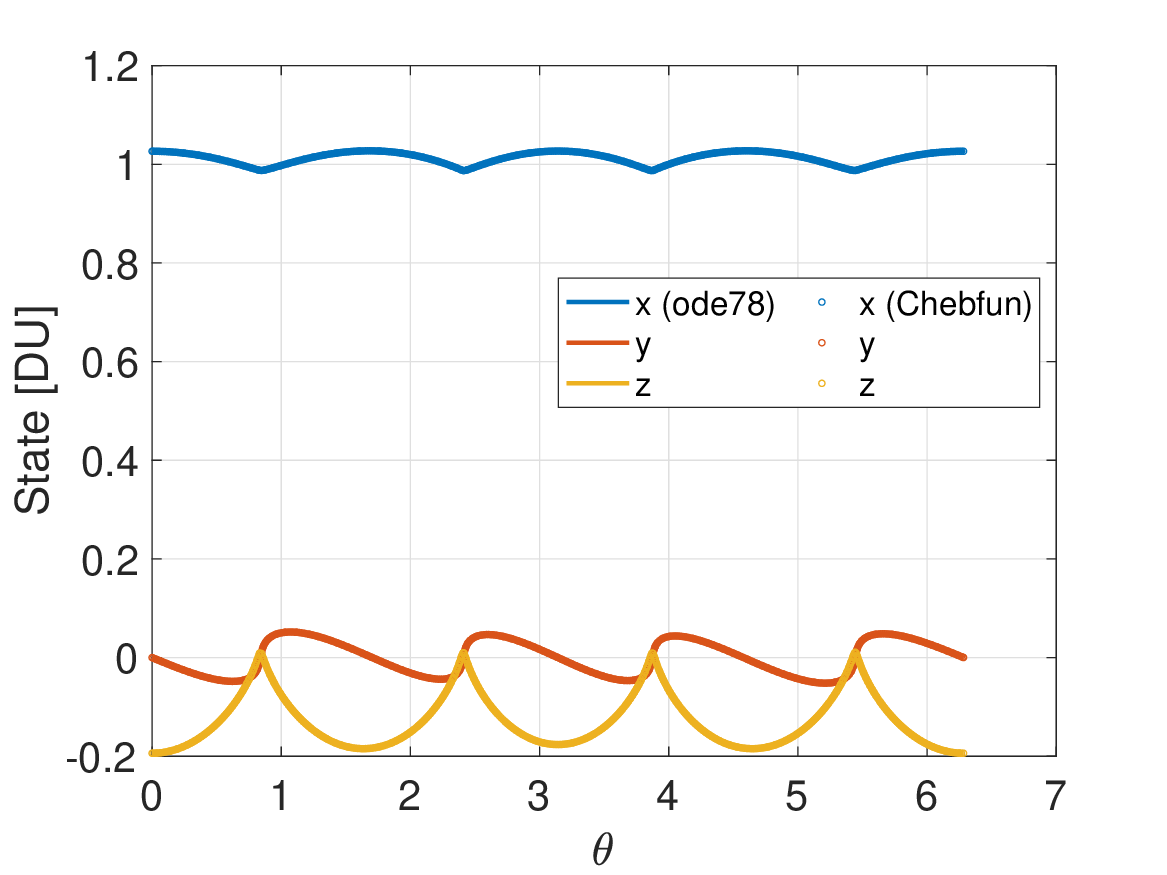}
    \label{fig:ChebyshevToODE_FinalOrbitPosStates}
}%
\\
\subfigure[NRHO velocity state variables]{%
    \includegraphics[width=0.7\columnwidth]{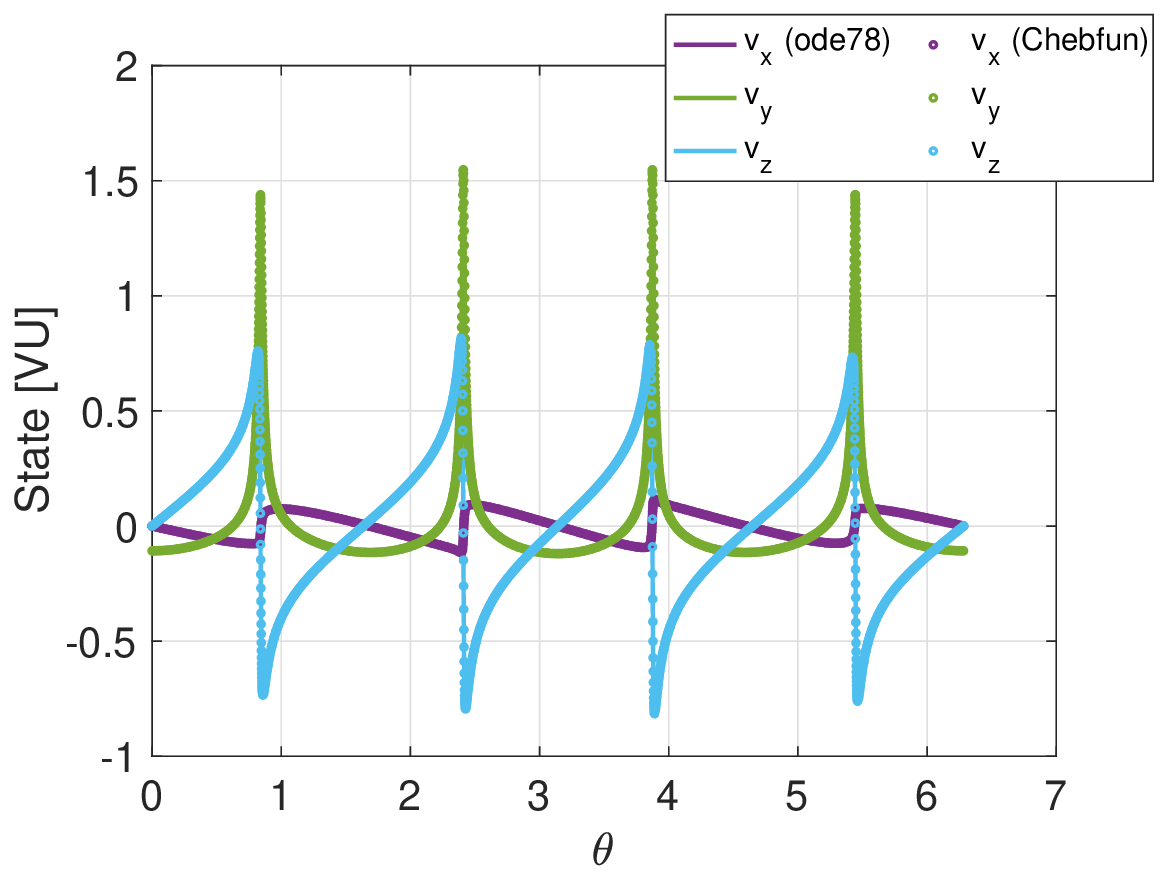}
    \label{fig:ChebyshevToODE_FinalOrbitVelStates.eps}
}
\caption{Chebyshev and RK78 solution for the 4:1 resonant southern NRHO shown in Fig.~\ref{fig:ER3BP_BoundaryOrbits}}
\label{fig:ChebyshevToODE_NRHO}
\end{figure}

Performing an exercise similar to that conducted for the 2:1 resonant $L_1$ Lyapunov orbit yields the plots shown in Fig.~\ref{fig:ChebyshevToODE_NRHO} for the NRHO.  That the Chebfun solution is indistinguishable from the RK78 propagation follows from the results of Fig.~\ref{fig:ER3BP_Chebyshev_OrbitStateError}.
%
%%%%%%%%%%%%%%%%%%%%%%%%%%%%%%%%%%%%%%%%%%%%%%%%%%%%%%%%%%%%%%%%%%%%%%%%%%%%%%
%%%%%%%%%%%%%%%%%%%%%%%%%%%%%%%%%%%%%%%%%%%%%%%%%%%%%%%%%%%%%%%%%%%%%%%%%%%%%%
%
Similar to \eqref{eq:x0=cheb} we represent the solution shown in Fig.~\ref{fig:ChebyshevToODE_NRHO} as
\begin{equation}\label{eq:xf=cheb}
\bx_f = \bcheb_f(\beta), \quad 0 \le \beta \le 2\pi
\end{equation}
Unlike the case of the initial orbit, \eqref{eq:xf=cheb} cannot be written with $\beta$ replaced by $\theta_f$.  This is because in the case of the initial condition we can take $\theta_0 \in [0, 2\pi]$ without loss in generality.  In the case of the final ``time'' condition, we cannot assume  $\theta_f \in [0, 2\pi]$.  Recall that $\theta$ is an independent variable and a proxy for time in the ER3BP. Because of this modulo $2\pi$ issue, the relationship between $\beta$ and $\theta_f$ can be written as
\begin{equation}\label{eq:beta=upperSemi}
\beta(\theta_f) := \theta_f - \left\lfloor\frac{\theta_f}{2\pi}\right\rfloor \times (2\pi)
\end{equation}
where we have abused notation in favor of clarity by reusing the symbol $\beta$ to also denote a function of $\theta_f$.
\begin{remark}
The function $\theta_f \mapsto \beta(\theta_f)$ is a discontinuous sawtooth function.  In fact, it is upper semicontinuous.  Nonetheless, unlike the case of \eqref{eq:DeltaProx4ER3BP}, there is no need to transform \eqref{eq:beta=upperSemi} to a smooth function for a numerical instantiation because the derivative $\partial_{\theta_f} \beta(\theta_f)$ is continuous.
\end{remark}
With this understanding, \eqref{eq:e2=0} is given by
\begin{equation}\label{eq:e2==cheb}
\be_2(\bx_f, \theta_f) = \bx_f - \bcheb_f(\beta(\theta_f)) = \bzero, \quad \theta_f > \theta_0
\end{equation}
In the case of $\bcheb_f$, $N = 9,778$ was necessary to produce the accuracy depicted in Fig.~\ref{fig:ER3BP_Chebyshev_OrbitStateError}.  The reason this value of $N$ is significantly higher than the $L_1$ Lyapunov orbit is quite apparent by inspection of Fig.~\ref{fig:ChebyshevToODE_FinalOrbitVelStates.eps}.  The four ``spikes'' in velocity components correspond to passage near the perilune.  From the physics of the problem (see \eqref{eq: DynamicsEOM}), these spikes are not discontinuities; rather, they are smooth and differentiable even though they do not appear to be by visual inspection of Fig.~\ref{fig:ChebyshevToODE_FinalOrbitVelStates.eps}.
\newpage
\begin{remark}
As noted in \cite{trefethen_approximation_2020,trefethen-myths-2011}, large values of $N$ are often stupefying to many engineers because of the widespread misinformation of polynomial approximation.  Perhaps a reader of this paper might be surprised to learn that $\bcheb_f$, corresponding  to $N \approx 10, 000$, was generated in less than one tenth of a second on a MacBook Pro.  See \cite{ross_hessianTN_2025} for further analysis of computational complexity.  More importantly, the compute times for function evaluations of \eqref{eq:e1==cheb} and \eqref{eq:e2==cheb} are not quite measurable; i.e., they are instantaneous.  These advances are part of the reason why it was announced in \cite{koeppen_fast_2019} that mesh refinement can now be performed at $\mathcal{O}(1)$ computational speed; i.e., instantaneously, even on a low end processor.
\end{remark}

\subsection{Orbit Targeting Constraints}
\begin{figure}[h!]
\centering
\includegraphics[width=0.7\columnwidth]{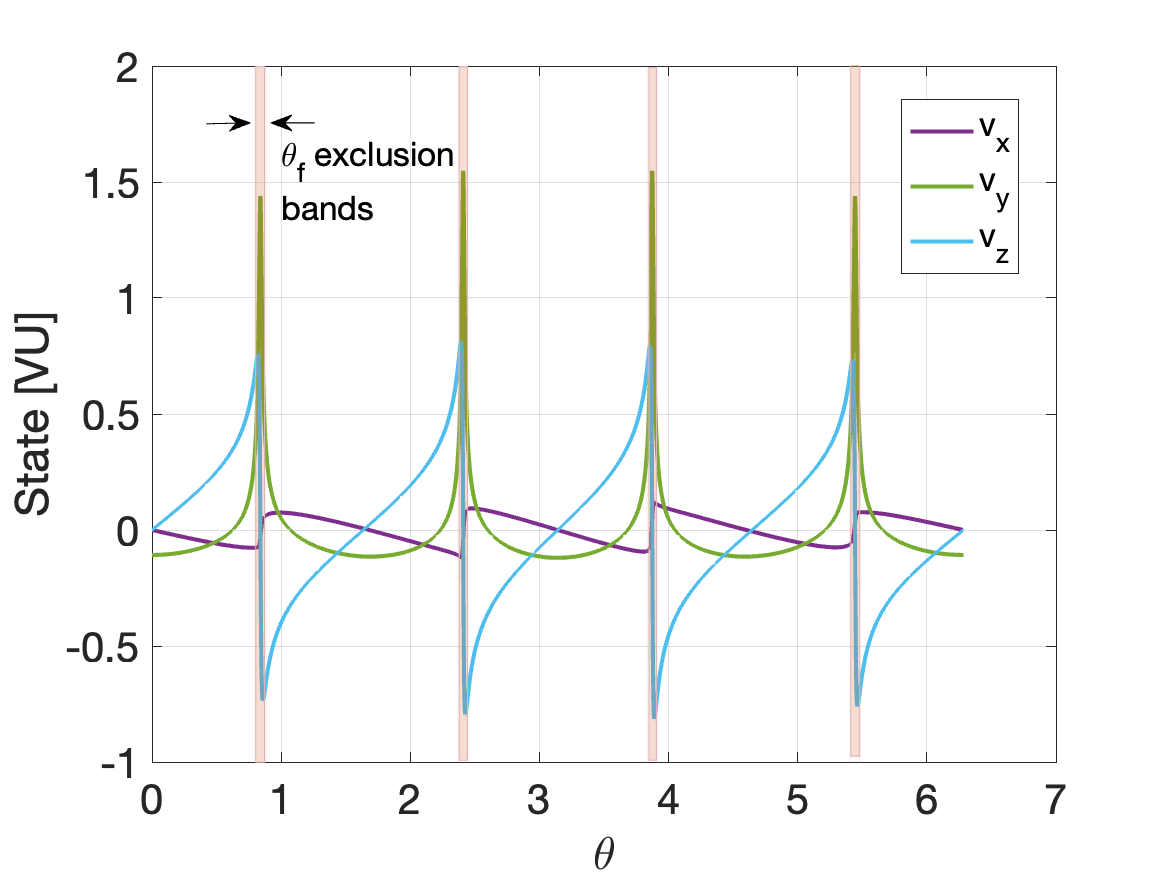}
\caption{True anomaly exclusion bands around areas of extreme sensitivities in velocity space for the NRHO}
\label{fig:NRHO_ZfMax_VelStates_withThetaBands}
\end{figure}

\begin{figure}[hbt!]
\centering
\subfigure[Final admissible position states]{%
    \includegraphics[width=0.7\columnwidth]{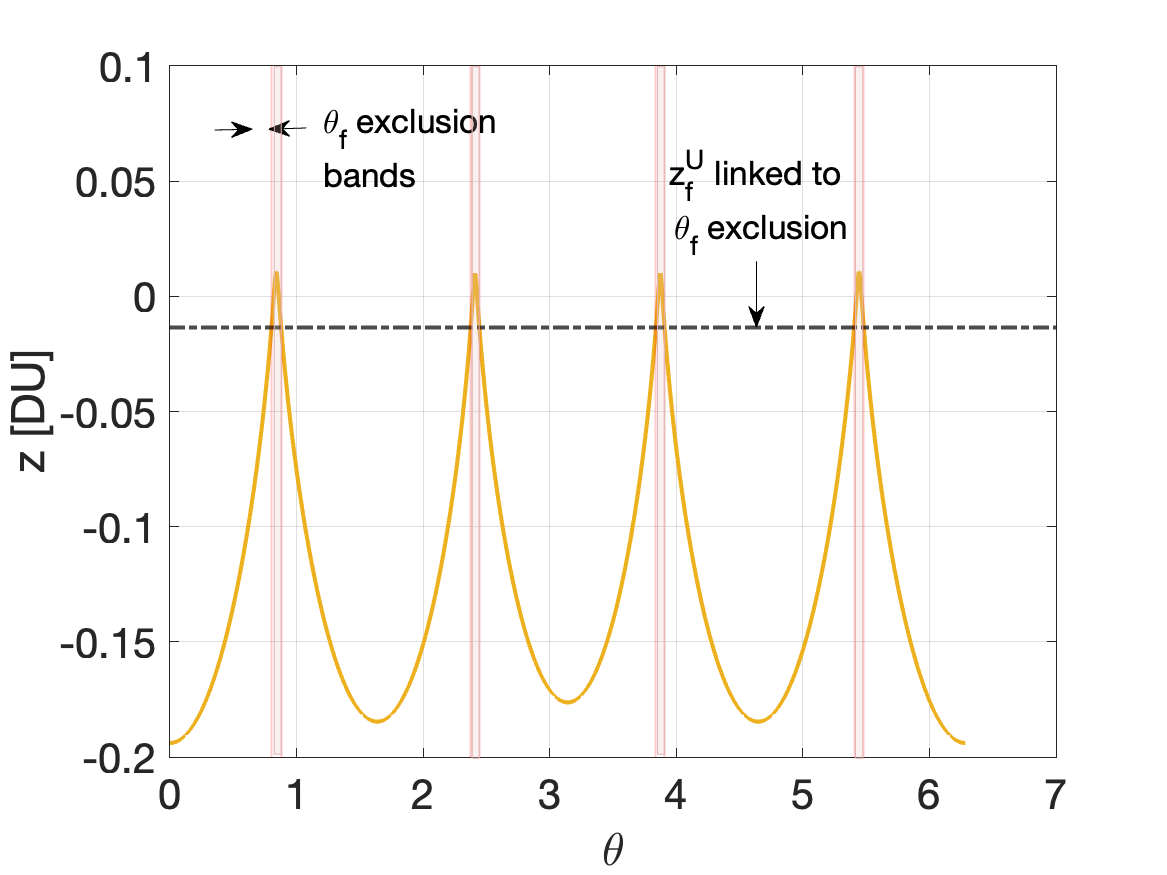}
    \label{fig:NRHO_ZfMax_PosStates_withThetaBands}
}%
\\
\subfigure[Physical interpretation]{%
    \includegraphics[width=0.7\columnwidth]{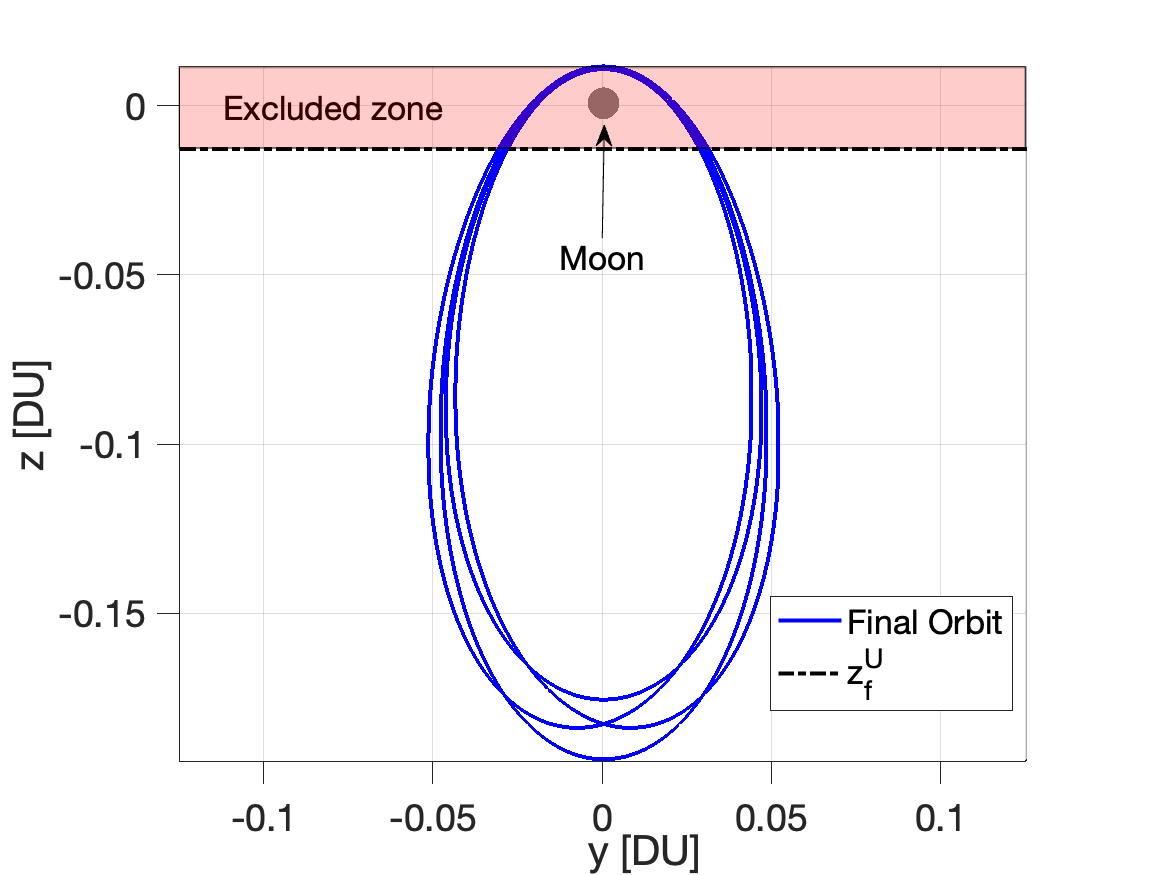}
    \label{fig:NRHO_ZfMax}
}
\caption{Orbit targeting constraints for the southern NRHO.}
\label{fig:NRHO_ZfStatesWithNRHO_ZfConstraint}
\end{figure}

As observed in the prior subsection, there are rapid changes in the velocity components at the perilune point of the NRHO. A quick examination of  Fig. \ref{fig:ChebyshevToODE_FinalOrbitVelStates.eps} shows that the most dramatic change is in the $y$-component of the velocity vector.  From a practical guidance perspective, it is desirable to avoid targeting a neighborhood around the perilune point of the NRHO. Otherwise, the natural sensitivity of this point would place a high targeting demand on a spacecraft guidance system leading to increased cost and additional risk. Consequently, practical considerations demand that neighborhoods around these high sensitivity points be excluded for targeting. This requirement is indicated in  Fig. \ref{fig:NRHO_ZfMax_VelStates_withThetaBands}.
%
%
%%%%%%%%%%%%%%%%%%%%%%%%%%%%%%%%%%%%%%%%%%%%%%%%%%%%%%%%%%%%%%%%%%%%%%%%%%%%%
%%%%%%%%%%%%%%%%%%%%%%%%%%%%%%%%%%%%%%%%%%%%%%%%%%%%%%%%%%%%%%%%%%%%%%%%%%%%%
%
Mathematically, this exclusion zone can be simply stated as,
\begin{equation}\label{eq:v-constr}
 v^L_\xi \le v_\xi \le v^U_\xi, \quad \xi = x, y, z
\end{equation}
where the lower bound $v^L_\xi$ and the upper bound $v^U_\xi$ can be chosen based on the selected ``width'' of the exclusion zone shown in Fig.~\ref{fig:NRHO_ZfMax_VelStates_withThetaBands}.  A more careful examination of Fig.~\ref{fig:NRHO_ZfMax_VelStates_withThetaBands} reveals that a low value of $v^U_y$ or $v^U_z$ would be sufficient to enforce all three velocity constraints indicated in  \eqref{eq:v-constr}.    Taking a step further, we examine the exclusion band in the position space. This is shown in Fig.~\ref{fig:NRHO_ZfMax_PosStates_withThetaBands} for just the $z$-position. From this figure it follows that a single $z$-position constraint at the final time given by
\begin{equation}
z(\theta_f) \leq z_f^{U} = -0.013556
    \label{eq: zF_constr}
\end{equation}
would be sufficient to enforce the targeting exclusion zone constraint.  Fig.~\ref{fig:NRHO_ZfMax} shows the simple physical meaning of this exclusion zone constraint.  This constraint effectively filters out the high-cost, high-risk targeting zones of the NRHO.
%
%%%%%%%%%%%%%%%%%%%%%%%%%%%%%%%%%%%%%%%%%%%%%%%%%%%%%%%%%%%%%%%%%%%%%%%%%%%%%%
%%%%%%%%%%%%%%%%%%%%%%%%%%%%%%%%%%%%%%%%%%%%%%%%%%%%%%%%%%%%%%%%%%%%%%%%%%%%%%
%
\subsection{Mathematical Statement of a Class of Problem Formulations}
Collecting and summarizing all of the ideas discussed so far, a class of trajectory optimization problems in the ER3BP is formulated as
\begin{equation}
\begin{aligned}
 \underset{[\bxf, \buf, \theta_0, \theta_f]}{\text{minimize}}
\quad &  J[\bxf, \buf, \theta_0, \theta_f] \\
 \text{subject to} \qquad
& \bx'(\theta) = \bff(\bx(\theta),\,\bu(\theta),\,\theta) \\
 & \be_1(\bx_0,\, \theta_0) = \bzero \\
 & \be_2(\bx_f,\, \theta_f) = \bzero \\
 & e_3(z_f, \theta_f) \leq z_f^U \\
 & 0 < \theta_f - \theta_0 \le \Delta\theta^U\\
 & \mathbf{0} \leq \bu^a \leq \bu^U\\
 & \mathbf{0} \leq \bu^b \leq \bu^U\\
\end{aligned}
\label{eq: OCP}
\end{equation}
%========================
%
A detailed explanation of \eqref{eq: OCP} is in order.
In this paper, we consider two different formulations of the cost functional $J$.  One that is given by \eqref{eq:prop-11-smooth} corresponding to a non-dimensional propellant consumption given in ``velocity units" (VU) and another given by ``minimum-time'':
\begin{equation}\label{eq: Cost_minTime}
J[\bxf, \buf, \theta_0, \theta_f] := \theta_f - \theta_0
\end{equation}
Note that $\theta_0$ is not necessarily zero; hence, the right-hand-side of \eqref{eq: Cost_minTime} cannot be replaced by $\theta_f$ only. In the case of the minimum-time problem, the $\Delta\theta^U$ bound in \eqref{eq: OCP} is removed to avoid the possibility of an empty feasible set.

For the $\Delta{prox}_1$ propellant consumption model, the EOM given in \eqref{eq: DynamicsEOM} are recast in state-space form according to
\begin{equation}
\begin{aligned}
{x}' &= v_x \\
{y}' &= v_y \\
{z}' &= v_z \\
{v}_x' &= 2v_y + \partial_x\Psi(x,y,z, \theta) + \frac{{u}_x^a - {u}_x^b}{(1+e\cos{\theta})^3} \\
{v}_y' &= -2v_x + \partial_y\Psi(x,y,z, \theta) + \frac{{u}_y^a - {u}_y^b}{(1+e\cos{\theta})^3}  \\
{v}_z' &= \partial_z\Psi(x,y,z, \theta) + \frac{{u}_z^a - {u}_z^b}{(1+e\cos{\theta})^3} \\
\end{aligned}
\label{eq: ERTBP_Controlled_StateSpace}
\end{equation}
The endpoint constraint functions $\be_1$ and $\be_2$ are given by \eqref{eq:e1==cheb} and \eqref{eq:e2==cheb} respectively.  The function $e_3$ in \eqref{eq: OCP} is simply $z(\theta_f)$ which is constrained according to \eqref{eq: zF_constr}.

\section{Development of The Necessary Conditions for Optimality}\label{sec:necessary}
The necessary conditions for optimality are first developed for the minimum-propellant case with certain key results subsequently reused for the other two cases denoted in Section~\ref{sec:OCPs}.  A fundamental reason for developing these necessary conditions is to identify those conditions that can be easily verified in a computational setting.  These verifiable conditions are boxed for easy reference.

\subsection{Minimum-Propellant Case}
The Pontryagin Hamiltonian\cite{ross_primer_2015} for the $\Delta{prox}_1$ cost functional is given by
\begin{equation}
    H(\blam, \bx, \bu, \theta) := \frac{\|\bu\|_1}{(1+e\cos{\theta})^2}+
    \boldsymbol{\lambda}^T\bff(\bx,\, \bu,\,\theta)
    \label{eq: Hamiltonian}
\end{equation}
where $\blam$ is a costate (adjoint covector) that satisfies the adjoint differential equation
\begin{equation}
    \boldsymbol{\lambda}' = -\frac{\partial{H}}{\partial{\bx}} = -\left[\frac{\partial{\bff(\bx,\,\bu,\,\theta)}}{\partial{\bx}}\right]^T \boldsymbol{\lambda}
    \label{eq: Adjoint}
\end{equation}
Expanding equation \eqref{eq: Adjoint}, the costate dynamics in scalar form are given by
\begin{equation}\label{eq:costates}
\begin{aligned}
    {\lambda}'_x&=-\partial^2_{xx}\Psi\, \lambda_{v_x} - \partial^2_{xy}\Psi\,\lambda_{v_y} - \partial^2_{xz}\Psi\,\lambda_{v_z} \\
    {\lambda}'_y&=-\partial^2_{yx}\Psi\, \lambda_{v_x} - \partial^2_{yy}\Psi\,\lambda_{v_y} - \partial^2_{yz}\Psi\,\lambda_{v_z}  \\
    {\lambda}'_z&=-\partial^2_{zx}\Psi\, \lambda_{v_x} - \partial^2_{zy}\Psi\,\lambda_{v_y} - \partial^2_{zz}\Psi\,\lambda_{v_z} \\
     {\lambda}'_{v_x} &= -\lambda_x + 2 \lambda_{v_y} \\
     {\lambda}'_{v_y} &= -\lambda_x - 2 \lambda_{v_x}\\
     {\lambda}'_{v_z} &= -\lambda_z
    \end{aligned}
\end{equation}
A quick examination of \eqref{eq:costates} shows that all of the co-position equations contain the Hessian of the pseudopotential function $\Psi$.  This Hessian is not explicitly expressed here (in terms of $x, y, z$ and $\theta$) for the purpose of brevity.  It suffices to note that the Hessian is not a constant or some other simple function that permits a closed-form integration of \eqref{eq:costates}.

Transforming the right-hand-side of \eqref{eq: Hamiltonian} in the manner described in Section~\ref{sec:OCPs}, we get a differentiable Lagrangian of the Hamiltonian\cite{ross_primer_2015},
$ \bar{H}$, such that
\begin{equation}
    \bar{H}(\bmu,\blam,\bx,\bu, \theta) = H(\blam, \bx, \bu, \theta)  + \boldsymbol{\mu}_{u^a}^T\,{\bu}^a + \boldsymbol{\mu}_{u^b}^T\,{\bu}^b
    \label{eq: LagrangianOfHamiltonian}
\end{equation}
where we have abused the $\bu$ notation in $\bar{H}$ and $H$ to imply the pair $(\bu^a, \bu^b)$.  See \eqref{eq:L1-to-smooth}.  The instantaneous path covectors $\bmu_{u^a}$ and $\bmu_{u^b}$ in \eqref{eq: LagrangianOfHamiltonian} satisfy the complementarity condition\cite{ross_primer_2015}
\begin{equation}
\boxed{
{\mu}_\xi^{a,b}
\begin{cases}
\begin{aligned}
&  \leq 0 \qquad\text{if }\qquad  u_\xi^{a,b} = 0 \\
&  = 0 \qquad\text{if }\qquad 0 < u_\xi^{a,b}< u_\xi^{U} \\
&  \geq 0 \qquad\text{if }\qquad u_\xi^{a,b} = u_\xi^{U}  \\
\end{aligned}
\end{cases}
} \quad\forall\ \xi \in \{x, y, z\}
\label{eq: Coplementarity}
\end{equation}
Equation \eqref{eq: Coplementarity} is boxed because it can be computationally verified. It dictates the switching structure for the optimal control\cite{ross_primer_2015}, and can therefore be used to verify the timing and magnitude of the maximum thrust arcs and singular arcs, if any\cite{singular-arc-2010-CDC}. The optimal control satisfies the stationarity conditions\cite{ross_primer_2015}
\begin{equation}
    \begin{aligned}
        \frac{\partial{\bar{H}}}{\partial{\bu^a}} &= \bzero \\
        \frac{\partial{\bar{H}}}{\partial{\bu^b}} &= \bzero
    \end{aligned}
    \label{eq: Stationarity}
\end{equation}
Let $\bar{E}$ be the Endpoint Lagrangian\cite{ross_primer_2015} given by
\begin{multline}
    \bar{E}(\bnu,\,\bx(\theta_0),\,\bx(\theta_f),\, \theta_0,\,\theta_f) = \boldsymbol{\nu}_0^T\,\be_1(\bx(\theta_0),\, \theta_0)  \\
    + \boldsymbol{\nu}_f^T\,\be_2(\bx(\theta_f),\, \theta_f) + \nu_{z_f}\,e_3(z(\theta_f),\,\theta_f)
    \label{eq: EnpointLagrangian_General}
\end{multline}
where $\boldsymbol{\nu}_0$ and $\boldsymbol{\nu}_f$ are the endpoint covectors.  The transversality conditions are then given by
\begin{equation}\label{eq:Costate_FinalValues}
\boxed{
\begin{aligned}
\boldsymbol{\lambda}(\theta_0) &= -\frac{\partial{\bar{E}}}{\partial{\bx(\theta_0)}} = -  \left[ \frac{\partial\,\be_1(\bx(\theta_0),\, \theta_0)}{\partial\bx(\theta_0)} \right]^T \boldsymbol{\nu}_0 = -\boldsymbol{\nu}_0
\\
\boldsymbol{\lambda}(\theta_f) & = \frac{\partial{\bar{E}}}{\partial{\bx(\theta_f)}} =   \left[ \frac{\partial\,\be_2(\bx(\theta_f),\, \theta_f)}{\partial\bx(\theta_f)} \right]^T \boldsymbol{\nu}_f \\
&\qquad\qquad\qquad  + \nu_{z_f} \left[ \frac{\partial\,{e}_3(z(\theta_f),\, \theta_f)}{\partial\bx(\theta_f)} \right] \\
&= \boldsymbol{\nu}_f + [0 \; 0 \: \nu_{z_f} \: 0 \: 0 \: 0]^T
\end{aligned}
}
\end{equation}
where $\nu_{z_f}$ satisfies the complementarity condition given by
\begin{equation}
\nu_{z_f}
\begin{cases}
\begin{aligned}
&  = 0 \quad  z(\theta_f) < z_f^{\text{U}}  \\
&  \geq 0 \quad z(\theta_f) = z_f^{\text{U}}  \\
\end{aligned}
\end{cases}
\label{eq: KKT_zf}
\end{equation}
Similarly, the Hamiltonian value condition\cite{ross_primer_2015} is given by
\begin{equation}
\boxed{
\begin{aligned}
        \mathcal{H}[@\theta_0] = \frac{\partial{\bar{E}}}{\partial{\theta_0}} &= \boldsymbol{\nu}_0^T\,\left[ \frac{\partial\,\be_1(\bx(\theta_0),\, \theta_0)}{\partial\theta_0} \right]\\
        &= -\boldsymbol{\nu}_0^T\,\left[ \frac{d\,\bcheb_0(\theta_0)}{d\theta_0} \right]
        \\[1em]
        \mathcal{H}[@\theta_f] = -\frac{\partial{\bar{E}}}{\partial{\theta_f}} &= - \boldsymbol{\nu}_f^T\left[\frac{\partial\,\be_2(\bx(\theta_f),\, \theta_f)}{\partial\theta_f}\right] \\
            &= \boldsymbol{\nu}_f^T\, \left[\frac{ d\,\bcheb_f(\beta(\theta_f))}{d\theta_f}\right]
\end{aligned}
}
    \label{eq:HVC_InitAndFinal}
\end{equation}
where $\mathcal{H}$ is the minimized Hamiltonian\cite{ross_primer_2015}.  The evaluation of the derivative of the Chebyshev function in \eqref{eq:HVC_InitAndFinal} at the initial and final true anomalies is easily computed because the derivative of a Chebyshev polynomial can be expressed in terms of Chebyshev polynomials\cite{ross_universal_2023}.  This equivalence relationship is used in Chebfun\cite{t_a_driscoll_chebfun_2014} to produce the derivatives indicated in \eqref{eq:HVC_InitAndFinal}.  These transversality conditions are boxed because the final values of the covectors are an output of the optimal control software, DIDO \cite{ross_enhancements_2020}.

As noted before, the dynamics of the ER3BP are non-autonomous.  As a result, the Hamiltonian evolution equation\cite{ross_primer_2015},
\begin{equation}
    \frac{d\mathcal{H}}{d\theta} =\frac{\partial{H}}{\partial{\theta}},
    \label{eq: HEE}
\end{equation}
does not provide an easily verifiable integral of motion.

\subsection{Time-Bounded Case}
For the time-bounded case, the necessary conditions follow the same derivation as given in the minimum-propellant case with the exception of transversality. The Endpoint Lagrangian is now given by
\begin{multline}
    \bar{E}(\bnu,\,\bx(\theta_0),\,\bx(\theta_f),\, \theta_0,\,\theta_f) = \boldsymbol{\nu}_0^T\,\be_1(\bx(\theta_0),\, \theta_0)  \\
    + \boldsymbol{\nu}_f^T\,\be_2(\bx(\theta_f),\, \theta_f) + \nu_{z_f}\,e_3(z(\theta_f),\,\theta_f) + \nu_{\theta}\,(\theta_f -\theta_0)
    \label{eq: EnpointLagrangian_timeBound}
\end{multline}
Note that ``time-bounded" means the time of flight (TOF) is restricted to an arbitrary duration instead of simply restricting the final true anomaly by itself.  This is a useful formulation since the transfer depends on the configuration of the primaries at the initiation point.  By taking the boundedness to mean TOF, the configurations of the primaries other than $\theta_0=0$ are admissible.  The initial and final values of the costates remain the same as given in equation \eqref{eq:Costate_FinalValues}, as does the evolution of the minimized Hamiltonian in equation \eqref{eq: HEE}.  The initial and final value of the minimized Hamiltonian, however, is now given by
\begin{equation}
\boxed{
\begin{aligned}
        \mathcal{H}[@\theta_0] &= \frac{\partial{\bar{E}}}{\partial{\theta_0}} =  -\boldsymbol{\nu}_0^T\, \left[ \frac{d\,\bcheb_0(\theta_0)}{d\theta_0} \right] - \nu_{\theta}
        \\
        \mathcal{H}[@\theta_f] &= -\frac{\partial{\bar{E}}}{\partial{\theta_f}} = \boldsymbol{\nu}_f^T \left[\frac{ d\,\bcheb_f(\beta(\theta_f))}{d\theta_f}\right] - \nu_{\theta}
\end{aligned}
}
    \label{eq:HVC_final_timeBound}
\end{equation}
where $\nu_\theta$ satisfies the complementarity condition given by
\begin{equation}
\nu_{\theta}
\begin{cases}
\begin{aligned}
&  \leq 0 \qquad \theta_f - \theta_0 = 0 \\
&  = 0 \quad 0 < \theta_f - \theta_0 < \Delta\theta^{U} \\
&  \geq 0 \qquad \theta_f - \theta_0 = \Delta\theta^{U}  \\
\end{aligned}
\end{cases}
\label{eq: KKT_theta}
\end{equation}
From \eqref{eq: KKT_theta}, it follows that we expect $\nu_\theta \nless 0$.  Alternatively, we expect $\nu_\theta \ge 0$.

\subsection{Minimum-Time Case}
For the minimum-time case, the development of the necessary conditions are carried out in the same manor as the minimum-propellant derivation except that the cost functional is given by \eqref{eq: Cost_minTime}.
As in the time-bounded case, the minimization of the TOF permits additional initial primary configurations other than $\theta_0 = 0$.  The Hamiltonian is now formulated as
\begin{equation}
    H(\blam, \bx, \bu, \theta) := \boldsymbol{\lambda}^T\bff(\bx,\, \bu,\,\theta)
    \label{eq: Hamiltonian_minTime}
\end{equation}
The optimal control is still governed by equations \eqref{eq: Stationarity} and \eqref{eq: Coplementarity}. Similarly, the costate dynamics given by equation \eqref{eq: Adjoint} remain valid.  The Endpoint Lagrangian, however, now becomes
\begin{multline}
     \bar{E}(\bnu,\,\bx(\theta_0),\,\bx(\theta_f),\, \theta_0,\,\theta_f) = \theta_f-\theta_0 + \boldsymbol{\nu}_0^T\,\be_1(\bx(\theta_0),\, \theta_0)  \\
     + \boldsymbol{\nu}_f^T\,\be_2(\bx(\theta_f),\, \theta_f) + \nu_{z_f}\,e_3(z(\theta_f),\,\theta_f)
    \label{eq: EnpointLagrangian_minTime}
\end{multline}
Though the new Endpoint Lagrangian includes the endpoint cost\cite{ross_primer_2015}, the initial and final values of the costates are still given by equation \eqref{eq:Costate_FinalValues}.  However, now the initial and final value of the minimized Hamiltonian are determined by
\begin{equation}
\boxed{
\begin{aligned}
        \mathcal{H}[@\theta_0] &= \frac{\partial{\bar{E}}}{\partial{\theta_0}} =  - \boldsymbol{\nu}_0^T\, \left[ \frac{d\,\bcheb_0(\theta_0)}{d\theta_0} \right] - 1
        \\
        \mathcal{H}[@\theta_f] &= -\frac{\partial{\bar{E}}}{\partial{\theta_f}} = \boldsymbol{\nu}_f^T \left[\frac{ d\,\bcheb_f(\beta(\theta_f))}{d\theta_f}\right] - 1
\end{aligned}
}
    \label{eq:HVC_final_minTime}
\end{equation}
As before, the partial derivatives of the endpoint functions are determined by computing the derivatives of the Chebyshev interpolants using Chebfun\cite{t_a_driscoll_chebfun_2014}.

\section{An Overveiw of the Universal Birkhoff Theory for Trajectory Optimization}\label{sec:BirkIntro}
A comprehensive description of the universal Birkhoff theory for trajectory optimization is presented in [\citen{ross_universal_2023-1}].  As noted in Section~\ref{sec:intro}, the theory can be applied to both ``direct'' and ``indirect'' methods for trajectory optimization\cite{ross_universal_2023-1,ross_universal_2023}; however, it is best used in conjunction with the fast spectral algorithm\cite{gong_spectral_2008,proulx_implementations_2023,ross_enhancements_2020,ross_pseudospectral_2004, auto-knots}.  In this section we briefly review this theory and its computational principles.

\subsection{Introduction to the Universal Birkhoff Interpolant}
Let $\pi^N :=\set{\theta_1, \ldots, \theta_N}$ be an arbitrary grid such that \footnote{The notation used in this section is different from [\citen{ross_universal_2023-1}]
and [\citen{ross_universal_2023}] but is consistent with the other sections of this paper.}
$$\theta_0 \le \theta_1 < \cdots < \theta_N \le \theta_f  $$
The universal Birkhoff theory requires that the grid $\pi^N$ be selected such that the state variable $\bx(\cdot): \theta \mapsto \real{6}$ be approximated by two equivalent $a$- and $b$-expansions given by \cite{ross_universal_2023-1, ross_universal_2023}
\begin{equation}\label{eq:birk-a}
\bx^N(\theta) =\bx_0 B_0^0(\theta) + \sum_{j=1}^N \bv_j B_j^a(\theta) = \sum_{j=1}^N \bv_j B_j^b(\theta) + \bx^b B_N^N(\theta)
\end{equation}
where, $\bx_0 = \bx(\theta_0)$, $\bx_f = \bx(\theta_f) $ and $\bv_j \in \real{N_x}, j = 1, \ldots, N$ are the unknown optimization variables, $B_j^a(\theta), B_j^b(\theta), j = 1, \ldots, N$ and $B^0_0(\theta)$ are the Birkhoff basis functions that satisfy the interpolation conditions
\begin{align}\label{eq:birk-conditions}
\begin{aligned}
B_0^0(\theta_0)           &= 1,            & {\dot B}_0^0(\theta_i)   &=  0,      &&   i = 1, \ldots, N \\
B_j^a(\theta_0)           &=  0,           &{\dot B}_j^a(\theta_i) &=  \delta_{ij},  &&   i = 1, \ldots, N, & j = 1, \ldots, N \\
B_N^N(\theta_f)           &= 1,            & {\dot B}_N^N(\theta_i)   &=  0,      &&   i = 1, \ldots, N \\
B_j^b(\theta_f)           &=  0,           &{\dot B}_j^b(\theta_i) &=  \delta_{ij},  &&   i = 1, \ldots, N, & j = 1, \ldots, N \, ,
\end{aligned}
\end{align}
and $\delta_{ij}$ is the Kronecker delta.   The $a$- and $b$-versions of the Birkhoff matrices are given by
\begin{equation}\label{eq:Ba-mat-def}
\bB^a := [B^a_{ij}] := B^a_j(\theta_i) \quad \bB^b := [B^b_{ij}] := B^b_j(\theta_i)
\end{equation}
Explicit formulas for the computation of the Birkhoff matrices are given in \cite{ross_universal_2023}. See also \cite{sandia_millionPt_2025} for an FFT-based fast computation of Birkhoff matrix-vector products.

\subsection{Discretization of the State and Adjoint Equations}

Let $\bX \in \real{N_x \times N}$ and $\bV \in \real{N_x \times N}$ be matrices whose columns are given by $[\bx^N(\theta_1), \ldots, \bx^N(\theta_N) ]$ and $[\bv_1, \ldots, \bv_N) ]$ respectively. Then, from \eqref{eq:birk-a} and \eqref{eq:Ba-mat-def} it follows that
\begin{equation}\label{eq:X=BV}
\bX^T = \bB^a \bV^T + \bb\, \bx_0^T
\end{equation}
where $\bb := [1, \ldots, 1]^T \in \real{N}$ is a vector of ones. In its pseudospectral version, the matrix $\bV$ satisfies the dynamics function according to
\begin{equation}
\bV = \bff(\bX, \bU, \pi^N)
\end{equation}
where $\bU:= [\bu_1, \ldots, \bu_N ] \in \real{N_x \times N}$ and $\bff$ is overloaded according to $\bff(\bx_i, \bu_i, \theta_i) = \bv_i, \ i = 1, \ldots, N$.  The endpoint conditions are satisfied according to
\begin{equation}
\be^L \le \be(\bx_0, \theta_0, \bx_f, \theta_f) \le \be^U
\end{equation}
where $\bx_f$ satisfies a linear constraint given by
\begin{equation}\label{eq:x0-xf}
\bx_0 - \bx_f  + \bV \bw_B= \bzero
\end{equation}
and $\bw_B := [w_1, \ldots, w_N]^T$ are the Birkhoff weights\cite{ross_universal_2023} associated with the grid $\pi^N$.

\begin{figure}[h!]
\centering
\includegraphics[width=0.7\columnwidth]{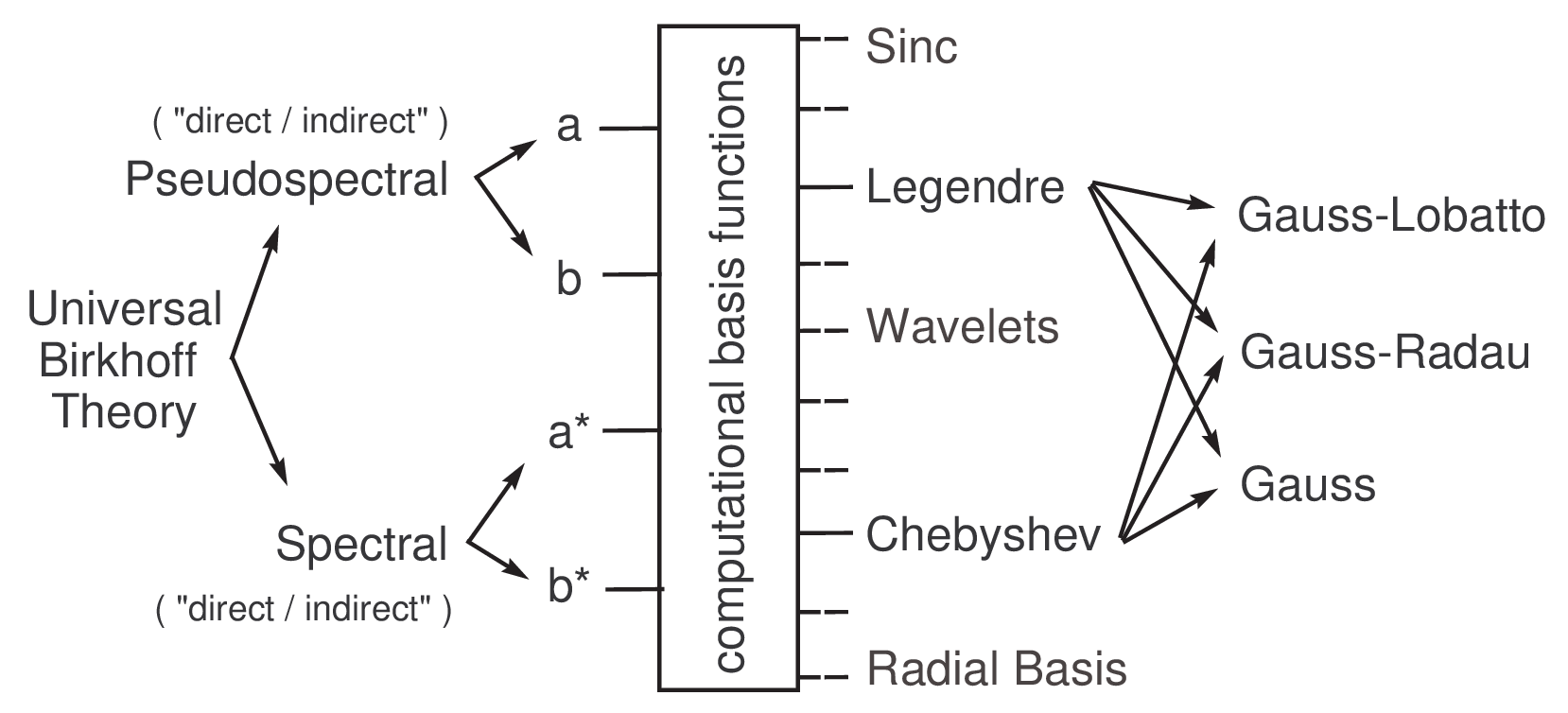}
\caption{A schematic for the flow a universal Birkhoff theory from inception to computation.  Figure adapted from \cite{ross_enhancements_2020}.}
\label{fig:BirkhoffTypes}
\end{figure}

Next, consider the approximation of the adjoint variable $\blam(\cdot): \theta \mapsto \real{N_x}$.  This equation can be written in exactly the same form as \eqref{eq:X=BV}\cite{ross_universal_2023-1}.  Nonetheless, to illustrate the interplay between the $a$- and $b$-versions, we consider the approximation of the adjoint variable by the $b$-version of the universal Birkhoff interpolant,
\begin{equation}\label{eq:birk-b}
\blam^N(\theta) :=\sum_{j=1}^N \bomega_j B_j^b(\theta)  + \blam_f B_N^N(\theta)
\end{equation}
Let $\bLam \in \real{N_x \times N}$ and $\bOmega \in \real{N_x \times N}$ be matrices whose columns are given by $[\blam^N(\theta_1), \ldots, \blam^N(\theta_N) ]$ and $[\bomega_1, \ldots, \bomega_N) ]$ respectively. Then, similar to \eqref{eq:x0-xf}, $\blam_f$ satisfies a linear constraint given by
\begin{equation}
\blam_0 - \blam_f  + \bw_B^T \bOmega^T = \bzero
\end{equation}
where $\blam_0$ and $\blam_f$ satisfy the transversality conditions
\begin{equation}
\blam_0 = -\partial_{\bx_0}\overline{E}(\bnu, \bx_0, \theta_0, \bx_f, \theta_f), \quad \blam_f = \partial_{\bx_f}\overline{E}(\bnu, \bx_0, \theta_0, \bx_f, \theta_f)
\end{equation}
In its spectral version\cite{ross_universal_2023-1}, $\bOmega$ satisfies the spectral analog of \eqref{eq:X=BV}:
\begin{equation}\label{eq:Lam=BOmega}
\bW^T\circ\bLam^T = \bW^T\circ\bB^b \bOmega^T + \bW^T\circ\bb \blam_f^T
\end{equation}
where $\circ$ denotes a Hadamard product and $\bW \in \real{N_x \times N}$ matrix with each row identically equal to $\bw_B^T$. Similarly, $\bOmega$ is chosen to satisfy the adjoint equation according to
$$\bW\circ \bOmega + \bW\circ \partial_{\bx} H(\bLam, \bX, \bU, \pi^N) = \bzero $$
where $\partial_{\bx} H $ is overloaded according to $\partial_{\bx} H(\blam_i,  \bx_i, \bu_i, t_i) = \bomega_i, \ i = 1, \ldots, N$.
\begin{remark}
An important point to note in all of the preceding equations is that there was no mention of a polynomial basis function. In fact, the universal Birkhoff theory does not require a polynomial basis function\cite{ross_universal_2023-1}.
\end{remark}
Now suppose we expand the Birkhoff basis functions in terms of functions that have known computational advantages.  If such a representation of a Birkhoff basis function satisfies the two hypotheses enunciated in \cite{ross_universal_2023-1}, then it constitutes a candidate Birkhoff method for trajectory optimization.   A flow of these ideas is depicted in Fig.~\ref{fig:BirkhoffTypes}.
%
%%%%%%%%%%%%%%%%%%%%%%%%%%%%%%%%%%%%%%%%%%%%%%%%%%%%%%%%%%%%%%%%%%%%%%%%%%%%%
%%%%%%%%%%%%%%%%%%%%%%%%%%%%%%%%%%%%%%%%%%%%%%%%%%%%%%%%%%%%%%%%%%%%%%%%%%%%%
%
As implied in Fig.~\ref{fig:BirkhoffTypes}, the Birkhoff theory while not relying on polynomials does not exclude it.  It was shown in [\citen{proulx_implementations_2023}] that a family of Gegenbauer polynomials (and their grids) satisfy all of the hypotheses needed for a Birkhoff theory to hold.
%----
\begin{remark}
Even when the basis functions are polynomials, the control function $\bu(\cdot)$ is not necessarily expressed in terms of polynomial expansion.  Moreover, it is frequently and naturally expressed in terms of a non-polynomial basis function.  This is why there is no Gibbs phenomenon if the controls are discontinuous.  See \cite{ross_universal_2023-1} and \cite{proulx_implementations_2023} for complete details. 
\end{remark}
%==============

\subsection{Spectral Algorithm for the Universal Birkhoff Method}
The original spectral algorithm\cite{gong_spectral_2008, ross_pseudospectral_2004,ross_review_2012,auto-knots} can be easily adapted to the Birkhoff-discretized optimal control problem. In fact, a Birkhoff discretization allows for a faster, simpler version of the spectral algorithm. This is because the condition number of a Birkhoff-discretized problem does not change as $N$ tends to infinity\cite{ross_universal_2023-1}. Hence, there is no need to subdivide the time interval into ``$h$'' elements and control the order ``$p$'' of the local polynomial just to maintain a low condition number. As a result, error analysis and error control\cite{ross_pseudospectral_2004,auto-knots} are substantially simplified to the core spectral algorithm\cite{gong_spectral_2008,auto-knots} together with the enhancements described in \cite{ross_enhancements_2020}.
In simple terms, the core spectral algorithm can be described as follows. First, note that $\bX$ and $\bLam$ in \eqref{eq:X=BV} and \eqref{eq:Lam=BOmega} are dependent upon $N$.  Hence, to describe the spectral algorithm, we relabel $\bX$ and $\bLam$ as $\bX^N$ and $\bLam^N$ respectively.  Second, for any fixed value of $N$, consider solving \eqref{eq:X=BV}, \eqref{eq:Lam=BOmega} and other related equations described in the previous subsection. This solution process is iterative.  For a fixed value of $N$, these iterations can be described as generating a sequence $(\bX^N_m, \bLam^N_m)_{m=0}^\infty $.  These simple ideas constitute the core ideas behind a spectral algorithm. Succinctly, the spectral algorithm produces a double-infinite sequence\cite{ross_enhancements_2020,gong_spectral_2008,ross:guess-free} of vectors
\begin{multline*}
  (\bX, \bLam)_{m_0}^{N_0},  (\bX, \bLam)_{m_1}^{N_1},\ldots, (\bX, \bLam)_{m_f}^{N_f}, \quad N_{k+1} > N_k, \\
  k = 0, 1, \ldots, N_f
\end{multline*}
such that for some $m_f$ and $N_f$, $(\bX, \bLam)_{m_f}^{N_f}$ is a solution within some tolerance criterion.  The production of this double-infinite sequence is built on three main algorithmic components illustrated in Fig.~\ref{fig:DIDOMain3}:
%
%%%%%%%%%%%%%%%%%%%%%%%%%%%%%%%%%%%%%%%%%%%%%%%%%%%%%%%%%%%%%%%%%%%%%%%%%%%%%
\begin{figure}[h!]
\centering
\includegraphics[width=0.7\columnwidth]{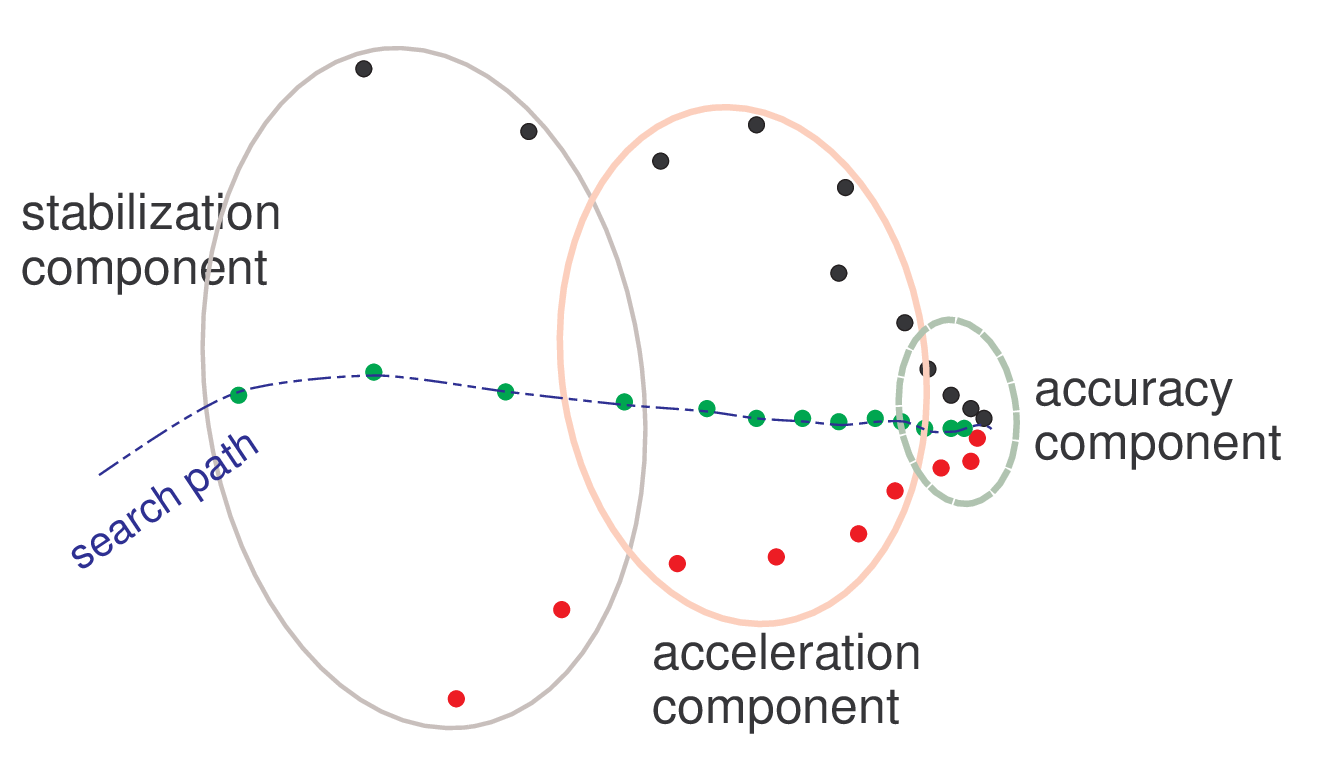}
\caption{Three main algorithmic components of DIDO.  Figure from \cite{ross_enhancements_2020}.}
\label{fig:DIDOMain3}
\end{figure}
%%%%%%%%%%%%%%%%%%%%%%%%%%%%%%%%%%%%%%%%%%%%%%%%%%%%%%%%%%%%%%%%%%%%%%%%%%%%%
%
\begin{enumerate}
\item The first sub-algorithm is the stabilization component.  This component is largely the guess-free algorithm described in \cite{ross:guess-free}. Its task is to start from an ``arbitrary point'' and generate a hand-off solution to the acceleration component\cite{ross_enhancements_2020}.
\item The main task of the acceleration component is to produce a sequence of numbers $N_1, N_2, \ldots$ that are coordinated with the spectral convergence rate\cite{gong_spectral_2008} and error tolerances $\epsilon^{N_k} > 0$ to produce an approximate but fast solution to the Birkhoff-discretized optimal control problem~\cite{proulx_implementations_2023}.
\item This approximate solution is handed off to the accuracy component to drive the solution to specified error tolerances\cite{ross_review_2012,ross_enhancements_2020} $\epsilon \ge \epsilon^{N_f} > \sqrt{\epsilon_M}$, where $\epsilon_M$ is machine precision.
\end{enumerate}
A search Lyapunov function\cite{rossJCAM-1,ross-CD,ross_enhancements_2020,ross-JNVA} coordinates these sequences so that the iterations transition ``desirably'' as illustrated by the search path shown in Fig.~\ref{fig:DIDOMain3}.  In general, backtracking in function space may be required to promote convergence\cite{ross_enhancements_2020}.

The preceding ideas adapted for the Birkhoff-Gegenbauer method\cite{proulx_implementations_2023} is implemented in an $\alpha$-version of DIDO\cite{ross_enhancements_2020}.  It is apparent from the preceding description that the spectral algorithm is neither a direct nor an indirect method\cite{conway_survey_2012,caillau_algorithmic_2023}. That is, although the spectral algorithm utilizes the costates, the Hamiltonian and the various other dual variables generated from Pontryagin's Principle, all of this information is autonomously generated\cite{ross_enhancements_2020} using mechanized formulas presented in \cite{ross_primer_2015}.  Such an enhanced $\alpha$-version of DIDO was used in all the computations reported in Section~\ref{sec:numerics}.

\subsection{Grid Selection}
The spectral algorithm is agnostic to the particulars of a grid. The main requirement for a successful implementation of the spectral algorithm is that the error must theoretically vanish as $N$ tends to infinity.
Figure~\ref{fig:BirkhoffTypes} implies that six Gegenbauer grid points meet this criteiron\cite{proulx_implementations_2023}.  These grid points are the Legendre/Chebyshev-Gauss-Lobatto (LGL/CGL), Legendre/Chebyshev-Gauss-Radau (LGR/CGR) and the Legendre/~Chebyshev-Gauss (LG/CG) nodes.  All of these grid points satisfy the hypotheses required for the Birkhoff theory to hold\cite{ross_universal_2023-1,proulx_implementations_2023}.  As a result, it is reasonable to ask the question if any one of these grid points is computationally superior to the other.

There is a widespread misconception within the engineering community that the Legendre grids are superior to Chebyshev and that the Gauss grid, in particular, offers the highest accuracy\cite{trefethen-myths-2011}. This misconception is based on misreading theorems on Gauss quadrature\cite{trefethen-myths-2011}.  In trajectory optimization, such misconceptions run deeper\cite{ross_universal_2023-1}.  The origin of the issue in trajectory optimization is conflating the difference between Gauss quadrature over known versus unknown points.  That is, textbook theorems on Gauss quadrature are over known points whereas trajectory optimization is the determination of unknown points. See \cite{ross_review_2012} and \cite{ross_universal_2023-1} for details. In any case, to prevent an occurrence of new misconceptions in Birkhoff-theoretic methods, we consider solving a sample optimal problem posed in \eqref{eq: OCP} using all six grid points.  Before discussing the numerical results, we briefly note that all six grid points inherently posses a spectral convergence rate; i.e., a near-exponential rate of convergence\cite{ross_universal_2023-1,proulx_implementations_2023}.  
\subsubsection{Primal Variables}

A comparison of the numerical study for an illustrative sample of primal variables is shown in Table~\ref{tab:GridComparePrimal}.
%
%%%%%%%%%%%%%%%%%%%%%%%%%%%%%%%%%%%%%%%%%%%%%%%%%%%%%%%%%%%%%%%%%%%%%%%%%%%%%%
\begin{table*}[h!]
	\fontsize{10}{10}\selectfont
    \caption{Sample values of primal variables for different choices of Birkhoff grids }
    \label{tab:GridComparePrimal}
    \centering
    \begin{tabular}{l  c  c  c  c c c} % Column formatting
      \hline\hline
       Parameter &LGL & LGR & LG & CGL & CGR & CG \\
      \hline
      cost   & $0.6096$   & $0.6096$  & $0.6096$ & $0.6096$ & $0.6096$ & $0.6096$\\
      $\theta_0$        & $0.0000$   & $0.0000$ & $0.0000$ & $0.0000$ & $0.0000$ & $0.0000$  \\
      $\theta_f$      & $7.0847$   & $7.0847$ & $7.0847$ & $7.0847$ & $7.0847$ & $7.0847$ \\
      \hline\hline
    \end{tabular}
\end{table*}
%%%%%%%%%%%%%%%%%%%%%%%%%%%%%%%%%%%%%%%%%%%%%%%%%%%%%%%%%%%%%%%%%%%%%%%%%%%%%%
%
It is abundantly clear from the value of the cost function and the sample values of the primal variables presented in Table~\ref{tab:GridComparePrimal} that there is no difference in relative accuracy across all of the six grid points up to four decimal places.  As noted elsewhere\cite{ross_enhancements_2020,ross_universal_2023-1} any comparison of number past six significant digits must be viewed with great suspicion (for any method) due to the fundamental connection between the condition number of the Hessian and the square root of machine precision.  Furthermore, in noting that a number past four significant digits is a comparison of one part in $10,000$, any claim of higher accuracy past four significant digits is a claim of minutia.  With this perspective in mind, it is apparent that the only conclusion that can be drawn from Table~\ref{tab:GridComparePrimal} is that all six grid points provide identical accuracy, at least as far as primal variables are concerned.  This numerical study is consistent with prior analyses\cite{proulx_implementations_2023,koeppen_fast_2019}.  We briefly note that the practical accuracy of primal variables in aerospace engineering must be paired with sensor accuracy and the corresponding granularity of error that can be achieved by a control action. 

\subsubsection{Dual Variables}
Consider a similar numerical comparison in dual space. For brevity, only the values of the first three components of the initial multiplier $\bnu_0$ are shown in Table~\ref{tab:GridCompareDual}.
%
%%%%%%%%%%%%%%%%%%%%%%%%%%%%%%%%%%%%%%%%%%%%%%%%%%%%%%%%%%%%%%%%%%%%%%%%%%%%%%
\begin{table*}[h!]
	\fontsize{10}{10}\selectfont
    \caption{Sample values of $\bnu_0$ for different choices of Birkhoff grids }
    \label{tab:GridCompareDual}
    \centering
    \begin{tabular}{l  c  c  c  c c c} % Column formatting
      \hline\hline
       $\bnu_0$ &LGL & LGR & LG & CGL & CGR & CG \\
      \hline
      $1$   & $-0.3691$   & $-0.3689$  & $-0.3697$ & $-0.3689$ & $-0.3687$ & $-0.3689$\\
      $2$        & $0.0639$   & $0.0638$ & $0.0640$ & $0.0637$ & $0.0635$ & $0.0634$  \\
      $3$      & $0.2390$   & $0.2388$ & $0.2389$ & $0.2391$ & $0.2390$ & $0.2397$ \\
      \hline\hline
    \end{tabular}
\end{table*}
%%%%%%%%%%%%%%%%%%%%%%%%%%%%%%%%%%%%%%%%%%%%%%%%%%%%%%%%%%%%%%%%%%%%%%%%%%%%%%
%
Clearly, there are small differences in the values of the components of $\bnu_0$.  First, note that all these differences are in a few parts per thousand. For instance, in the second row of Table~\ref{tab:GridCompareDual}, the values of the components of $\bnu_0$ are $0.0639$ and $0.0635$ for the LGL and CGR grids respectively.  The difference between these numbers is $0.0004$. Second, what is not shown in Table~\ref{tab:GridCompareDual} is the corresponding value of $\blam_0$. The corresponding components of $\blam_0$ for the LGL and CGR grids were computed to be $0.0639$ and $0.0636$ respectively.  Thus, the difference between the components of $\blam_0$ and $\bnu_0$ within the CGR grid differ by $0.0001$ even though the difference across the LGL and CGR grid points is larger at $0.0004$. That is, even though a particular component of $\bnu_0$ might be slightly different across the grid points, the self-consistency between $\blam_0$ and $\bnu_0$ is greater within a grid point selection than their consistency across the grid points. Note also that because the covectors are not necessarily unique\cite{ross_primer_2015}, self consistency is more valuable than any perceived inconsistency across different grid points. Regardless, we observe that these differences are numerically negligible.  

Note that the preceding tests were performed for the current suite of numerically sensitive problems under consideration and not some other sample problem where an ``exact'' solution is known.  The issue in comparative studies with toy problems is that it can easily obfuscate certain numerical characteristics that may be critical to practical problems where closed-form solutions are unavailable.  Finally, we note that none of these numerical analyses are of any value without a sound theoretical mathematical backing.

\subsubsection{The Case for Chebyshev Grids}
The main takeaway from all of these numerical comparisons is to emphasize that the computations (when performed correctly) are consistent with the mathematical predictions that there are virtually no differences in accuracy across the six Gegenbauer grid points.  We caveat the preceding statement that there might be corner cases that amplify the differences\cite{ross_review_2012}.  Given the totality of all these studies and analyses, a reasonable question to ask at this juncture is if there are other criteria for a selection of grid points.  As noted previously\cite{proulx_implementations_2023,fahroo:cheb-jgcd,cheb-costate}, of all of the six grid points, the Chebyshev nodes are the easiest to implement because they can be computed simply and ``instantaneously'' using a single line of code (namely, using a cosine function).  In contrast, a computation of the Legendre grids require sophisticated algorithms particularly for large values of $N$.  Finally, we note that a Birkhoff-Chebyshev matrix-vector product can be computed at $O(N\log(N))$ computational speed using an FFT\cite{sandia_millionPt_2025}.  For all of these reasons we choose the Chebyshev version of the Birkhoff theory for all computations presented in Section~\ref{sec:numerics}.

\section{Sample Numerical Results and Analysis}\label{sec:numerics}
It is apparent from the discussions of the preceding sections that the universal Birkhoff theory is theoretically capable of providing verifiable extremal solutions to a large class of trajectory optimization problems in the ER3BP. To demonstrate these ideas we consider a variety of non-coplanar orbit transfer problems starting from an ${L_1}$ Lyapunov orbit and ending in a southern NRHO as depicted in Fig.~\ref{fig:ER3BP_BoundaryOrbits}. For all cases considered, we use the $\Delta{prox}_1$ model for propellant consumption given by \eqref{eq:L1-l1} with an upper control bound of $\bu^U=(0.1,0.1,0.1)$.  

In this section, we also perform an independent validation of the numerical results.  To assist such numerical analyses, we present results up to six decimal places caveated by the discussions of the preceding section. We first consider the minimum-time solution as a matter of standard practice\cite{ross_primer_2015}.  This is because the minimum-time problem establishes the lower bound on transfer time for all other cases.

\subsection{Minimum-Time Solution}
\begin{figure}[h!]
 \centering
 \includegraphics[width=0.7\columnwidth]{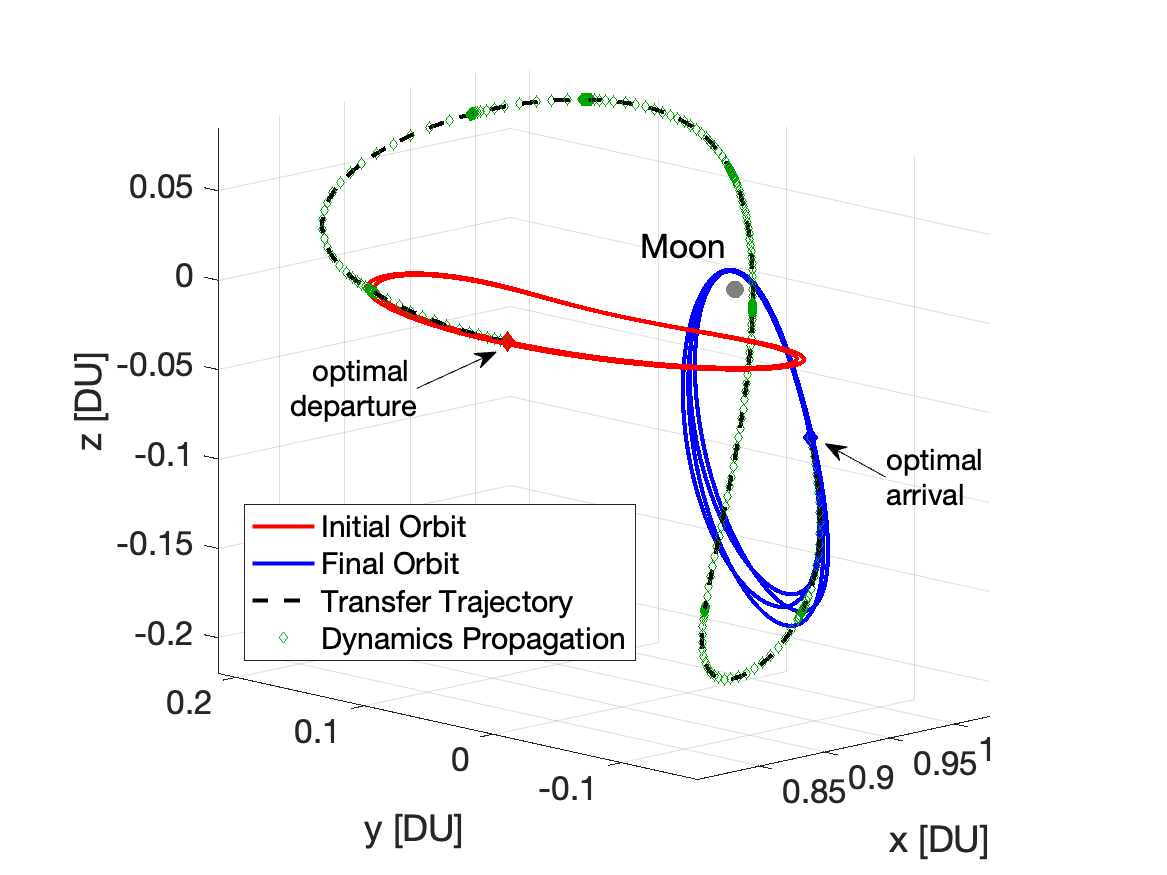}
 \caption{Minimum-time transfer between an $L_1$ Lyapunov orbit and $L_2$ southern NRHO}
 \label{fig:L1_Lyap_to_SNRHO_minTime_Trajectory}
 \end{figure}

Fig.~\ref{fig:L1_Lyap_to_SNRHO_minTime_Trajectory} shows a candidate minimum-time transfer from the 2:1 resonant $L_1$ Lyapunov orbit to the 4:1 resonant southern NRHO (see also Fig.~\ref{fig:ER3BP_BoundaryOrbits}).  The Birkhoff-theoretic spectral algorithm (implemented in an $\alpha$-version of DIDO\cite{ross_enhancements_2020}) converged for $N = 1000$.
%
%%%%%%%%%%%%%%%%%%%%%%%%%%%%%%%%%%%%%%%%%%%%%%%%%%%%%%%%%%%%%%%%%%%%%%%%%%%%%%
%%%%%%%%%%%%%%%%%%%%%%%%%%%%%%%%%%%%%%%%%%%%%%%%%%%%%%%%%%%%%%%%%%%%%%%%%%%%%
%
Also shown in Fig.~\ref{fig:L1_Lyap_to_SNRHO_minTime_Trajectory} are the locations of the candidate optimal departure and arrival points.  We first discuss the extremality of these points via the transversality conditions derived in Section~\ref{sec:necessary}.

\subsubsection{Extremality of the Departure and Arrival Points}

The candidate optimal departure and arrival points were found (by $\alpha$-DIDO) to be the values given in Table~\ref{tab: optimalStateTheta_minTime}.
%
%%%%%%%%%%%%%%%%%%%%%%%%%%%%%%%%%%%%%%%%%%%%%%%%%%%%%%%%%%%%%%%%%%%%%%%%%%%%%%
\begin{table}[h!]
	\fontsize{10}{10}\selectfont
    \caption{Candidate optimal true anomaly and state vector for departure and arrival points for minimum-time transfer }
    \label{tab: optimalStateTheta_minTime}
    \centering
    \begin{tabular}{l  c  c } % Column formatting
      \hline\hline
       coordinate & optimal departure $(x_0^*)$ & optimal arrival $(x_f^*)$ \\
      \hline
      $\theta$   & $3.120787$   & $6.961777$  \\
      $x$        & $0.807280$   & $0.999370$  \\
      $y$        & $-0.006668$  & $-0.047280$  \\
      $z$        & $0.000000$   & $-0.077317$  \\
      $v_x$      & $-0.006912$  & $-0.073385$  \\
      $v_y$      & $0.320285$   & $0.037572$  \\
      $v_z$      & $0.000000$   & $0.405848$  \\
      \hline\hline
    \end{tabular}
\end{table}
%%%%%%%%%%%%%%%%%%%%%%%%%%%%%%%%%%%%%%%%%%%%%%%%%%%%%%%%%%%%%%%%%%%%%%%%%%%%%%
%
First, note that $\theta_0 \ne 0$.  In fact, $\theta_0$ is closer to $\pi$ than $0$.
To test the optimality of the initial and final true anomalies, we use the Hamiltonian value conditions given by \eqref{eq:HVC_final_minTime}.  Note that the Hamiltonian is not the negative of unity for minimum-time as is the case of an autonomous system. Leveraging Chebfun's ability to determine  derivatives\cite{t_a_driscoll_chebfun_2014}, the computed left- and right-hand-sides of \eqref{eq:HVC_final_minTime} are presented in Table \ref{tab: compHVC_minTime}.
%
%%%%%%%%%%%%%%%%%%%%%%%%%%%%%%%%%%%%%%%%%%%%%%%%%%%%%%%%%%%%%%%%%%%%%%%%%%%%%%
\begin{table}[h!]
	\fontsize{10}{10}\selectfont
    \caption{Comparison of left and right-hand-sides of the Hamiltonian value condition given by \eqref{eq:HVC_final_minTime}}
    \label{tab: compHVC_minTime}
    \centering
    \begin{tabular}{c  r  r } % Column formatting
      \hline\hline
      Hamiltonian Value Condition & LHS  & RHS \\
      \hline
      Initial     & $-0.009312$   & $-0.009344$     \\
      Final       & $-0.002699$   & $-0.002824$     \\
      \hline\hline
    \end{tabular}
\end{table}
%%%%%%%%%%%%%%%%%%%%%%%%%%%%%%%%%%%%%%%%%%%%%%%%%%%%%%%%%%%%%%%%%%%%%%%%%%%%%%
%
It is apparent from this table that the Hamiltonian value condition is satisfied up to three decimal places.

To test the optimality of the remainder of the points presented in Table~\ref{tab: optimalStateTheta_minTime} we use the transversality conditions given by \eqref{eq:Costate_FinalValues}.  That is, we must have $\blam(\theta_0) = -\bnu_0$ and $\blam(\theta_f) = \bnu_f + (0, 0, \nu_{z_f}, 0, 0, 0)^T $.  From Table~\ref{tab: optimalStateTheta_minTime} we have $z_f = -0.077317 < z_f^U = -0.013556$. Hence, according to \eqref{eq: KKT_zf}, we must have $\nu_{z_f} = 0$.  This was, in fact, true.  Furthermore, because $\nu_{z_f} = 0$,  we must have $\blam(\theta_f) = \bnu_f$.  The satisfaction of these conditions is shown in Table \ref{tab: compCostateToCovector_minTime}.
%
%%%%%%%%%%%%%%%%%%%%%%%%%%%%%%%%%%%%%%%%%%%%%%%%%%%%%%%%%%%%%%%%%%%%%%%%%%%%%%
\begin{table}[h!]
	\fontsize{10}{10}\selectfont
    \caption{Comparison of initial and final costates, $\boldsymbol{\lambda}$, to corresponding covectors, $\bnu$}
    \label{tab: compCostateToCovector_minTime}
    \centering
    \begin{tabular}{ r  r  r  r } % Column formatting
      \hline\hline
      $\boldsymbol{\lambda}(\theta_0)$ & $-\bnu_0$ & $\boldsymbol{\lambda}(\theta_f)$ & $\bnu_f$ \\
      \hline
      $0.255241$   & $ 0.255236$  & $0.648531$   & $0.648517$  \\
      $0.000034$   & $0.000031$   & $-1.306714$  & $-1.306689$ \\
      $0.389260$   & $0.389259$   & $-1.204315$  & $-1.204302$ \\
      $0.326263$   & $0.326261$   & $0.349382$   & $0.349383$  \\
      $-0.466118$  & $-0.466120$  & $-0.131831$  & $-0.131832$ \\
      $-0.051572$  & $-0.051573$  & $0.651108$   & $0.651106$  \\
      \hline\hline
    \end{tabular}
    \\
\end{table}
%%%%%%%%%%%%%%%%%%%%%%%%%%%%%%%%%%%%%%%%%%%%%%%%%%%%%%%%%%%%%%%%%%%%%%%%%%%%%%
%
It is clear that the transversality conditions have been met to a sufficiently high accuracy. As a result of all of these numerics, the exit and entry points given by Table \ref{tab: optimalStateTheta_minTime} are declared extremal points.

\subsubsection{Dynamic Feasibility}
Figure~\ref{fig:L1_Lyap_to_SNRHO_minTime_Trajectory} also includes a dynamics propagation of the EOM using the Birkhoff-generated control history.
%
%%%%%%%%%%%%%%%%%%%%%%%%%%%%%%%%%%%%%%%%%%%%%%%%%%%%%%%%%%%%%%%%%%%%%%%%%%%%%%
\begin{figure}[hbt!]
\centering
\subfigure[$u_x$ control history]{%
    \includegraphics[width=0.5\columnwidth]{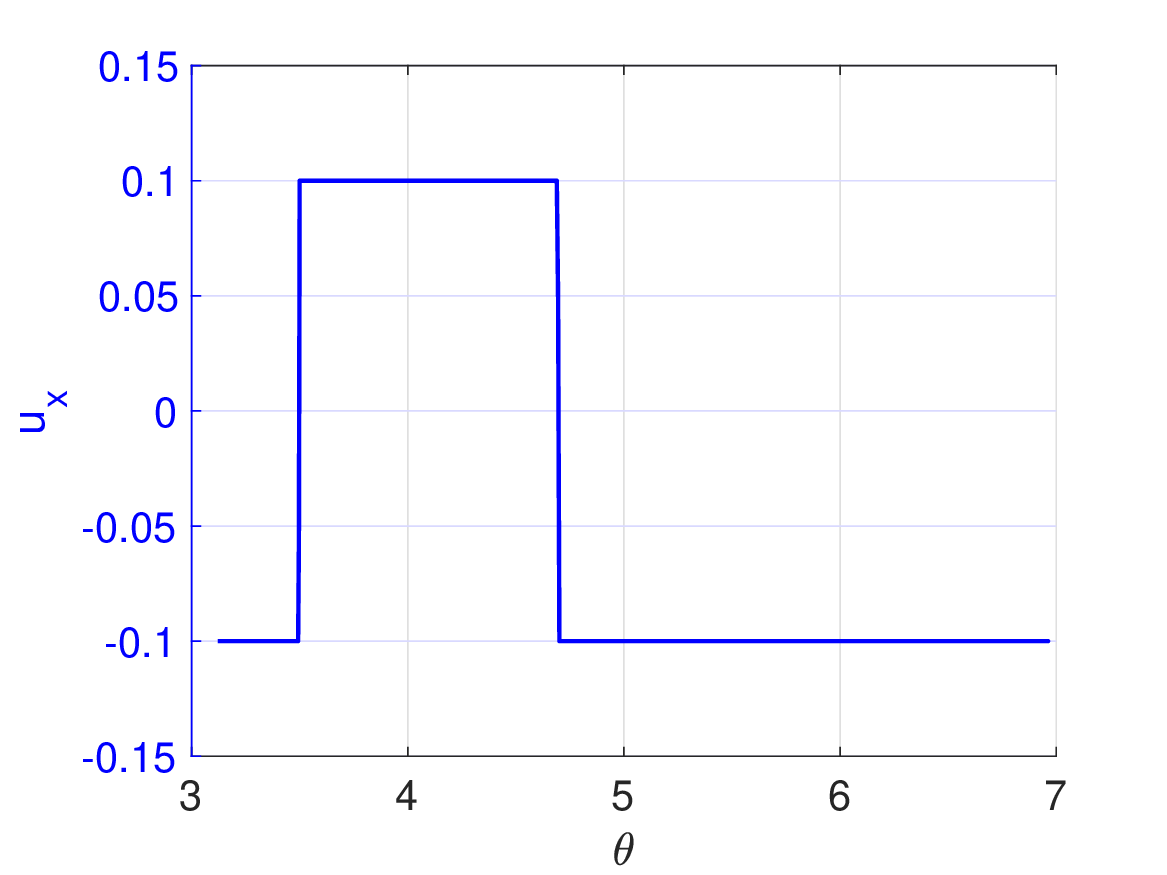}
    \label{fig:L1_Lyap_to_SNRHO_minTime_Ux}
}%
\\
\subfigure[$u_y$ control history]{%
    \includegraphics[width=0.5\columnwidth]{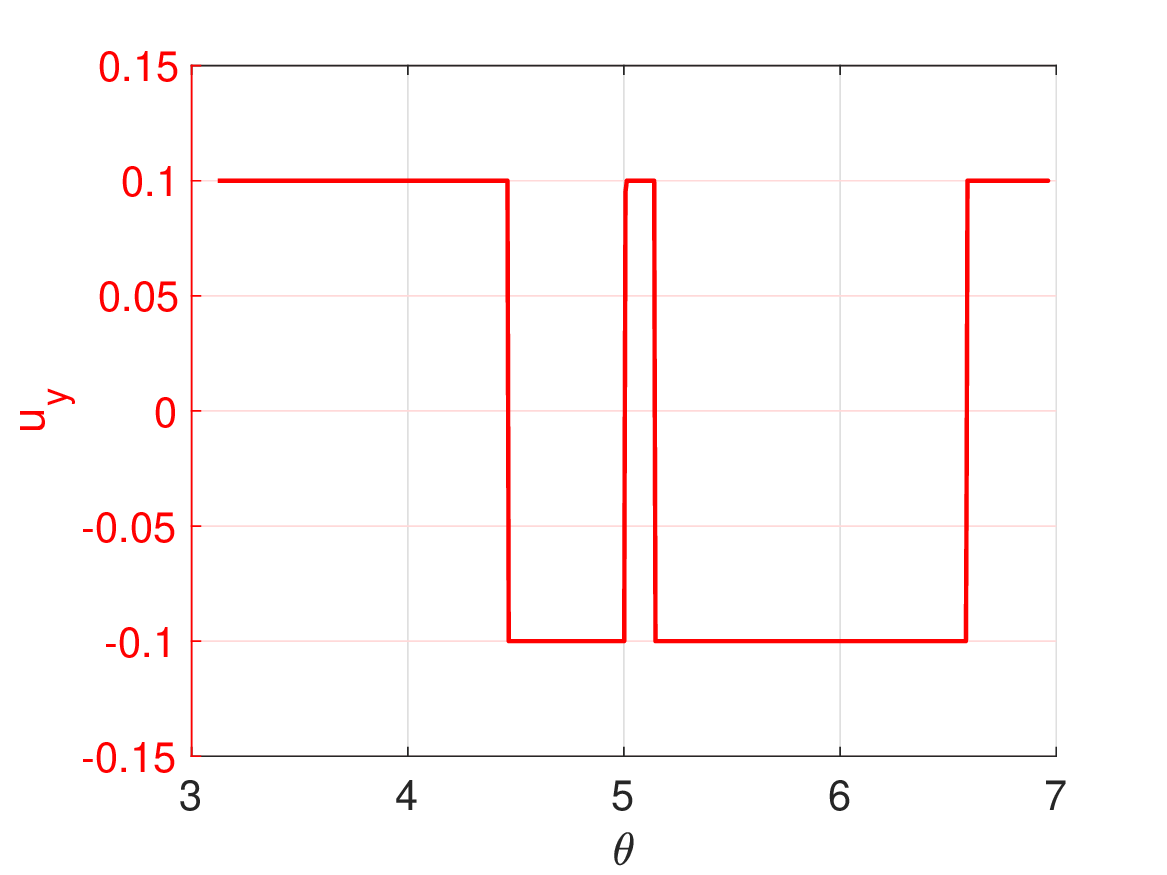}
    \label{fig:L1_Lyap_to_SNRHO_minTime_Uy}
}%
\\
\subfigure[$u_z$ control history]{%
    \includegraphics[width=0.5\columnwidth]{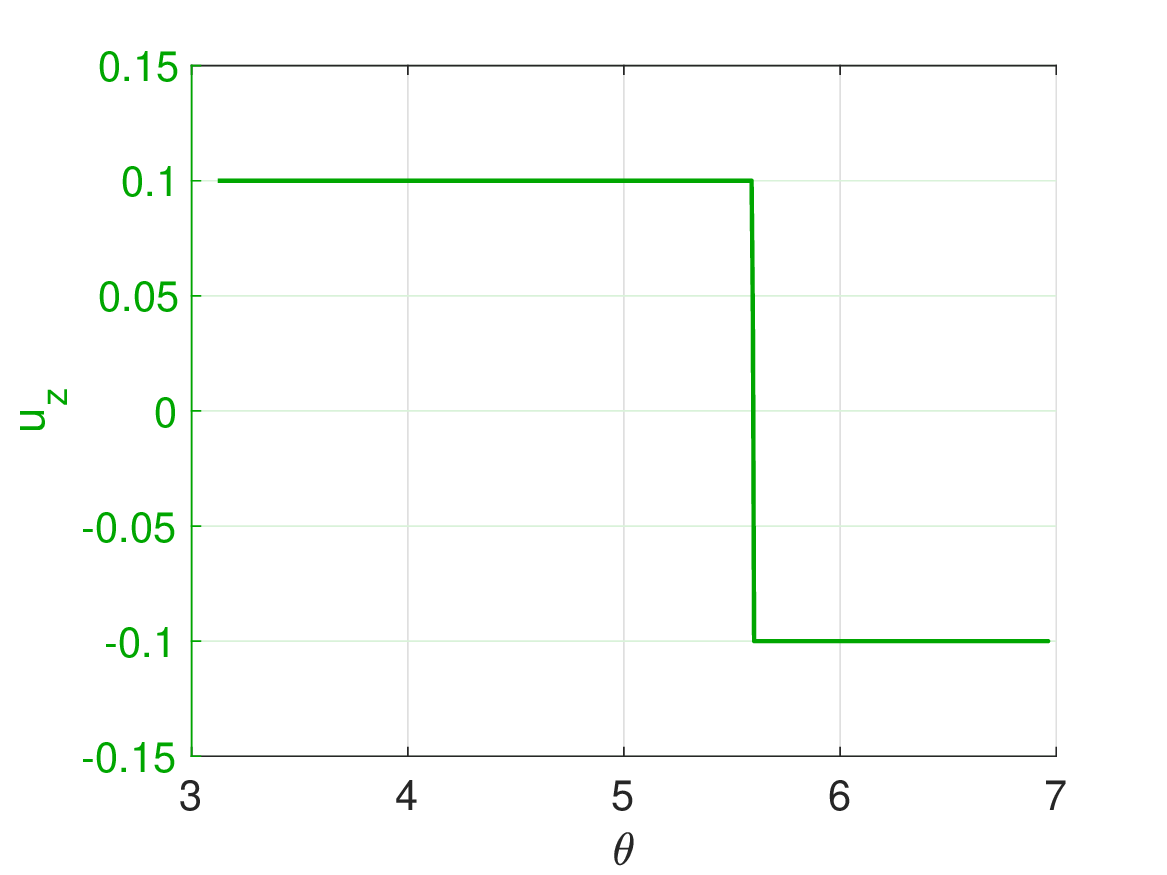}
    \label{fig:L1_Lyap_to_SNRHO_minTime_Uz}
}
\caption{Control trajectory for a minimum-time transfer between an $L_1$ Lyapunov and $L_2$ southern NRHO}
\label{fig:ControlHist_L1LyapToL2SNRHO_minTime}
\end{figure}
%%%%%%%%%%%%%%%%%%%%%%%%%%%%%%%%%%%%%%%%%%%%%%%%%%%%%%%%%%%%%%%%%%%%%%%%%%%%%%
%
The candidate optimal control is shown in Fig.~\ref{fig:ControlHist_L1LyapToL2SNRHO_minTime}.  We first note the following:
\begin{enumerate}
\item The controls plotted in Fig.~\ref{fig:ControlHist_L1LyapToL2SNRHO_minTime} are $u_x := u_x^a - u_x^b$, $u_y := u_y^a - u_y^b $ and $u_z := u_z^a - u_z^b$.
\item The controls are ``off'' prior to $ \theta_0 = 3.120787$  per the extremal departure point found in Table~\ref{tab: optimalStateTheta_minTime}.
\item Clearly, the controls are all bang-bang and consistent with typical minimum-time solutions\cite{bryson_applied_1978,longuski_optimal_2014,ross_primer_2015}.
\item There is no Gibbs phenomena in the generation of the controls.  This is because (see Sections~\ref{sec:intro} and \ref{sec:BirkIntro})  the control interpolant in a Birkhoff method is a non-polynomial.  See \cite{ross_universal_2023-1} and \cite{proulx_implementations_2023} for a complete explanation.
\end{enumerate}
It is apparent from Fig.~\ref{fig:L1_Lyap_to_SNRHO_minTime_Trajectory} that the continuous time propagation closely tracks the Birkhoff solution thereby providing an independent verification of feasibility.

\subsubsection{Extremality of the Bang-Bang Controls}

\begin{figure}[hbt!]
\centering
\subfigure[$u_x$ control history]{%
    \includegraphics[width=0.5\columnwidth]{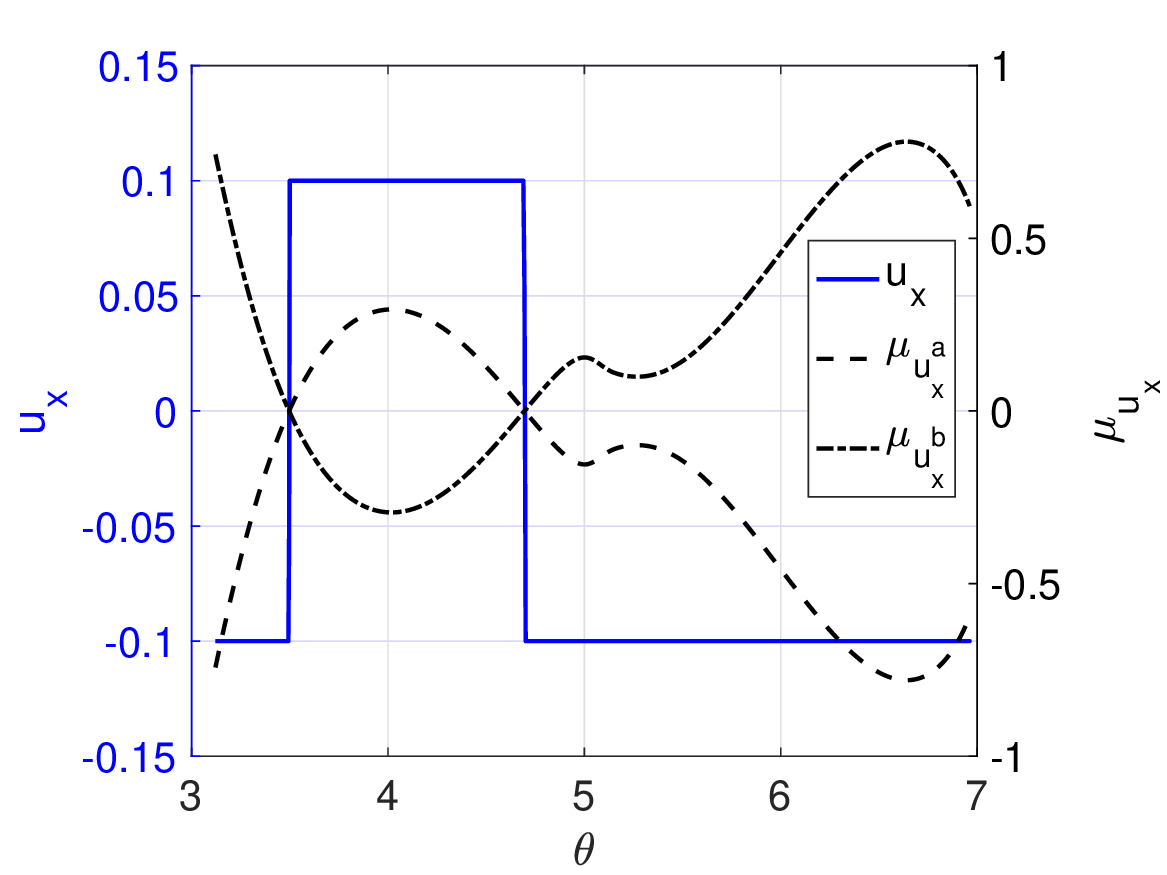}
    \label{fig:L1_Lyap_to_SNRHO_minTime_UxAndMu}
}%
\\
\subfigure[$u_y$ control history]{%
    \includegraphics[width=0.5\columnwidth]{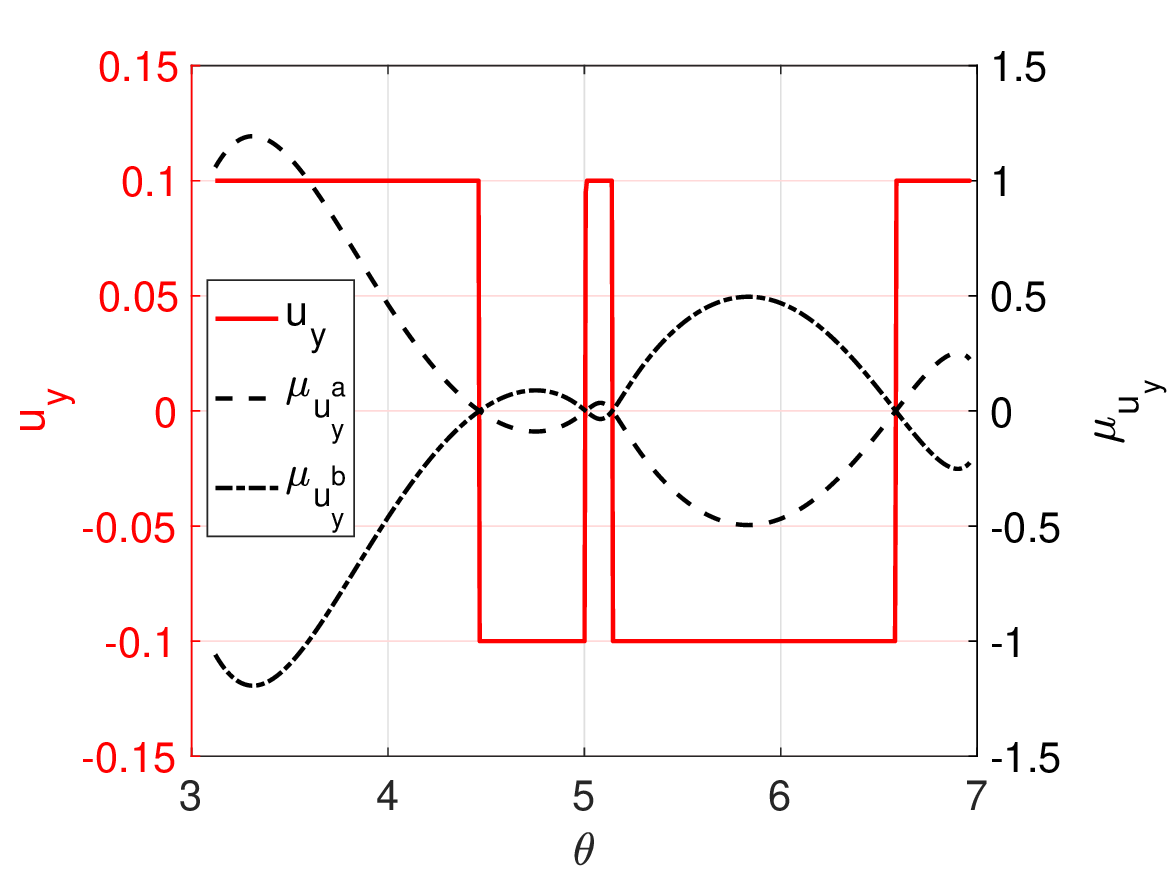}
    \label{fig:L1_Lyap_to_SNRHO_minTime_UyAndMu}
}%
\\
\subfigure[$u_z$ control history]{%
    \includegraphics[width=0.5\columnwidth]{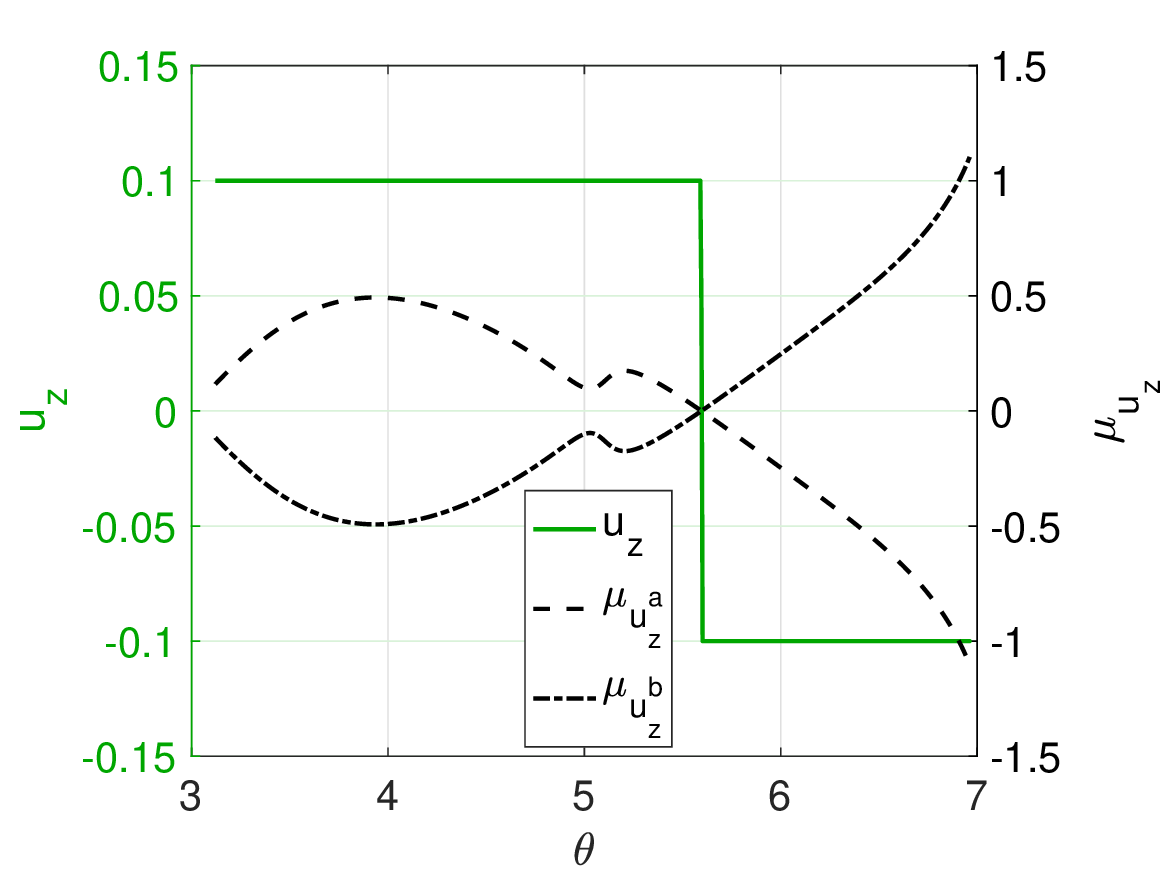}
    \label{fig:L1_Lyap_to_SNRHO_minTime_UzAndMu}
}
\caption{Control bound-constraint covector components for a minimum-time transfer between an $L_1$ Lyapunov and $L_2$ southern NRHO}
\label{fig: ControlHistAndMu_minTime}
\end{figure}

From \eqref{eq: Coplementarity}, it follows that if a bang-bang control solution is an extremal, we must have
\begin{equation}
u_{\xi} =
\begin{cases}
\begin{aligned}
 -u_{\xi}^{U}  &&\text{ if } && \mu_\xi^a \le 0 \text{ and } \mu_\xi^b \ge 0 \\
  u_{\xi}^{U}  &&\text{ if } && \mu_\xi^a \ge 0 \text{ and } \mu_\xi^b \le 0  \\
\end{aligned}
\end{cases}
\end{equation}
for $\xi = x, y, z$.
Figure~\ref{fig: ControlHistAndMu_minTime} shows that this condition is indeed satisfied to a high accuracy; hence we declare $\buf$ shown in Fig.~\ref{fig:ControlHist_L1LyapToL2SNRHO_minTime} to be an extremal control.
%
%%%%%%%%%%%%%%%%%%%%%%%%%%%%%%%%%%%%%%%%%%%%%%%%%%%%%%%%%%%%%%%%%%%%%%%%%%%%%%
%%%%%%%%%%%%%%%%%%%%%%%%%%%%%%%%%%%%%%%%%%%%%%%%%%%%%%%%%%%%%%%%%%%%%%%%%%%%%%

Finally, as a matter of completeness, we note that the computed minimum TOF with regard to true anomaly was $3.841000$ TU. The total non-dimensional $\Delta{prox}_1$ propellant consumption was $1.137681$ VU.

\subsection{Time-Limited Minimum-Propellant Solution}
\begin{figure}[h!]
 \centering
 \includegraphics[width=0.7\columnwidth]
 {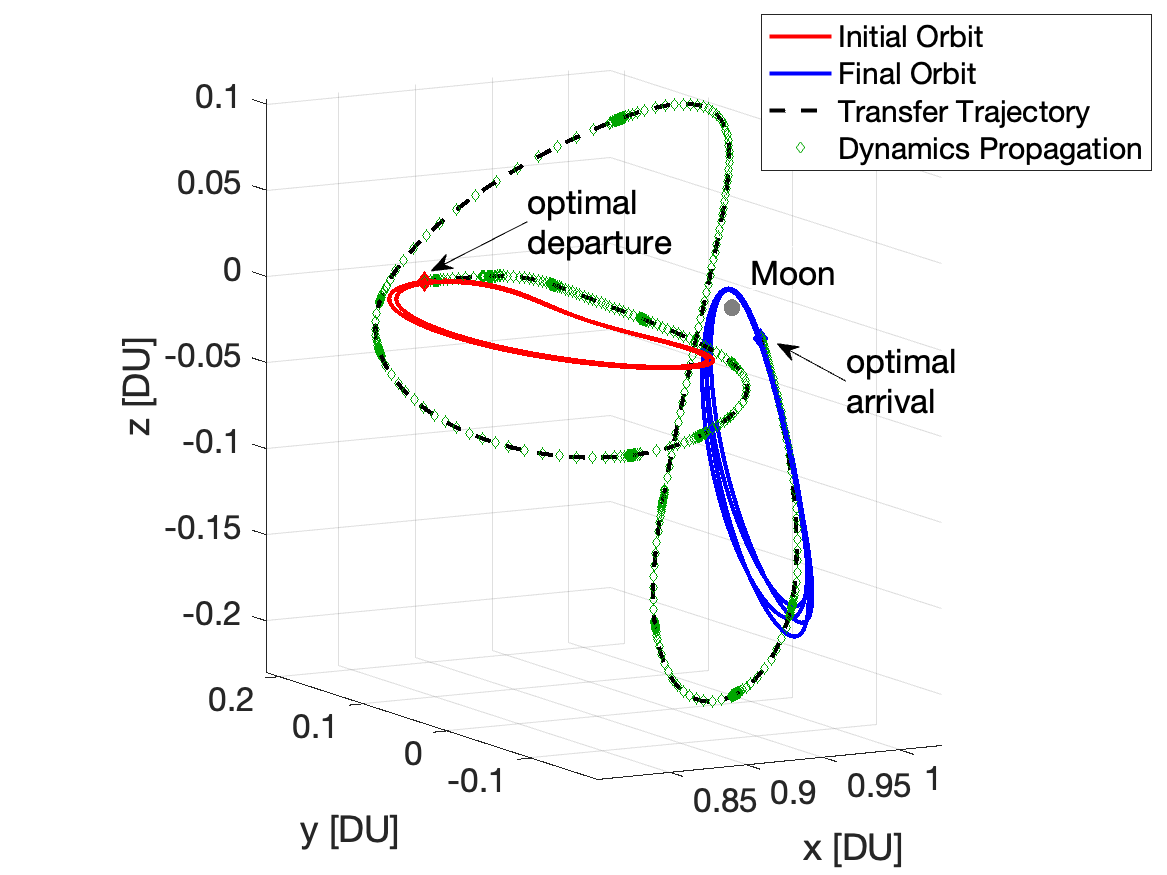}
     \caption{Time-limited, minimum-propellant transfer between an $L_1$ Lyapunov and $L_2$ southern NRHO.}
 \label{fig:L1_Lyap_to_SNRHO_timeBounded_Trajectory}
 \end{figure}

We now consider a practical, operational case wherein the cislunar spacecraft is required to complete the orbital transfer within a specified amount of time. If this specified time is less than $3.84$ TU, namely the minimum-time value, then it would be impossible to perform such a maneuver.  Thus, any bound on time must be greater than $3.84 \approx 1.2\pi$~TU. In this spirit, consider the orbit transfer problem with a TOF bounded by $2\pi$ TU; that is, a time bound equal to the orbital period of the primaries. The objective is to minimize propellant consumption with a constraint on transfer time.  The Birkhoff solution for this time-limited problem is shown in Fig. \ref{fig:L1_Lyap_to_SNRHO_timeBounded_Trajectory}. This solution converged at $N=1100$.
%
%%%%%%%%%%%%%%%%%%%%%%%%%%%%%%%%%%%%%%%%%%%%%%%%%%%%%%%%%%%%%%%%%%%%%%%%%%%%%%%
%%%%%%%%%%%%%%%%%%%%%%%%%%%%%%%%%%%%%%%%%%%%%%%%%%%%%%%%%%%%%%%%%%%%%%%%%%%%%
%
Comparing  Fig.~\ref{fig:L1_Lyap_to_SNRHO_timeBounded_Trajectory} with Fig.~\ref{fig:L1_Lyap_to_SNRHO_minTime_Trajectory} two points are immediately obvious: (a) the trajectories are starkly different and (b) the candidate optimal departure and arrival points are completely different.

\subsubsection{Extremality of the New Departure and Arrival Points}
The new departure and arrival points are shown in Table \ref{tab: optimalStateTheta_timeBound}.
%
%
%%%%%%%%%%%%%%%%%%%%%%%%%%%%%%%%%%%%%%%%%%%%%%%%%%%%%%%%%%%%%%%%%%%%%%%%%%%%%%
\begin{table}[h!]
	\fontsize{10}{10}\selectfont
    \caption{Optimal true anomaly and state vector for departure and arrival points for time-limited minimum-propellant transfer }
    \label{tab: optimalStateTheta_timeBound}
    \centering
    \begin{tabular}{l  c  c } % Column formatting
      \hline\hline
       coordinate & optimal departure $(x_0^*)$ & optimal arrival $(x_f^*)$ \\
      \hline
      $\theta$   & $0.801535$   & $7.084721$  \\
      $x$        & $0.873884$   & $0.990078$  \\
      $y$        & $0.158275$   & $-0.029906$  \\
      $z$        & $0.000000$   & $-0.013556$  \\
      $v_x$      & $0.113782$   & $-0.076072$  \\
      $v_y$      & $0.009185$   & $0.377910$  \\
      $v_z$      & $0.000000$   & $0.694526$  \\
      \hline\hline
    \end{tabular}
\end{table}
%%%%%%%%%%%%%%%%%%%%%%%%%%%%%%%%%%%%%%%%%%%%%%%%%%%%%%%%%%%%%%%%%%%%%%%%%%%%%%
From Table~\ref{tab: optimalStateTheta_timeBound} we have $\theta_f -\theta_0 = 2\pi$.  Obviously, the Birkhoff-theoretic spectral algorithm did not pick $\theta_0 = 0$ or the value of $\theta_0$ corresponding to the minimum-time solution presented in Table \ref{tab: optimalStateTheta_minTime}.  In performing the same optimality tests as before, we compute the left- and right-hand side of the Hamiltonian value conditions given by \eqref{eq:HVC_final_timeBound}. The value of $\nu_\theta$ was found to be equal to $0.002763$ (and hence satisfying the complementarity condition $\nu_\theta > 0$).  Carrying out the remainder of the computations yield the results shown in Table \ref{tab: compHVC_timeBound}.
%
%%%%%%%%%%%%%%%%%%%%%%%%%%%%%%%%%%%%%%%%%%%%%%%%%%%%%%%%%%%%%%%%%%%%%%%%%%%%%%
\begin{table}[h!]
	\fontsize{10}{10}\selectfont
    \caption{Comparison of left and right-hand-sides of the Hamiltonian value condition given by \eqref{eq:HVC_final_timeBound} for time-limited, minimum-propellant transfer}
    \label{tab: compHVC_timeBound}
    \centering
    \begin{tabular}{c  r  r } % Column formatting
      \hline\hline
      Hamiltonian Value Condition & LHS  & RHS \\
      \hline
      Initial     & $-0.011734$   & $-0.011727$     \\
      Final       & $-0.015880$   & $-0.016008$     \\
      \hline\hline
    \end{tabular}
\end{table}
%%%%%%%%%%%%%%%%%%%%%%%%%%%%%%%%%%%%%%%%%%%%%%%%%%%%%%%%%%%%%%%%%%%%%%%%%%%%%%
%
It is obvious that the Hamiltonian value condition is satisfied to a reasonable precision.

From the third column of Table \ref{tab: optimalStateTheta_timeBound}, we have $z_f = z_f^U = -0.013556$.  Hence, according to \eqref{eq: KKT_zf}, we must have $\nu_{z_f} > 0$.  The computed value of  $\nu_{z_f}$ is $ = 0.026670$.  Furthermore, according to \eqref{eq:Costate_FinalValues} we must have $\blam(\theta_0) = -\bnu_0$ and $\blam(\theta_f) = \bnu_f + (0, 0, \nu_{z_f}, 0, 0, 0)^T $.  The satisfaction of these conditions is shown in Table \ref{tab: compCostateToCovector_timeBound}.
%
%%%%%%%%%%%%%%%%%%%%%%%%%%%%%%%%%%%%%%%%%%%%%%%%%%%%%%%%%%%%%%%%%%%%%%%%%%%%%%
\begin{table}[h!]
	\fontsize{10}{10}\selectfont
    \caption{Comparison of initial and final costates, $\boldsymbol{\lambda}$, to corresponding covectors, $\bnu$, for time-limited minimum-propellant transfer}
    \label{tab: compCostateToCovector_timeBound}
    \centering
    \begin{tabular}{ r  r  r  r } % Column formatting
      \hline\hline
      $\boldsymbol{\lambda}(\theta_0)$ & $-\bnu_0$ & $\boldsymbol{\lambda}(\theta_f)$ & $\bnu_f$ \\
      \hline
      $-0.202685$ & $-0.202683$  & $2.269915$   & $2.269745$   \\
      $0.028382$  & $0.028381$   & $-1.311925$  & $-1.312024$  \\
      $0.000000$  & $0.002973$   & $2.387257$   & $2.387003^*$ \\
      $-0.116953$  & $-0.116952$  & $0.216632$  & $0.216638$   \\
      $-0.000430$  & $-0.000430$  & $-0.171882$  & $-0.171884$ \\
      $-0.121679$  & $-0.121678$  & $0.161905$   & $0.161911$  \\

      \hline\hline
    \end{tabular}
    \\
    \footnotesize{Note ($*$): Includes $\nu_{z_f} = 0.026670$}
\end{table}
%%%%%%%%%%%%%%%%%%%%%%%%%%%%%%%%%%%%%%%%%%%%%%%%%%%%%%%%%%%%%%%%%%%%%%%%%%%%%%
%

It is apparent from these computational results that the Birkhoff-theoretic spectral algorithm has found a solution that meets all of the transversality conditions within a reasonable numerical precision.

\subsubsection{Dynamic Feasibility and Extremality of the Controls}
\begin{figure}[hbt!]
\centering
\subfigure[$u_x$ control history]{%
    \includegraphics[width=0.5\columnwidth]{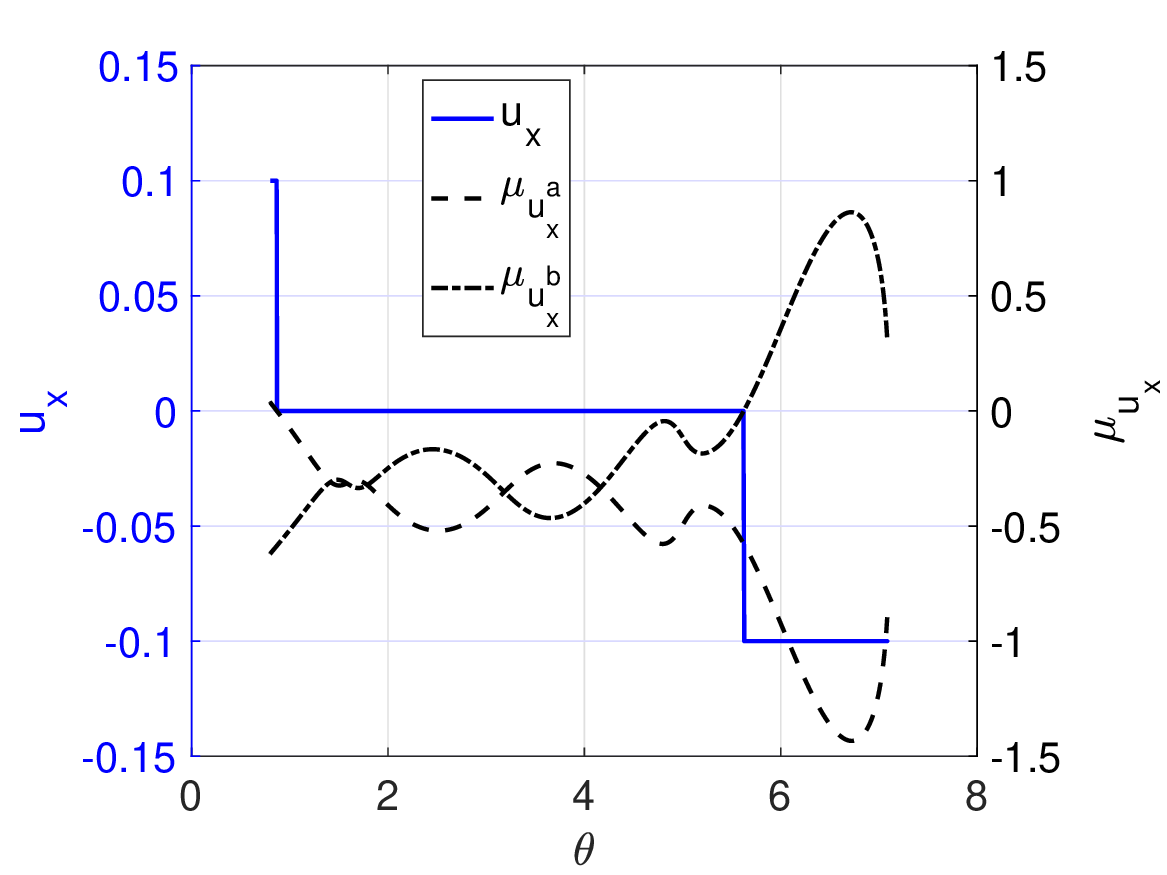}
    \label{fig:L1_Lyap_to_SNRHO_timeBounded_UxAndMu}
}%
\\
\subfigure[$u_y$ control history]{%
    \includegraphics[width=0.5\columnwidth]{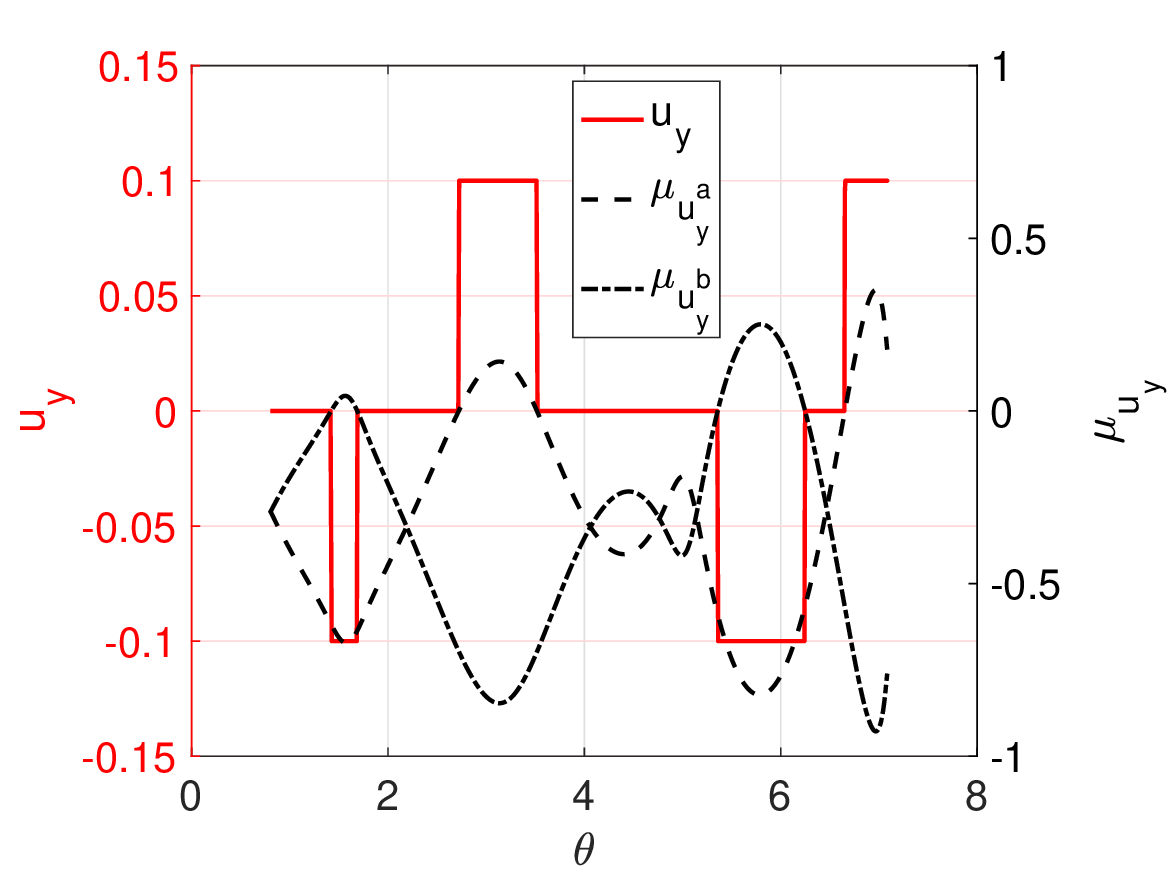}
    \label{fig:L1_Lyap_to_SNRHO_timeBounded_UyAndMu}
}%
\\
\subfigure[$u_z$ control history]{%
    \includegraphics[width=0.5\columnwidth]{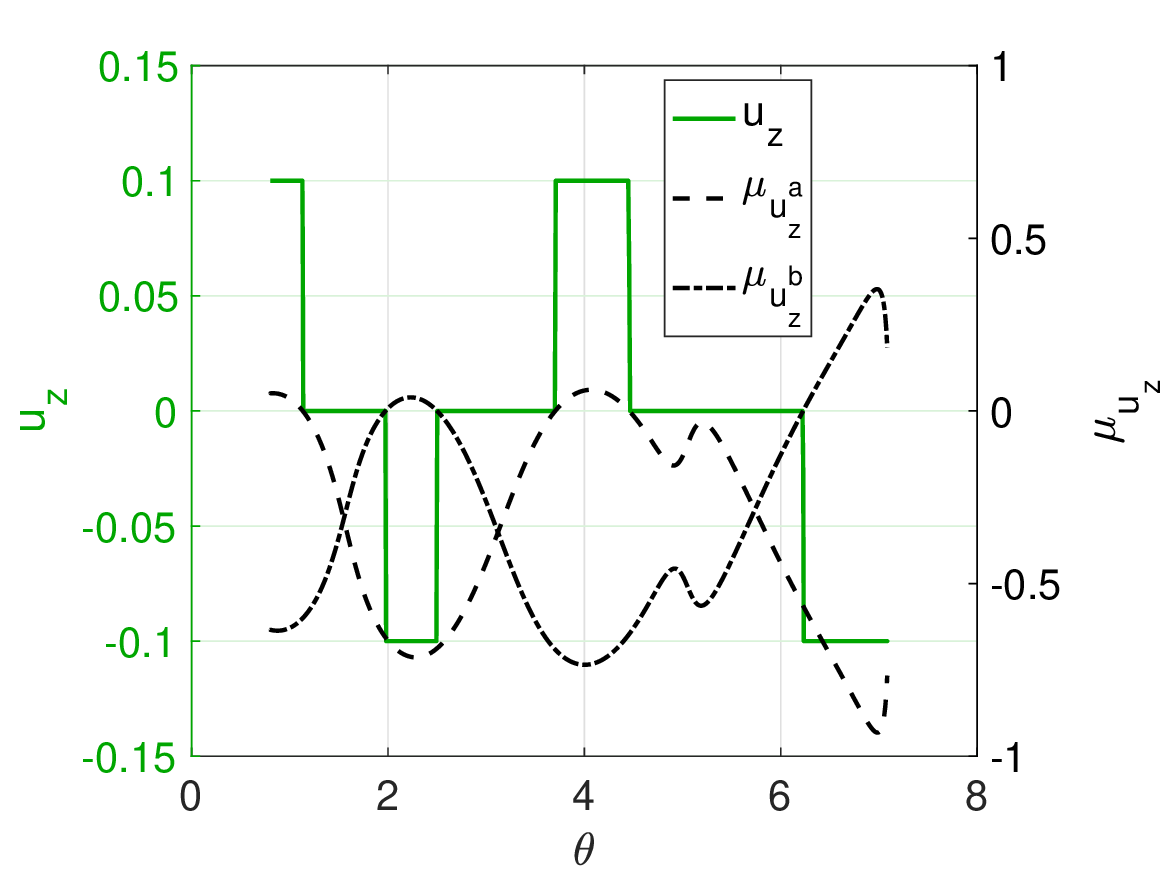}
    \label{fig:L1_Lyap_to_SNRHO_timeBounded_UzAndMu}
}
\caption{Control trajectory, and control bound constraint covector components for a time-bounded transfer between an $L_1$ Lyapunov and $L_2$ southern NRHO}
\label{fig:TrajectoryControlHistAndMu_L1LyapToL2SNRHO_timeBounded}
\end{figure}

For the purpose of brevity, we discuss in this subsection both topics of dynamic feasibility and extremality of the controls.  First, note that Fig.~\ref{fig:L1_Lyap_to_SNRHO_timeBounded_Trajectory} also includes a dynamics propagation of the EOM using the Birkhoff-generated control history. This control history is shown in  Fig.~\ref{fig:TrajectoryControlHistAndMu_L1LyapToL2SNRHO_timeBounded}.
%%%%%%%%%%%%%%%%%%%%%%%%%%%%%%%%%%%%%%%%%%%%%%%%%%%%%%%%%%%%%%%%%%%%%%%%%%%%%%
%%%%%%%%%%%%%%%%%%%%%%%%%%%%%%%%%%%%%%%%%%%%%%%%%%%%%%%%%%%%%%%%%%%%%%%%%%%%%%
%
From Figs.~\ref{fig:L1_Lyap_to_SNRHO_timeBounded_Trajectory} and \ref{fig:TrajectoryControlHistAndMu_L1LyapToL2SNRHO_timeBounded} we observe the following:
\begin{enumerate}
\item The continuous time propagation (shown in Fig~\ref{fig:L1_Lyap_to_SNRHO_timeBounded_Trajectory}) closely tracks the Birkhoff-computed solution thereby providing an independent verification of feasibility.
\item  Unlike the minimum-time case, the control histories have a bang-off-bang-structure.
\item For the bang-off-bang structure to be optimal, we must have
\begin{equation}
u_{\xi} =
\begin{cases}
\begin{aligned}
 -u_{\xi}^{U}  &&\text{ if } && \mu_\xi^a \le 0 \text{ and } \mu_\xi^b \ge 0 \\
 0  && \text{ if } && \mu_\xi^a \le 0 \text{ and } \mu_\xi^b \le 0  &&\text{ or }\quad \mu_\xi^a \ge 0 \text{ and } \mu_\xi^b \ge 0 \\
  u_{\xi}^{U}  &&\text{ if } && \mu_\xi^a \ge 0 \text{ and } \mu_\xi^b \le 0  \\
\end{aligned}
\end{cases}
\label{eq:complementarity_BangOffBang}
\end{equation}
for all $\xi \in \{ x, y, z \}$.
\item Figure~\ref{fig:TrajectoryControlHistAndMu_L1LyapToL2SNRHO_timeBounded} shows that the conditions presented in \eqref{eq:complementarity_BangOffBang} are indeed satisfied.  Note, in particular, the small initial ``impulse'' in Fig.~\ref{fig:L1_Lyap_to_SNRHO_timeBounded_UxAndMu}.  Obviously, a Birkhoff-theoretic method is well-positioned to capture this type of phenomenon.
\item  The condition  $\mu_\xi^a \ge 0 \text{ and } \mu_\xi^b \ge 0$ is allowed according to \eqref{eq:complementarity_BangOffBang} but it does not occur in Fig.~\ref{fig:TrajectoryControlHistAndMu_L1LyapToL2SNRHO_timeBounded}.  This situation implies that $u_{\xi} =  u_{\xi}^{U} - u_{\xi}^{U} = 0$. The reason this does not occur is because the
$\Delta{prop}_1$ cost would increase resulting in a non-optimal solution.  In other words, the fact that this does not happen is an indirect indication that the computed solution is indeed a time-limited minimum-propellant solution.
\end{enumerate}
%==================

The $\Delta{prox}_1$ cost for this time-limited problem was found to be $0.619465$ VU, a number less than the minimum-time value of $1.137681$ VU. Obviously, these numbers are consistent with optimal control theory.

\subsection{Time-Bounded Time-Free Minimum-Propellant Solution}
We now consider the case when the time bound is increased beyond the time-limited case of the previous subsection. That is instead of setting $\theta_f - \theta_0 \le 2\pi$, we consider a relaxation of this bound to $3\pi$.   The motivation of this time bound is to investigate the reduction in propellant, if any, with an allowable increase in transfer time. At this point, we do not know if the actual transfer time will indeed take up all of the allowable time because from a purely mathematical point of view the optimal transfer time is allowed to be $\theta_f -\theta_0 < 3\pi $.  In principle, this was true in the previous case as well but since the time bound here is larger than $2\pi$, this point takes on an added level of importance. Furthermore, as a precursor to the discussions to follow, we note that if the transfer time is satisfied with strict inequality, then $\nu_\theta$ must vanish per \eqref{eq: KKT_theta}.  In this case, the Hamiltonian value condition given by \eqref{eq:HVC_final_timeBound} reduces to that given by \eqref{eq:HVC_InitAndFinal}.  In other words, it is impossible to distinguish a time-free solution where $\theta_f - \theta_0 < \infty$ from $\theta_f -\theta_0 < 3\pi $.  Hence we refer to the latter as a time-bounded, time-free solution.

\begin{figure}[h!]
 \centering
 \includegraphics[width=0.7\columnwidth]{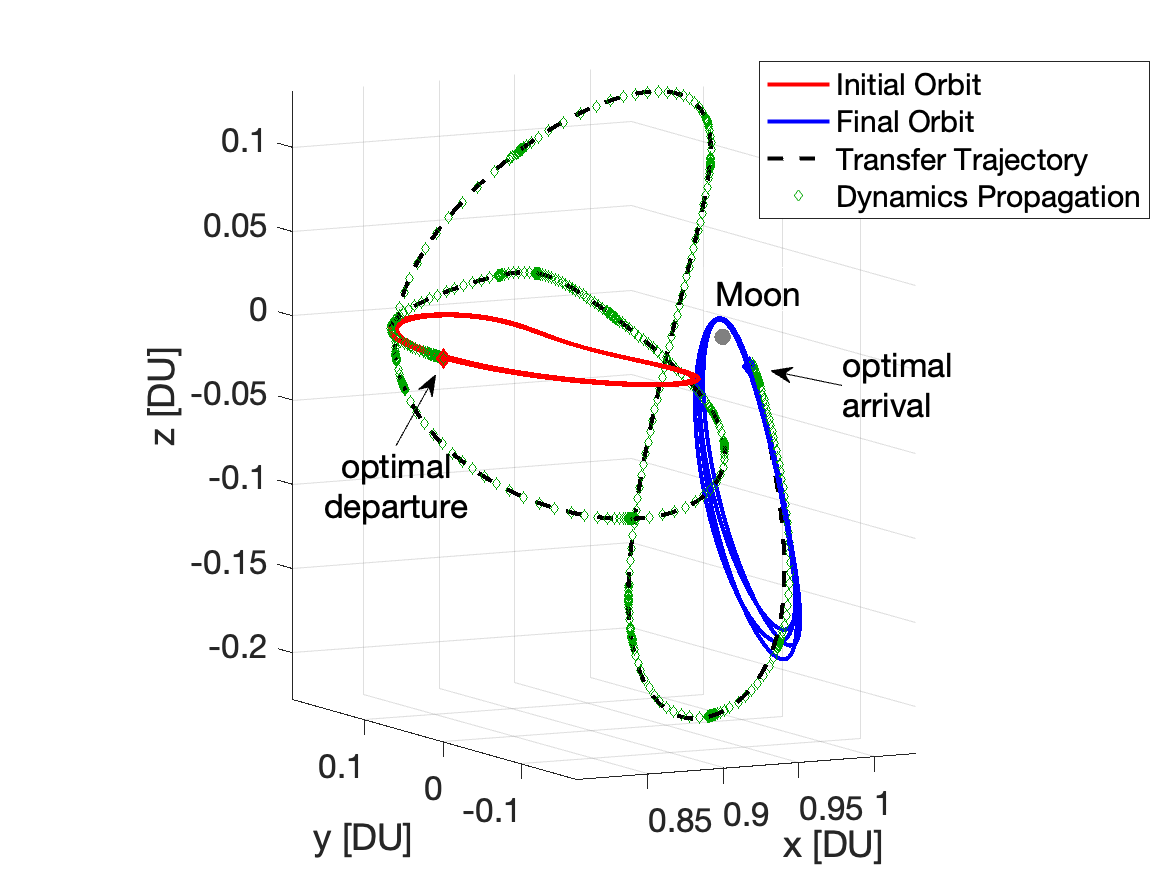}
 \caption{Minimum-propellant transfer between an $L_1$ Lyapunov and $L_2$ southern NRHO.}
 \label{fig:L1_Lyap_to_SNRHO_minProp_Trajectory}
 \end{figure}

A solution to this problem is shown in Figure~\ref{fig:L1_Lyap_to_SNRHO_minProp_Trajectory}. This solution converged at $N = 1400$.
%
%%%%%%%%%%%%%%%%%%%%%%%%%%%%%%%%%%%%%%%%%%%%%%%%%%%%%%%%%%%%%%%%%%%%%%%%%%%%%%
%%%%%%%%%%%%%%%%%%%%%%%%%%%%%%%%%%%%%%%%%%%%%%%%%%%%%%%%%%%%%%%%%%%%%%%%%%%%%
%
Comparing Fig.~\ref{fig:L1_Lyap_to_SNRHO_minProp_Trajectory} with Fig.~\ref{fig:L1_Lyap_to_SNRHO_timeBounded_Trajectory}, it is obvious that not only are the trajectories quite different but so are the departure points.  On the other hand, the candidate arrival points appear to be the same.  In fact, as shown in Table~\ref{tab: optimalStateTheta_minProp}, the candidate optimal arrival point is exactly the same as in Table~\ref{tab: optimalStateTheta_timeBound}.
%
%%%%%%%%%%%%%%%%%%%%%%%%%%%%%%%%%%%%%%%%%%%%%%%%%%%%%%%%%%%%%%%%%%%%%%%%%%%%%%
\begin{table}[h!]
	\fontsize{10}{10}\selectfont
    \caption{Optimal true anomaly and state vector for departure and arrival points for minimum-propellant transfer }
    \label{tab: optimalStateTheta_minProp}
    \centering
    \begin{tabular}{l  c  c } % Column formatting
      \hline\hline
       coordinate & optimal departure $(x_0^*)$ & optimal arrival $(x_f^*)$ \\
      \hline
      $\theta$   & $0.000000$   & $7.084721$  \\
      $x$        & $0.803019$   & $0.990078$  \\
      $y$        & $0.000000$   & $-0.029906$  \\
      $z$        & $0.000000$   & $-0.013556$  \\
      $v_x$      & $0.000000$   & $-0.076072$  \\
      $v_y$      & $0.316692$   & $0.377910$  \\
      $v_z$      & $0.000000$   & $0.694526$  \\
      \hline\hline
    \end{tabular}
\end{table}
%%%%%%%%%%%%%%%%%%%%%%%%%%%%%%%%%%%%%%%%%%%%%%%%%%%%%%%%%%%%%%%%%%%%%%%%%%%%%%
%
This is because the spacecraft arrives at the targeting point that is bounded by $z_f^U = -0.013556$. Furthermore, note that $\theta_0$ is now zero.  In other words, the Birkhoff-theoretic spectral algorithm has adjusted the departure point to an earlier point (relative to the time-limited maneuver discussed in the previous subsection) while fixing the arrival point to meet the targeting constraint. To test the optimality of these points, we once again invoke the Hamiltonian value conditions given by \eqref{eq:HVC_final_timeBound}.  The value of $\nu_\theta$ was computed (by the Birkhoff-theoretic spectral algorithm) to be zero suggesting that $\theta_f - \theta_0 < 3\pi$.  This is indeed the case as is apparent from the first row of Table~\ref{tab: optimalStateTheta_minProp}.  Carrying out the remainder of the computations indicated in \eqref{eq:HVC_final_timeBound} yields the results shown in  Table~\ref{tab: compHVC_minProp}.
%
%%%%%%%%%%%%%%%%%%%%%%%%%%%%%%%%%%%%%%%%%%%%%%%%%%%%%%%%%%%%%%%%%%%%%%%%%%%%%%
\begin{table}[h!]
	\fontsize{10}{10}\selectfont
    \caption{Comparison of left and right-hand-sides of the Hamiltonian value condition given by \eqref{eq:HVC_final_timeBound}}
    \label{tab: compHVC_minProp}
    \centering
    \begin{tabular}{c  r  r } % Column formatting
      \hline\hline
      Hamiltonian Value Condition & LHS  & RHS \\
      \hline
      Initial     & $-0.015566$   & $-0.015521$     \\
      Final       & $-0.021272$   & $-0.021471$     \\
      \hline\hline
    \end{tabular}
\end{table}
%%%%%%%%%%%%%%%%%%%%%%%%%%%%%%%%%%%%%%%%%%%%%%%%%%%%%%%%%%%%%%%%%%%%%%%%%%%%%%
%

To test the optimality of the remainder of the points presented in Table~\ref{tab: optimalStateTheta_minProp} we once again invoke the transversality conditions.
From Table~\ref{tab: optimalStateTheta_minProp} we have $z_f = z_f^U = -0.013556$.  Hence, according to \eqref{eq: KKT_zf}, we must have $\nu_{z_f} > 0$.  The computed value of  $\nu_{z_f}$ is $ 0.044024$.  Note that this value of $\nu_{z_f}$ is not the same as the one obtained earlier (Cf.~Table~\ref{tab: compCostateToCovector_timeBound}).  In any case, according to \eqref{eq:Costate_FinalValues} we must have $\blam(\theta_0) = -\bnu_0$ and $\blam(\theta_f) = \bnu_f + (0, 0, \nu_{z_f}, 0, 0, 0)^T $.  The satisfaction of these conditions is shown in Table \ref{tab: compCostateToCovector_minProp}.
%
%%%%%%%%%%%%%%%%%%%%%%%%%%%%%%%%%%%%%%%%%%%%%%%%%%%%%%%%%%%%%%%%%%%%%%%%%%%%%%
\begin{table}[h!]
	\fontsize{10}{10}\selectfont
    \caption{Comparison of initial and final costates, $\boldsymbol{\lambda}$, to corresponding covectors, $\bnu$, for minimum-propellant transfer}
    \label{tab: compCostateToCovector_minProp}
    \centering
    \begin{tabular}{ r  r  r  r } % Column formatting
      \hline\hline
      $\boldsymbol{\lambda}(\theta_0)$ & $-\bnu_0$ & $\boldsymbol{\lambda}(\theta_f)$ & $\bnu_f$ \\
      \hline
      $-0.368809 $  & $-0.368807$  & $2.508892$   & $2.508753$  \\
      $0.063417$   & $0.063417$   & $0.189979$   & $0.189748$   \\
      $0.238949$   & $0.238949$   & $3.654012$   & $3.653795^*$ \\
      $-0.107514$  & $-0.107514$  & $0.279175$   & $0.279180$   \\
      $-0.114873$  & $-0.114873$  & $-0.270395$  & $-0.270393$  \\
      $-0.040893$  & $-0.040894$  & $0.062188$   & $0.062195$   \\
      \hline\hline
    \end{tabular}
    \\
    \footnotesize{Note ($*$): Includes $\nu_{z_f} = 0.044024$}
\end{table}
%%%%%%%%%%%%%%%%%%%%%%%%%%%%%%%%%%%%%%%%%%%%%%%%%%%%%%%%%%%%%%%%%%%%%%%%%%%%%%

\begin{figure}[hbt!]
\centering
\subfigure[$u_x$ control history]{%
    \includegraphics[width=0.5\columnwidth]{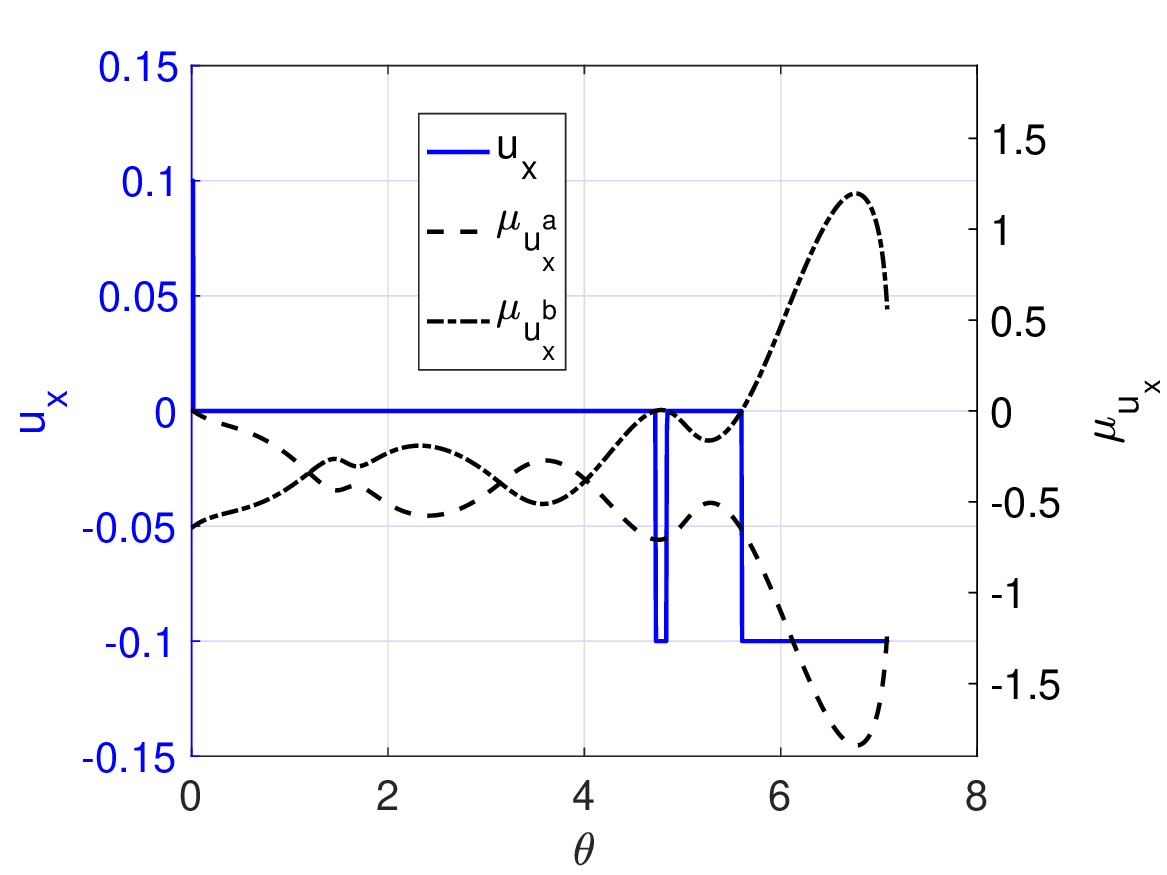}
    \label{fig:L1_Lyap_to_SNRHO_minProp_UxAndMu}
}%
\\
\subfigure[$u_y$ control history]{%
    \includegraphics[width=0.5\columnwidth]{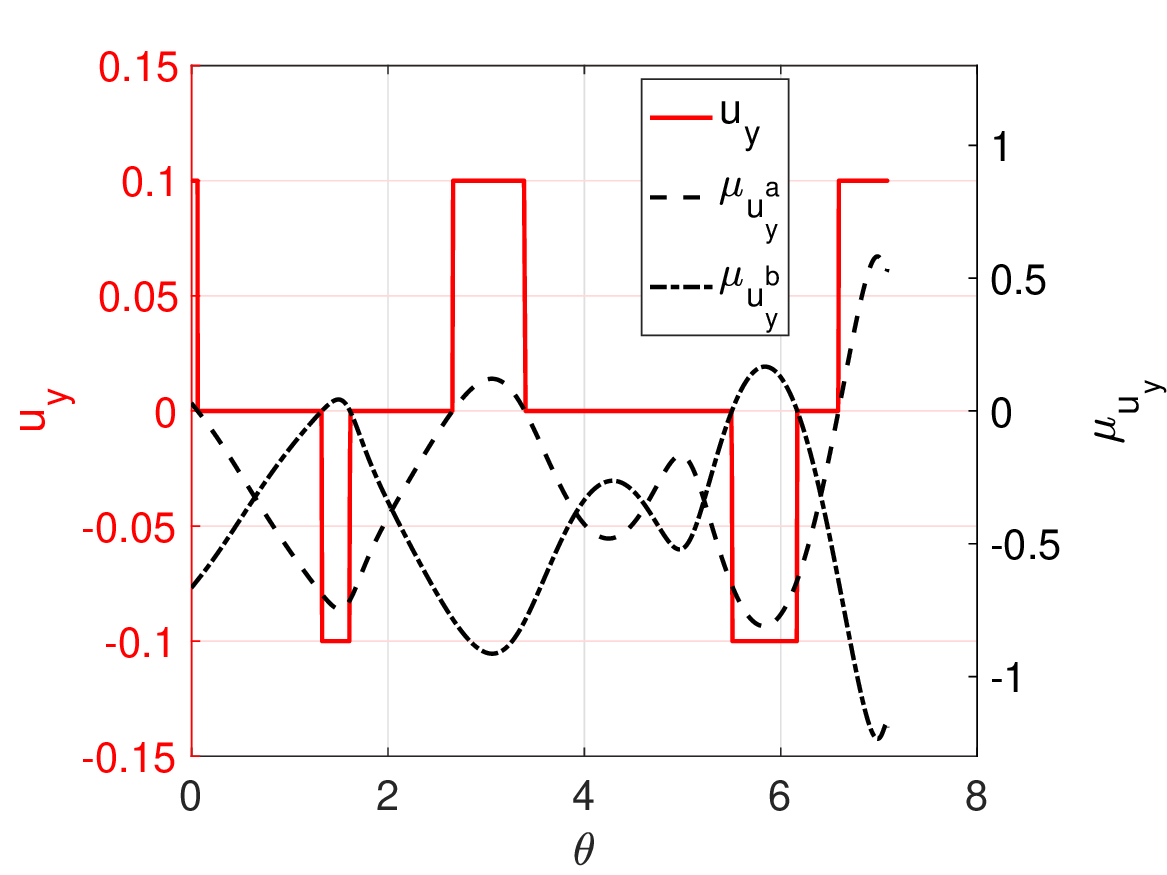}
    \label{fig:L1_Lyap_to_SNRHO_minProp_UyAndMu}
}%
\\
\subfigure[$u_z$ control history]{%
    \includegraphics[width=0.5\columnwidth]{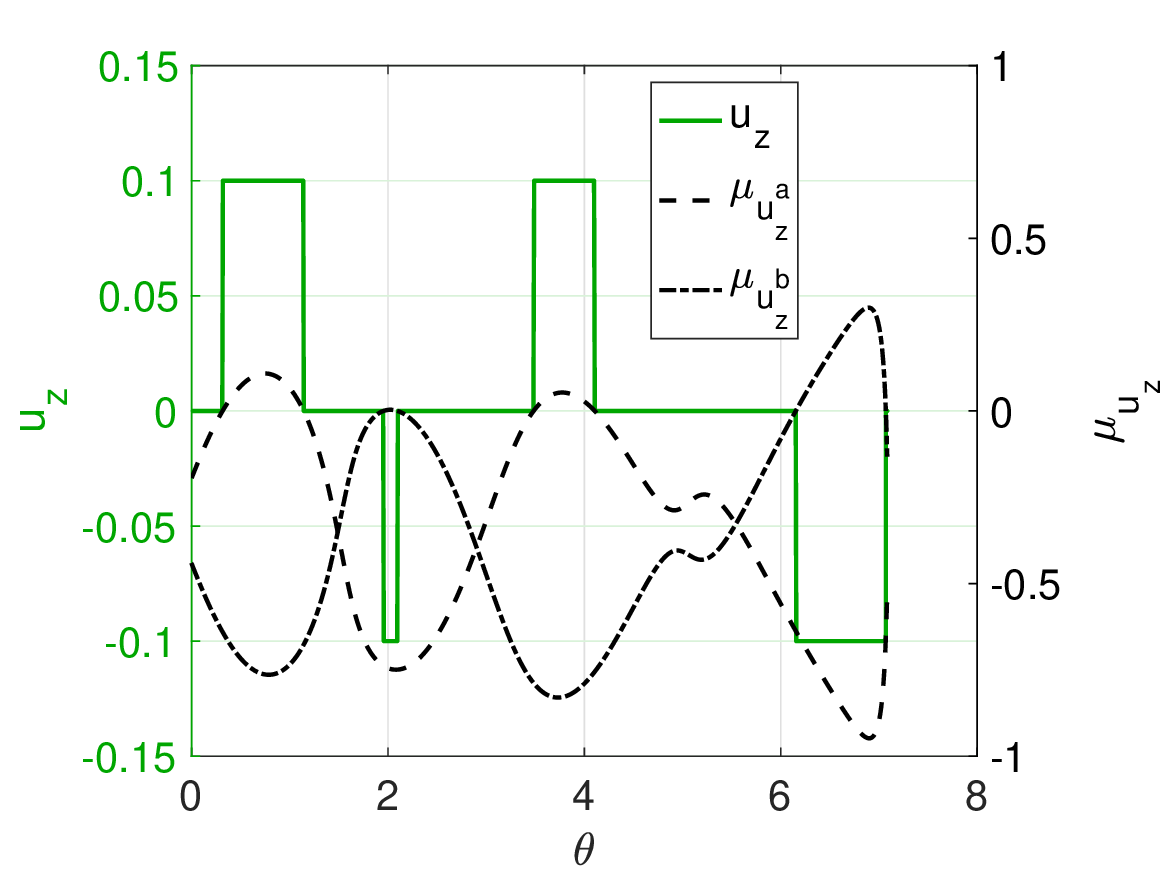}
    \label{fig:L1_Lyap_to_SNRHO_minProp_UzAndMu}
}
\caption{Control trajectory and control bound constraint covector components for a minimum-propellant transfer between an $L_1$ Lyapunov and $L_2$ southern NRHO}
\label{fig:ControlHistAndMu_L1LyapToL2SNRHO_minProp}
\end{figure}

As a final point of discussion, we note that Fig.~\ref{fig:L1_Lyap_to_SNRHO_minProp_Trajectory} also includes a dynamics propagation of the EOM using the Birkhoff-generated control history. This control history is shown in Fig.~\ref{fig:ControlHistAndMu_L1LyapToL2SNRHO_minProp}.
%
%%%%%%%%%%%%%%%%%%%%%%%%%%%%%%%%%%%%%%%%%%%%%%%%%%%%%%%%%%%%%%%%%%%%%%%%%%%%%%
%%%%%%%%%%%%%%%%%%%%%%%%%%%%%%%%%%%%%%%%%%%%%%%%%%%%%%%%%%%%%%%%%%%%%%%%%%%%%%
%
From Figs.~\ref{fig:L1_Lyap_to_SNRHO_minProp_Trajectory} and \ref{fig:ControlHistAndMu_L1LyapToL2SNRHO_minProp} we observe the following:
\begin{enumerate}
\item The continuous time propagation (shown in Fig~\ref{fig:L1_Lyap_to_SNRHO_minProp_Trajectory}) closely tracks the Birkhoff-computed solution thereby providing an independent verification of feasibility.
\item Of particular note is the $u_x$ control history.  There is a very small finite-burn thrust-arc near $\theta = 5$. Obviously, the Birkhoff-theoretic spectral algorithm had no trouble in locating and computing this small-duration finite thrust-arc.
\item Figure~\ref{fig:ControlHistAndMu_L1LyapToL2SNRHO_minProp} shows that the conditions presented in \eqref{eq:complementarity_BangOffBang} are indeed satisfied. Note, in particular, Fig.~\ref{fig:L1_Lyap_to_SNRHO_minProp_UxAndMu}.
\item  As in the case of the time-limited situation, the condition  $\mu_\xi^a \ge 0 \text{ and } \mu_\xi^b \ge 0$ is allowed according to \eqref{eq:complementarity_BangOffBang} but does not occur in Fig.~\ref{fig:ControlHistAndMu_L1LyapToL2SNRHO_minProp}.  As noted earlier, this is because the propellant consumption would increase resulting in a non-optimal solution.
\end{enumerate}
The $\Delta{prox}_1$ cost for this problem was found to be $0.609603$ VU, a number that is approximately $1.6$\% less than the time-limed problem presented earlier.  Given that the time of flight for this problem is approximately $2.26 \pi$ or about $13$\% greater than the time-limited problem, a mission designer can now make an informed decision if the added time of flight is worth the savings in propellant.  Alternatively, a mission designer may seek a lower propellant solution by increasing the maximum time of flight to, say, $4\pi$. It turns out that a further lowering of the propellant consumption is possible with an increased maximum time of flight, the details of which are discussed in \cite{dixon-ross-TN_2025}.

%%%%%%%%%%
%Conclusion
%%%%%%%%%%
\section{Conclusion}

On a present-day laptop, a Chebyshev representation of a libration point orbit can be performed in fractions of a second even when the order of the polynomial is over ten thousand.  This rapid computation is made possible through the use of a fast Fourier transform. As a result, the challenge associated in a representation of a libration point orbit can be overcome through the high-accuracy capability of Chebyshev polynomials.  Note also that computer memory requirements for storing a ten-thousand-order Chebyshev polynomial is only about $0.08$ MB.  Clearly the challenge is not technological but an irrational fear of very high-order polynomials.

The use of a recently-developed $\Delta{prox}$-functional allowed us to avoid selecting a particular value of the specific impulse of the propellant for trajectory optimization. Thus, the results presented in this paper are agnostic to the propulsion parameters except for the finite value of the maximum thrust capability.  By directly modeling finite thrust values, there is no need to account for gravity losses and other corrections that are necessary to map impulse solutions to practical engine capabilities.  In other words, the universal Birkhoff theory for trajectory optimization allows one to perform in one shot an end-to-end trajectory optimization in the in the elliptic restricted three-body problem.  Furthermore, the Birkhoff-theoretic method is capable of generating accurate solutions over thousands of grid points without using any high-performance computational devices other than a present-day ordinary laptop.

The EOM in the elliptic restricted three-body problem are non-autonomous.  This fact poses a substantial challenge in using dynamical systems theory for trajectory optimization, particularly for time-bound finite-thrust maneuvers.  Obviously, these issues are not a source of hinderance for a Birkhoff-theoretic method.  Furthermore, when coupled with the guess-free fast spectral algorithm, a Birkhoff-theoretic method does not require a good guess to converge to an extremal solution.  The independent feasibility tests show that these solutions can be implemented in a manner that is completely agnostic to the computational method used for trajectory optimization, namely the universal Birkhoff theory. Consequently, flight operators need not posses any knowledge of Birkhoff-theoretic methods for a practical implementation. As a result, the methods proposed in this paper can be deployed for operational use where rapid solutions are needed to make informed decisions on time-critical mission planning.

%%%%%%%%%%
%Acknowledgment
%%%%%%%%%%

\section{Acknowledgments}
We gratefully acknowledge funding for this research provided by the Air Force Office of Scientific Research (AFOSR). The views and conclusions contained herein are those of the authors and should not be interpreted as necessarily representing the official policies or endorsements, either expressed or implied, of the AFOSR or the U.S. Government.  The authors would also like to thank Dr. Ronald J. Proulx of the Naval Postgraduate School for his substantial help in understanding certain intricacies associated with the computational implementation of the Birkhoff theoretic method for trajectory optimization.

%%%%%%%%
%Bibliography
%%%%%%%%


\begin{thebibliography}{10}
\newcommand{\enquote}[1]{``#1''}

\bibitem{gmez_study_1993}
G\'omez, G., Jorba, A., Masdemont, J., and Sim\'o, C., \enquote{Study of the
  Transfer from the {Earth} to a Halo Orbit Around the Equilibrium {Point L}1,}
  {\em Celestial Mechanics \& Dynamical Astronomy\/}, Vol.~56, No.~4, Aug.
  1993, pp.~541--562. \\
  https://doi.org/10.1007/bf00696185

\bibitem{belbruno_capture_2004}
Belbruno, E., {\em Capture Dynamics and Chaotic Motions in Celestial Mechanics:
  With Applications to the Construction of Low Energy Transfers\/}, Princeton
  University Press, Princeton, 2004.\\
  https://doi.org/10.2307/j.ctv301g1f

\bibitem{koon_dynamical_2011}
Koon, W., Lo, M., Marsden, J., and Ross, S.~D., {\em Dynamical {Systems}, the
  {Three}-{Body} {Problem} and {Space} {Mission} {Design}\/}, Marsden Books,
  2011.

\bibitem{parker_lowenergy_2014}
Parker, J.~S. and Anderson, R.~L., {\em Low‐{Energy} {Lunar} {Trajectory}
  {Design}\/}, Wiley, 1st ed., June 2014. \\
  https://doi.org/10.1002/9781118855065

\bibitem{howell_application_1997}
Howell, K.~C., Barden, B.~T., and Lo, M.~W., \enquote{Application of
  {Dynamical} {Systems} {Theory} to {Trajectory} {Design} for a {Libration}
  {Point} {Mission},} {\em The Journal of the Astronautical Sciences\/},
  Vol.~45, No.~2, June 1997, pp.~161--178. \\
  https://doi.org/10.2514/6.1996-3602

\bibitem{longuski_optimal_2014}
Longuski, J., Guzmán, J.~J., and Prussing, J.~E., {\em Optimal Control with
  Aerospace Applications\/}, Space technology library, Springer ; Published
  jointly by Microcosm Press, New York, New York : El Segundo, California,
  2014. \\
  https://doi.org/10.1007/978-1-4614-8945-0

\bibitem{clarke_nonsmooth_1998}
Clarke, F.~H., Ledyaev, Y.~S., Stern, R.~J., and Wolenski, {\em Nonsmooth
  Analysis and Control theory\/}, No. 178 in Graduate Texts in Mathematics,
  Springer, New York Heidelberg, 1998. \\
https://doi.org/10.1007/b97650

\bibitem{bryson_applied_1978}
Bryson, A.~E. and Ho, Y.-C., {\em Applied Optimal Control: Optimization,
  Estimation, and Control\/}, Cambridge University Press, Cambridge, rev. print
  ed., 1978. \\
  https://doi.org/10.1201/9781315137667

\bibitem{ross_primer_2015}
Ross, I.~M., \textit{A Primer on Pontryagin's Principle in Optimal Control}, Second Edition, Collegiate Publishers, San Francisco, CA, 2015. https://www.amazon.com/Primer-Pontryagins-Principle-Optimal-Control/dp/0984357114.

\bibitem{park_leveraging_2022}
Park, B. and Howell, K., \enquote{Leveraging {Intermediate} {Dynamical}
  {Models} for {Transitioning} from the {Circular} {Restricted} {Three}-{Body}
  {Problem} to an {Ephemeris} {Model},} Ammerican Astronautical Society,
  Charlotte, North Carolina, Aug. 2022.


\bibitem{finite-burn-1966}
Robbins, H. M., ``An Analytical Study of Impulsive Approximation,''  \textit{AIAA Journal}, 4/8, 1417--1423. \\
https://doi.org/10.2514/3.3687

\bibitem{finite-burn-2004}
Ocampo, C., ``Finite Burn Maneuver Modeling for a Generalized Spacecraft Trajectory
Design and Optimization System,'' \textit{Annals of the N.Y. Academy of  Sciences},
1017: 210--233 (2004). \\
https://doi.org/10.1196/annals.1311.013

\bibitem{mingotti_low-energy_2009}
Mingotti, G., Topputo, F., and Bernelli-Zazzera, F., \enquote{Low-Energy,
  Low-Thrust Transfers to the {Moon},} {\em Celestial Mechanics and Dynamical
  Astronomy\/}, Vol.~105, No. 1-3, Nov. 2009, pp.~61--74. \\
  https://doi.org/10.1007/s10569-009-9220-7

\bibitem{chupin_transfer_2018}
Chupin, M., Haberkorn, T., and Tr\'{e}lat, E., \enquote{Transfer {Between}
  {Invariant} {Manifolds}: {From} {Impulse} {Transfer} to {Low}-{Thrust}
  {Transfer},} {\em Journal of Guidance, Control, and Dynamics\/}, Vol.~41,
  No.~3, March 2018, pp.~658--672. \\
  https://doi.org/10.2514/1.G002922

\bibitem{pritchett_impulsive_2018}
Pritchett, R.~E., Zimovan, E., and Howell, K., \enquote{Impulsive and
  {Low}-{Thrust} {Transfer} {Design} {Between} {Stable} and {Nearly}-{Stable}
  {Periodic} {Orbits} in the {Restricted} {Problem},} {\em 2018 {Space}
  {Flight} {Mechanics} {Meeting}\/}, American Institute of Aeronautics and
  Astronautics, Kissimmee, Florida, Jan. 2018. \\
  https://doi.org/10.2514/6.2018-1690


\bibitem{senent_low-thrust_2005}
Senent, J., Ocampo, C., and Capella, A., \enquote{Low-{Thrust}
  {Variable}-{Specific}-{Impulse} {Transfers} and {Guidance} to {Unstable}
  {Periodic} {Orbits},} {\em Journal of Guidance, Control, and Dynamics\/},
  Vol.~28, No.~2, March 2005, pp.~280--290. \\
  https://doi.org/10.2514/1.6398

\bibitem{caillau_minimum_2012}
Caillau, J.-B. and Daoud, B., \enquote{Minimum {Time} {Control} of the
  {Restricted} {Three}-{Body} {Problem},} {\em SIAM Journal on Control and
  Optimization\/}, Vol.~50, No.~6, Jan. 2012, pp.~3178--3202. \\
  https://doi.org/10.1137/110847299


\bibitem{bonnard_time-minimum_2016}
Chyba, M., Patterson, G., and Picot, G., \enquote{Time-{Minimum} {Control} of
  the {Restricted} {Elliptic} {Three}-{Body} {Problem} {Applied} to {Space}
  {Transfer},} {\em Recent {Advances} in {Celestial} and {Space}
  {Mechanics}\/}, edited by B.~Bonnard and M.~Chyba, Vol.~23, Springer
  International Publishing, Cham, 2016, pp. 179--208, Series Title: Mathematics
  for Industry. \\
  https://doi.org/10.1007/978-3-319-27464-5


\bibitem{yoon_minimum-fuel_2023}
Yoon, S.~W. and Petukhov, V., \enquote{Minimum-Fuel Low-Thrust Trajectories to
  the {Moon},} {\em Acta Astronautica\/}, Vol.~210, Sept. 2023, pp.~102--116. \\
  https://doi.org/10.1016/j.actaastro.2023.05.006

\bibitem{singh_exploiting_2021}
Singh, S.~K., Anderson, B.~D., Taheri, E., and Junkins, J.~L.,
  \enquote{Exploiting Manifolds of {L}1 Halo Orbits for End-to-End
  {Earth}–{Moon} Low-Thrust Trajectory Design,} {\em Acta Astronautica\/},
  Vol.~183, June 2021, pp.~255--272. \\
  https://doi.org/10.1016/j.actaastro.2021.03.017

\bibitem{beolchi_multiple-arc_2024}
 Beolchi, A., Pontani, M., Pozzi, C., and Fantino, E., \enquote{Multiple-Arc
   Optimization of Low-Thrust {Earth}-{Moon} Orbit Transfers Leveraging
   Implicit Costate Transformation,} {\em Acta Astronautica\/}, Vol.~220, July
   2024, pp.~330--344. \\
   https://doi.org/10.1016/j.actaastro.2024.04.030

\bibitem{chupin_low-thrust_2017}
Chupin, M., Haberkorn, T., and Tr\'{e}lat, E., \enquote{Low-Thrust {Lyapunov} to
  {Lyapunov} and {Halo} to {Halo} Missions with {L}$^{\textrm{2}}$
  -Minimization,} {\em ESAIM: Mathematical Modelling and Numerical Analysis\/},
  Vol.~51, No.~3, May 2017, pp.~965--996. \\
  https://doi.org/10.1051/m2an/2016044

\bibitem{du_transfer_2022}
Du, C., Starinova, O.~L., and Liu, Y., \enquote{Transfer between the Planar
  {Lyapunov} Orbits Around the {Earth}–{Moon} {L2} Point Using Low-Thrust
  Engine,} {\em Acta Astronautica\/}, Vol.~201, Dec. 2022, pp.~513--525. \\
  https://doi.org/10.1016/j.actaastro.2022.09.056

\bibitem{du_low-thrust_2023}
Du, C., Starinova, O., and Liu, Y., \enquote{Low-thrust Transfer Trajectory
  Planning and Tracking in the {Earth}–{Moon} Elliptic Restricted Three-Body
  Problem,} {\em Nonlinear Dynamics\/}, Vol.~111, No.~11, June 2023,
  pp.~10201--10216. \\
  https://doi.org/10.1007/s11071-023-08383-0

\bibitem{kelly_orthogonal_2023}
Kelly, P., Junkins, J.~L., and Majji, M., \enquote{Orthogonal {Approximation}
  of {Invariant} {Manifolds} in the {Circular} {Restricted} {Three}-{Body}
  {Problem},} {\em Journal of Guidance, Control, and Dynamics\/}, April 2023,
  pp.~1--12. \\
  https://doi.org/10.2514/1.g007304

\bibitem{mingotti_optimal_2011}
Mingotti, G., Topputo, F., and Bernelli-Zazzera, F., \enquote{Optimal Low-Thrust Invariant Manifold Trajectories via Attainable Sets,}, {\em{Journal} {of} {Guidance}, {Control}, {and} {Dynamics}}, Vol. 34, No. 6, pp. 1644–1656, November 2011. \\
https://doi.org/10.2514/1.52493


\bibitem{ozimek_low-thrust_2010}
Ozimek, M.~T. and Howell, K.~C., \enquote{Low-{Thrust} {Transfers} in the
  {Earth}-{Moon} {System}, {Including} {Applications} to {Libration} {Point}
  {Orbits},} {\em Journal of Guidance, Control, and Dynamics\/}, Vol.~33,
  No.~2, March 2010, pp.~533--549. \\
  https://doi.org/10.2514/1.43179

\bibitem{perez-palau_fuel_2018}
Pérez-Palau, D. and Epenoy, R., \enquote{Fuel Optimization for Low-thrust
  {Earth}–{Moon} Transfer via Indirect Optimal Control,} {\em Celestial
  Mechanics and Dynamical Astronomy\/}, Vol.~130, No.~2, Feb. 2018, pp.~21. \\
  https://doi.org/10.1007/s10569-017-9808-2

\bibitem{conway_survey_2012}
Conway, B.~A., \enquote{A {Survey} of {Methods} {Available} for the {Numerical}
  {Optimization} of {Continuous} {Dynamic} {Systems},} {\em Journal of
  Optimization Theory and Applications\/}, Vol.~152, No.~2, Feb. 2012,
  pp.~271--306. \\
  https://doi.org/10.1007/s10957-011-9918-z

\bibitem{caillau_algorithmic_2023}
Caillau, J.-B., Ferretti, R., Tr\'{e}lat, E., and Zidani, H., \enquote{An
  Algorithmic Guide for Finite-Dimensional Optimal Control Problems,} {\em
  Handbook of {Numerical} {Analysis}\/}, Vol.~24, Elsevier, 2023, pp. 559--626. \\
  https://doi.org/10.1016/bs.hna.2022.11.006

\bibitem{mingotti_combined_2007}
Mingotti, G., Topputo, F., and Bernelli-Zazzera, F., \enquote{Combined
  {Optimal} {Low}-{Thrust} and {Stable}-{Manifold} {Trajectories} to the
  {Earth}-{Moon} {Halo} {Orbits},} {\em {AIP} {Conference} {Proceedings}\/},
  Vol. 886, AIP, Princeton, New Jersey (USA), 2007, pp. 100--112. \\
  https://doi.org/10.1063/1.2710047

\bibitem{parrish_low-thrust_2016}
Parrish, N.~L., Parker, J.~S., Hughes, S.~P., and Heiligers, J.,
  \enquote{Low-{Thrust} {Transfers} from {Distant} {Retrograde} {Orbits} to
  {L2} {Halo} {Orbits} in the {Earth}-{Moon} {System},} Darmstadt, Germany,
  March 2016.

\bibitem{pritchett_low-thrust_2017}
Pritchett, R., Howell, K., and Grebow, D., \enquote{Low-Thrust Transfer Design
  Based on Collocation Techniques: Applications in the Restricted Three-Body
  Problem,} AAS 17-62, \textit{AIAA/AAS Astrodynamics Specialist Conference}, Columbia River Gorge, Stevenson, WA, August 21--24, 2017.

\bibitem{ams-1913}
Macmillan, W. D., ``On Poincar\'{e}'s Correction to Bruns' Theorem,'' \textit{Bulletin of the American Mathematical Society}, 19, 1913, 349--355. %\\

\bibitem{moser-1971}
Siegel, C. L. and Moser, J. K., \textit{Lectures on Celestial Mechanics}, Springer-Verlag, New York, N.Y., 1971 (Translated by C. I. Kalme). %\\


\bibitem{ferrari_periodic_2018}
Ferrari, F. and Lavagna, M., \enquote{Periodic Motion Around Libration Points
  in the {Elliptic} {Restricted} {Three}-{Body} {Problem},} {\em Nonlinear
  Dynamics\/}, Vol.~93, No.~2, July 2018, pp.~453--462. \\
  https://doi.org/10.1007/s11071-018-4203-4

\bibitem{caillau_minimum_2012-1}
Caillau, J.-B., Daoud, B., and Gergaud, J., \enquote{Minimum Fuel Control of
  the Planar Circular Restricted Three-Body Problem,} {\em Celestial Mechanics
  and Dynamical Astronomy\/}, Vol.~114, No. 1-2, Oct. 2012, pp.~137--150. \\
  https://doi.org/10.1007/s10569-012-9443-x

\bibitem{oshima_global_2017}
Oshima, K., Campagnola, S., and Yanao, T., \enquote{Global Search for
  Low-Thrust Transfers to the {Moon} in the Planar Circular Restricted
  Three-Body Problem,} {\em Celestial Mechanics and Dynamical Astronomy\/},
  Vol.~128, No. 2-3, June 2017, pp.~303--322. \\
  https://doi.org/10.1007/s10569-016-9748-2

 \bibitem{campagnola_subregions_2008}
Campagnola, S., Lo, M., and Newton, P., \enquote{Subregions of Motion and
  Elliptic Halo Orbits in the Elliptic Restricted Three-Body Problem,} AAS 08-200, \textit{AIAA/AAS Spaceflight Mechanics Meeting}, Galveston, TX, January 27--31, 2008.

\bibitem{celletti_dynamics_2024}
Celletti, A., Lhotka, C., and Pucacco, G., \enquote{The Dynamics Around the
  Collinear Points of the Elliptic Three-Body Problem: {A} Normal Form
  Approach,} {\em Physica D: Nonlinear Phenomena\/}, Vol.~468, Nov. 2024,
  pp.~134302. \\
  https://doi.org/10.1016/j.physd.2024.134302

\bibitem{szebehely_theory_1967}
Szebehely, V.~G., {\em Theory of Orbits: {The} Restricted Problem of Three
  Bodies\/}, Academic Press, New York, 1967.

\bibitem{broucke_periodic_1969}
Broucke, R.~A., \enquote{Periodic Orbits in the Elliptic Restricted Three-Body
  Problem,} Tech. rep. 32-1168, Jet Propulsion Laboratory, California Institute of
  Technology, Pasadena, CA, Feb, 1968.

\bibitem{broucke_stability_1969}
Broucke, R., \enquote{Stability of Periodic Orbits in the Elliptic, Restricted
  Three-Body Problem.} {\em AIAA Journal\/}, Vol.~7, No.~6, June 1969,
  pp.~1003--1009. \\
  https://doi.org/10.2514/3.5267

\bibitem{hyeraci_method_2010}
Hyeraci, N. and Topputo, F., \enquote{Method to {Design} {Ballistic} {Capture}
  in the {Elliptic} {Restricted} {Three}-{Body} {Problem},} {\em Journal of
  Guidance, Control, and Dynamics\/}, Vol.~33, No.~6, Nov. 2010,
  pp.~1814--1823. \\
  https://doi.org/10.2514/1.49263

\bibitem{dixon-ross-TN_2025}
Dixon, M. J. and Ross, I. M., ``Nuances in Propellant Computation in the Elliptic Restricted Three-Body Problem,'' \textit{Journal of Guidance, Control and Dynamics}, Vol. 49, No. 1, Jan. 2026. \\
https://doi.org/10.2514/1.G009053

\bibitem{hiday-johnston_transfers_1994}
Hiday-Johnston, L.~A. and Howell, K.~C., \enquote{Transfers Between
  Libration-Point Orbits in the Elliptic Restricted Problem,} {\em Celestial
  Mechanics \& Dynamical Astronomy\/}, Vol.~58, No.~4, April 1994,
  pp.~317--337. \\
  https://doi.org/10.1007/bf00692008

\bibitem{neelakantan_two-impulse_2022}
Neelakantan, R. and Ramanan, R.~V., \enquote{Two-Impulse Transfer to
  Multi-Revolution Halo Orbits in the {Earth}-{Moon} Elliptic Restricted
  Three-Body Problem Framework,} {\em Journal of Astrophysics and Astronomy\/},
  Vol.~43, No.~2, Aug. 2022, pp.~50. \\
  https://doi.org/10.1007/s12036-022-09830-x

\bibitem{scantamburlo_interplanetary_2022}
Scantamburlo, E., Guzzo, M., and Paez, R.~I., \enquote{Interplanetary Transfers
  Using Stable and Unstable Manifold Tubes Originating at {L}1 and {L}2 in
  the Elliptic Restricted Three-Body Problems,} {\em Acta Astronautica\/},
  Vol.~200, Nov. 2022, pp.~97--110. \\
  https://doi.org/10.1016/j.actaastro.2022.07.036

\bibitem{ross_how_2004}
Ross, I. M., \enquote{How to {Find} {Minimum}-{Fuel} {Controllers},} {\em {AIAA}
  {Guidance}, {Navigation}, and {Control} {Conference} and {Exhibit}\/},
  American Institute of Aeronautics and Astronautics, Providence, Rhode Island,
  Aug. 2004. \\
  https://doi.org/10.2514/6.2004-5346

\bibitem{ross_space_2006}
Ross, I.~M., \enquote{Space {Trajectory} {Optimization} and {$L^1$}-{Optimal}
  {Control} {Problems},} {\em Elsevier {Astrodynamics} {Series}\/}, Vol.~1,
  Elsevier, 2006, pp. 155--VIII. \\
  https://doi.org/10.1016/s1874-9305(07)80008-2

\bibitem{ross_universal_2023}
Ross, I. M., Proulx, R. J., and Borges, C. F., ``A Universal Birkhoff Pseudospectral Method for Solving
Boundary Value Problems,'' \textit{Applied Mathematics and Computation}, Vol.~454, 128101, October, 2023. \\ https://doi.org/10.1016/j.amc.2023.128101

\bibitem{ross_universal_2023-1}
Ross, I. M., ``A Universal Birkhoff Theory for Fast Trajectory Optimization,'' \textit{Journal of Guidance, Control and Dynamics},  vol. 47, no. 12, pp.~2468--2481, Dec. 2024. https://doi.org/10.2514/1.G007737

\bibitem{proulx_implementations_2023}
Proulx R. J. and Ross, I. M., ``Implementations of the Universal Birkhoff Theory for Fast Trajectory Optimization,'' \textit{Journal of Guidance, Control and Dynamics},  vol. 47, no. 12, pp. 2482--2496, Dec. 2024. https://doi.org/10.2514/1.G007738

\bibitem{ross_hessianTN_2025}
Ross, I. M., ``Hessians in Birkhoff-Theoretic Trajectory Optimizaiton,'' \textit{Journal of Guidance, Control and Dynamics}, Vol.~48, No.~9, 2025, pp.~2105--2112. https://doi.org/10.2514/1.G008778.

\bibitem{sandia_millionPt_2025}
Javeed, A., Kouri, D. P., Ridzal, D., Steinman, J. D., and Ross, I. M., ``A Million-Point Fast Trajectory Optimization Solver,'' AAS Paper 25-689, August 2025. arXiv preprint arXiv:2509.01855.\\ https://arxiv.org/abs/2509.01855

\bibitem{singular-arc-2010-CDC}
Park, C., Ross, I. M., Gong, Q., and Yan, H., ``Necessary Conditions for Singular Arcs for General Restricted Multi-Body Problem,'' \textit{49th IEEE Conference on Decision and Control}, IEEE, Piscataway, NJ, 2010, pp. 5816-5821. \\
https://doi.org/10.1109/cdc.2010.5717305

\bibitem{ross_enhancements_2020}
Ross, I.~M., \enquote{Enhancements to the {DIDO} {Optimal} {Control}
  {Toolbox},} June 2020, arXiv:2004.13112.\\
  https://arxiv.org/abs/2004.13112

\bibitem{gong_spectral_2008}
Gong, Q., Fahroo, F., and Ross, I.~M., \enquote{Spectral {Algorithm} for
  {Pseudospectral} {Methods} in {Optimal} {Control},} {\em Journal of Guidance,
  Control, and Dynamics\/}, Vol.~31, No.~3, May 2008, pp.~460--471. \\
  https://doi.org/10.2514/1.32908

\bibitem{scaling}
Ross, I. M., Gong, Q., Karpenko, M. and Proulx, R. J., ``Scaling and Balancing for High-Performance Computation of Optimal Controls,''  \textit{Journal of Guidance, Control and Dynamics}, Vol.~41, No.~10, 2018, pp.~2086--2097.\\
https://doi.org/10.2514/1.g003382

\bibitem{ross:guess-free}
Ross, I. M. and Gong, Q., ``Guess-Free Trajectory Optimization,'' \textit{AIAA/AAS Astrodynamics Specialist Conference and Exhibit}, 18--21 August 2008, Honolulu, Hawaii, AIAA 2008-6273. \\
https://doi.org/10.2514/6.2008-6273
%https://core.ac.uk/download/pdf/36722475.pdf

\bibitem{LRO-CSM}
M. Karpenko, J. K. Halverson, Y. Tedla,
T. Lippman, K. Han, I. M. Ross, M. K. Barker,
M. S. Robinson, W. Carter, A. K. Boyd,
G. W. Patterson, C. M. O'Shea, F. S. Turner,
R. T. Poffenbarger, C. J. Dennehy,
T. Vanzwieten, and A. A. Wolf, ``Fast Attitude Maneuvers for NASA's Lunar Reconnaissance Orbiter: Practical Flight Application of Attitude Guidance Using Birkhoff Pseudospectral Theory and Hamiltonian Programming,''
\textit{IEEE Control Systems Magazine}, Vol.~44, No.~2, 2024, pp.~26--54. \\
https://doi.org/10.1109/mcs.2024.3358593

\bibitem{ross_pseudospectral_2004}
Ross, I. M. and Fahroo, F., ``Pseudospectral Knotting Methods for
Solving Optimal Control Problems,'' {\it Journal of Guidance,
Control and Dynamics}, Vol. 27, No. 3, pp. 397-405, 2004. \\
https://doi.org/10.2514/1.3426

\bibitem{auto-knots}
Gong, Q. and Ross, I. M., ``Autonomous Pseudospectral Knotting Methods for Space Mission Optimization,'' \textit{Advances in the Astronatuical Sciences}, Vol.~124, 2006, AAS 06-151, pp.~779--794.


\bibitem{Farq-1968}
Farquhar, R. W., “The Control and Use of Libration-Point Satellites,” Stanford University, CA, 1968. %Accessed: Aug. 29, 2024. [Online]. Available: https://ntrs.nasa.gov/api/citations/19680021946/downloads/19680021946.pdf

\bibitem{NRHO-2017}
Zimovan, E. M., Howell, K. C. and Davis, D. C., ``Near Rectilinear Halo Orbits and
Their Application in Cis-Lunar Space'' \textit{3rd IAA Conference on Dynamics and Control of Space Systems}, May 2017. IAA-AAS-DyCoSS3-125.

\bibitem{trefethen_approximation_2020}
Trefethen, L.~N., {\em Approximation Theory and Approximation Practice\/},
  SIAM, Society for Industrial and Applied Mathematics, Philadelphia, extended
  edition ed., 2020. \\
  https://doi.org/10.1137/1.9781611975949

\bibitem{boyd_chebyshev_2000}
Boyd, J.~P., {\em Chebyshev and {Fourier} {Spectral} {Methods}\/}.
Mineola, New York: Dover Publications, 2000.

\bibitem{fahroo_advances_2008}
Fahroo, F. and Ross, I.~M., \enquote{Advances in {Pseudospectral} {Methods} for
  {Optimal} {Control},} {\em {AIAA} {Guidance}, {Navigation} and {Control}
  {Conference} and {Exhibit}\/}, American Institute of Aeronautics and
  Astronautics, Honolulu, Hawaii, Aug. 2008. \\
  https://doi.org/10.2514/6.2008-7309

\bibitem{ross_review_2012}
Ross, I.~M. and Karpenko, M., \enquote{A Review of Pseudospectral Optimal
  Control: {From} Theory to Flight,} {\em Annual Reviews in Control\/},
  Vol.~36, No.~2, Dec. 2012, pp.~182--197. \\
  https://doi.org/10.1016/j.arcontrol.2012.09.002

\bibitem{fahroo:cheb-jgcd}
Fahroo, F. and Ross, I. M., ``Direct Trajectory Optimization by a
Chebyshev Pseudospectral Method,'' \textit{Journal of Guidance,
Control and Dynamics}, Vol. 25, No. 1, 2002, pp.160-166. \\
https://doi.org/10.2514/2.4862

\bibitem{cheb-costate}
Gong, Q., Ross, I. M. and Fahroo, F.,  ``Costate Computation by a Chebyshev Pseudospectral Method,''  \textit{Journal of Guidance, Control and Dynamics}, Vol. 33, No. 2, pp. 623-628, 2010. \\
https://doi.org/10.2514/1.45154

\bibitem{t_a_driscoll_chebfun_2014}
{T. A. Driscoll}, Hale, N., and Trefethen, L.~N., editors, {\em Chebfun
  {Guide}\/}, Pafnuty Publications, Oxford, 2014.

\bibitem{trefethen-myths-2011}
L. N. Trefethen, ``Six Myths of Polynomial Interpolation and Quadrature,'' \textit{Mathematics Today}, 47
(2011), pp. 184--188.

\bibitem{koeppen_fast_2019}
Koeppen, N., Ross, I.~M., Wilcox, L.~C., and Proulx, R.~J., \enquote{Fast
  {Mesh} {Refinement} in {Pseudospectral} {Optimal} {Control},} {\em Journal of
  Guidance, Control, and Dynamics\/}, Vol.~42, No.~4, April 2019, pp.~711--722. \\
  https://doi.org/10.2514/1.g003904


\bibitem{rossJCAM-1}
Ross, I. M., ``An Optimal Control Theory for Nonlinear Optimization,'' \textit{Journal of Computational and Applied Mathematics}, Vol.~354, July 2019, 39--51. \\
http://dx.doi.org/10.1016/j.cam.2018.12.044

\bibitem{ross-CD}
Ross, I. M., ``Derivation of Coordinate Descent Algorithms from Optimal Control Theory'' \textit{Operartions Research Forum}, Vol.~4, No.~2, 2023. \\
https://doi.org/10.1007/s43069-023-00215-6

\bibitem{ross-JNVA}
Ross, I. M., ``On the Mathematics of the Natural Physics of Optimization,'' \emph{Journal of Nonlinear and Variational Analysis}, to appear. 

\end{thebibliography}
\end{document}